\magnification=\magstep 1
\baselineskip=15pt     %normal baselineskip = 12 pt
%\baselineskip=18pt     %normal baselineskip = 12 pt
%\hsize=6.5truein
%\hoffset=.25truein
\parskip=3pt plus1pt minus.5pt
\overfullrule=0pt
\font\hd=cmbx10 scaled\magstep1

\def\Pone{{\bf P}^1}
\def\P{{\bf P}}
\def\mod{\mathop{\rm mod}}

\def\C{{\cal C}}

\def\O{{\cal O}}

\def\T{{\cal T}}
\def\H{{\cal H}}

\def\F{{\cal F}}

\def\M{{\cal M}}

\def\R{{\cal R}}

\def\Pic{{\rm Pic}}

\def\Im{\mathop{\rm Im}}

\def\re{\mathop{\rm Re}}
\def\im{\mathop{\rm Im}}

\def\boldz{{\bf Z}}

\def\exact#1#2#3{0\rightarrow#1\rightarrow#2\rightarrow#3\rightarrow0}

\def\mapright#1{\smash{
  \mathop{\longrightarrow}\limits^{#1}}}

\def\u{{\frak u}}

\input amssym.tex

%\input prepictex.tex
%\input pictex.tex
%\font\thinlinefont=cmr5

\def\C{{\bf C}}
\def\P{{\bf P}}
\def\R{{\bf R}}
\def\Z{{\bf Z}}
\def\O{{\cal O}}

\def\T{{\cal T}}
\def\Re{{\rm Re}}
\def\Im{{\rm Im}}
\def\px{\partial_x}
\def\pxc{\partial_{\bar x}}
\def\py{\partial_y}
\def\pyc{\partial_{\bar y}}
\def\pv{\partial_{\rm v}}
\def\pvc{\bar\partial_{\rm v}}
\def\ph{\partial_{\rm h}}
\def\phc{\bar\partial_{\rm h}}
\def\pc{\bar\partial}
\def\th{\vartheta_{\rm h}}
\def\tv{\vartheta_{\rm v}}
\def\thc{\bar\vartheta_{\rm h}}
\def\tvc{\bar\vartheta_{\rm v}}

\centerline{\hd Large Complex Structure Limits of K3 Surfaces.}
\bigskip
\centerline{\it Mark Gross\footnote{*}{Supported in part by NSF grant
DMS-9700761, Trinity College, Cambridge, and EPSRC}\hfil P.M.H. Wilson}
\medskip
\centerline{August 1st, 2000; revised September 7th, 2000 and April 2nd, 2001}
\medskip
{\settabs 3 \columns
\+Mathematics Institute&&Department of Pure Mathematics\cr
\+University of Warwick&&University of Cambridge\cr
\+Coventry, CV4 7AL, U.K.&&Cambridge, CB3 0WB, U.K.\cr
\+mgross@maths.warwick.ac.uk&&pmhw@dpmms.cam.ac.uk\cr}
\bigskip

{\hd \S 0. Introduction.}

The notion of large complex structure limit plays a special role in
the theory of mirror symmetry. If $X$ is a Calabi--Yau manifold,
a large complex structure limit point is a point in a compactified moduli
space of complex structures $\overline{\M}_X$ 
on $X$ which, in some sense, represents the ``worst
possible degeneration'' of the complex structure. This notion was given
a precise Hodge-theoretic meaning in [27]. The basic example to keep in mind
of this sort of degeneration is the degeneration of a hypersurface of degree
$n+1$ in $\P^n$ to a union of the $n+1$ coordinate hyperplanes. Mirror
symmetry posits the existence of a mirror to $X$ associated to each
large complex structure limit point of $X$. To first approximation, this means
that if $p\in\overline{\M}_X$ is a large complex structure limit point
in a compactification of the complex moduli space of $X$, then there
exists a mirror $\check X$ and an isomorphism between a neighbourhood
of $p$ in $\overline{\M}_X$ and the complexified K\"ahler moduli space of
$\check X$ which preserves certain additional information, 
such as the Yukawa couplings (which will not concern us in this paper).
This isomorphism is known as the mirror map.

Now the Strominger-Yau-Zaslow conjecture [32] suggests that mirror
symmetry can be explained by the existence of a special Lagrangian fibration
on $X$ when the complex structure on $X$ is near a large complex structure
limit point. The mirror $\check X$ is then expected to be constructed
as the dual of this special Lagrangian fibration. The notion of special
Lagrangian is a metric one: it depends on both the complex structure
(determined by a holomorphic $n$-form $\Omega$ on $X$, where 
$n=\dim_{\bf C} X$), and a Ricci-flat K\"ahler metric, determined by
its K\"ahler form $\omega$. Thus we expect the existence of special
Lagrangian fibrations will depend a great deal on the metric properties
of Calabi--Yau manifolds near large complex structure limit points. 

The simplest example of such a situation occurs for elliptic curves.
Consider the family of elliptic curves $E_{\alpha}={\bf C}/\langle 1, 
i\alpha\rangle$, with $\alpha\rightarrow\infty$. We also choose a 
Ricci-flat, i.e. flat, metric $g$, which we will take to be the standard
Euclidean metric. As $\alpha\rightarrow\infty$, the complex structure
approaches the large complex structure limit point in the moduli space
of elliptic curves; the period $i\alpha$ is approaching the cusp
point of the compactification of $\H/SL_2(\boldz)$.

Now given the metric $g$, as $\alpha\rightarrow\infty$ it is clear that these
elliptic curves converge to an infinitely long cylinder. However, if
we rescale the metric, with $g_{\alpha}=g/\alpha$, then
$Vol(E_{\alpha})=1$ in this metric. With this metric, we can instead
view $E_{\alpha}$ as ${\bf C}/\langle 1/\sqrt{\alpha}, i\sqrt{\alpha}\rangle$
with the standard Euclidean metric, and then $E_{\alpha}$ converges
to a line as $\alpha\rightarrow\infty$.

Finally, we may renormalize the metric again so that the diameter
of $E_{\alpha}$ remains bounded, with the metric $g_{\alpha}=g/\alpha^2$.
Then $E_{\alpha}$ can be identified with ${\bf C}/\langle 1/\alpha,i\rangle$
with the Euclidean metric, and $E_{\alpha}$ converges to a circle.

Of course, in this situation the special Lagrangian $T^1$-fibration on 
$E_{\alpha}$ is $E_{\alpha}\rightarrow S^1$ obtained by projection onto
the imaginary axis. So with the second and third choices of normalization,
the special Lagrangian fibres collapse.

This is a rather trivial example, but forms a good basis for speculating
about what might happen in higher dimensions. Intuitively, if we normalize
the metric so as to keep the volume of the manifold bounded, we expect to
see the fibres of the hypothetical special Lagrangian fibration contracting
down to points; if furthermore we normalize so as to have bounded diameter, 
we expect the Calabi--Yau manifold to ``converge'' (in a sense we will make more 
explicit in \S 6) to a sphere of dimension $n$. 

To test this picture, and to improve our understanding of Ricci-flat
metrics, we have chosen to study the metric on K3 surfaces approaching
large complex structure limit points. This is made easier by the fact
that special Lagrangian fibrations are known to exist on K3 surfaces
by a standard trick of performing a hyperk\"ahler rotation of the complex
structure, so that one reduces the problem of finding a special Lagrangian
fibration to that of finding an elliptic fibration. Using this, we show
in \S 1 that after performing this hyperk\"ahler rotation,
approaching a large complex structure limit is more or less the same
as fixing the complex structure on
a K3 elliptic fibration $f:X\rightarrow \Pone$, and letting the K\"ahler
form $\omega$ on $X$ vary in such a way that the area of the fibres 
approaches zero. Thus we ask the question: what does a Ricci-flat metric
on an elliptic K3 surface look like when the area of the fibres is
very small? 

This is an interesting question even if one is
not interested in mirror symmetry. 
In [1], M. Anderson studied
degenerations of Ricci-flat metrics on K3 surfaces. If the volume of the
surface is fixed and the diameter remains bounded, then the metrics
converge to an orbifold metric (corresponding to degeneration to a K3
with rational double points). This picture of the moduli space of
K3 surfaces with orbifold metric was originally studied in [21].
If the diameter is unbounded, Anderson
proved collapsing must occur, but gave no more detailed information.
The case under consideration in this paper can be considered to be
the most extreme degeneration of metric. In particular, the orbifold case
and the elliptic fibration case are the only K\"ahler degenerations,
in which the complex structure of the K3 surface is held fixed.
We will in fact
consider a slightly more general situation, where the complex structure
still varies to some extent. This is described more clearly
in \S 1.

We assume the generic case, so that $f$ has 24 singular fibres, each
of Kodaira type $I_1$ (a pinched torus). If $X_0$ denotes the complement
of these 24 singular fibres, then it is possible to write down a family
of explicit Ricci-flat metrics which we refer to as {\it semi-flat}:
these metrics are in fact flat when restricted to the fibres. The semi-flat
metric was first introduced in [12]. There, it was used to get a first
approximation to a complete Ricci-flat metric on the complement of a fibre of
a rational elliptic surface. In [12], an arbitrary metric was then glued
in to take care of the singular fibres so that techniques of [33,34]
could  be applied to obtain a complete Ricci-flat metric on this manifold.
While we follow this idea in spirit, we have here
a new ingredient we can take advantage of.
There is an explicit Ricci-flat metric 
defined in a neighbourhood of each singular fibre,
first written down by Ooguri and Vafa in [29].
It is not semi-flat, but it in fact decays to a
semi-flat metric exponentially. We can glue 24 copies of the Ooguri--Vafa
metric in to the semi-flat metric, and thereby obtain a metric which
is remarkably close to being Ricci-flat: in fact, the Ricci curvature
is bounded in absolute value by $O(e^{-C/\epsilon})$, where $\epsilon$
denotes the area of a fibre. Thus as $\epsilon\rightarrow 0$, the 
Ricci curvature of this glued metric approaches zero very rapidly.

We then follow standard techniques to show that the genuine Ricci-flat
metric representing the same K\"ahler class is very close to the glued metric,
hence showing the explicit metric we constructed is a very good approximation
to the genuine metric. We follow the proof of Kobayashi in [20], based on the
original methods of Yau [35] --- cf. also [7,33,34]. In [20] 
Kobayashi proves that near a Kummer surface, the Ricci-flat metric on
a K3 surface is close to the flat orbifold metric on the Kummer surface.
While the techniques are the same, it is perhaps surprising that they
apply in our circumstances. Indeed, if the volume of the K3 surface is
held fixed, then as $\epsilon\rightarrow 0$, the
diameter of our metric approaches $\infty$. Thus the relevant Sobolev
constant approaches zero, and so it will be important to control
this precisely. It turns out that everything works because the starting
glued metric is already extremely close to being Ricci-flat.

More explicitly, for K\"ahler classes $[\omega_{\epsilon}]$ on $X$, where
$\epsilon$ denotes the volume of a fibre of $f$, we construct
a representative K\"ahler metric $\omega_{\epsilon}$ with very
small Ricci curvature. Yau's proof [35] of the Calabi conjecture
yields a solution $u_{\epsilon}$ to the equations
$$\eqalign{ (\omega_{\epsilon}+i\partial\bar\partial u_{\epsilon})^2
&= e^{F_{\epsilon}}\omega_{\epsilon}^2\cr
\int_X u_{\epsilon}\omega_{\epsilon}^2&=0\cr}$$
with $F_{\epsilon}=\log\left({\Omega\wedge\bar\Omega/2\over \omega_{\epsilon}^2}
\right)$. The metric $\omega_{\epsilon}+i\partial\bar\partial u_{\epsilon}$
is the desired Ricci-flat metric.
We obtain a global $C^2$-estimate (Lemma 5.3), namely that for some
positive constant $C$,
$$ C^{-1} \omega _{\epsilon} \le \omega _{\epsilon} + i \partial
\bar\partial u_{\epsilon} \le C \omega _{\epsilon}.$$
Moreover, the main theorem of the paper (Theorem 5.6) states that
for any simply connected open set $U \subset B$ whose closure 
is disjoint from the discriminant locus of $f$, and for any
$k\ge 2$, $0<\alpha <1$, there exist positive constants $C_1 , C_2 ,
\epsilon _0$ such that, for all $\epsilon < \epsilon _0$,
$$ \| u_{\epsilon} \| _{C ^{k, \alpha}} \le C_1 e^{-C_2 / \epsilon },$$
where the $C^{k,\alpha}$ norm is on the set $f^{-1}(U)$. Thus, away
from the singular fibres, $\omega_{\epsilon}$ is a very good approximation to
the actual Ricci-flat metric. See Theorem 5.6 for a more precise
statement, which requires some care in the choice of the K\"ahler class
$[\omega_{\epsilon}]$.

The information obtained gives a clear picture of the metric behaviour
as $\epsilon\rightarrow 0$. Using the above results,
we prove the fibres are collapsing to points, and
that
away from the singular fibres, the metric
approaches the semi-flat metric.
In fact we will compute the Gromov--Hausdorff
limit of a sequence of K3 surfaces with $\epsilon\rightarrow 0$ and
the metrics renormalized so that the diameter remains bounded. This limit
is indeed an $S^2$, but the metric on the $S^2$ is singular at precisely
24 points corresponding to the singular fibres. See \S 6 for more precise
statements. There, we state a conjecture, also made independently
by Kontsevich, Soibelman, and Todorov, about
the Gromov--Hausdorff limit of Calabi--Yau manifolds approaching large
complex structure limit points. The above results prove
this conjecture in the two dimensional case.

The structure of the paper is as follows. In \S 1 we briefly review
mirror symmetry for K3 surfaces, so as to reduce the problem
to one of understanding elliptic fibrations. In \S\S 2 and 3, we introduce
various ways of thinking about Ricci-flat metrics on elliptic fibrations,
and then discuss required properties of the semi-flat and Ooguri--Vafa
metrics. In \S 4, we build the glued metric. In \S 5, we run through
the standard program to obtain estimates for Ricci-flat metrics, proving
the main result of the paper, Theorem 5.6. Finally, in \S 6, we 
relate these results to Gromov--Hausdorff convergence, and speculate as to
what kind of results in this direction might be expected and useful
in higher dimensions.
\medskip

{\it Acknowledgements:} We would like to thank G. Gibbons, N. Hitchin,
N.C. Leung, M. Micallef, M. Singer, Y. Soibelman, and R. Thomas.
The first author would
especially like to thank S.-T. Yau and E. Zaslow, with whom he held extensive
discussions about the Ooguri--Vafa metric in 1998.
\bigskip
  
{\hd \S 1. Identification of large complex structure limits.}

There are a number of variants of mirror symmetry for K3 surfaces: see
especially [10] for mirror symmetry between algebraic families of K3
surfaces and [4] for a more general version. We will use
an intermediate version here, following [14], \S 7, which
highlights the role of the special Lagrangian fibration. See also [17], \S 1.
We review
this point of view here. This will serve as motivation for Question 1.2
below, which will be
addressed in the remainder of the paper. However, the setup of mirror
symmetry will not be used again in this paper.

Let $L$ be the K3 lattice, $L=H^2(X,\boldz)$ for $X$ a K3 surface. Fix
a sublattice of $L$ isomorphic to the hyperbolic plane
$H$ generated by $E$ and $\sigma_0$,
with $E^2=0,\sigma_0^2=-2$, and $E.\sigma_0=1$. We will view mirror
symmetry as an involution
acting on the moduli space of triples $(X,{\bf B}+i\omega,\Omega)$
where $X$ is a marked K3 surface, $\Omega$ is the class of a holomorphic
2-form on $X$, $\omega\in E^{\perp}\otimes{\bf R}$ a
K\"ahler class on $X$, and
the $B$-field ${\bf B}$ lies in $E^{\perp}/E\otimes {\bf R}$.
In addition $\Omega$ is normalised
so that $\im\Omega\in E^{\perp}\otimes{\bf R}$ and 
$\omega^2=(\re\Omega)^2=(\im\Omega)^2$. 
Mirror symmetry interchanges $(X,{\bf B}+i\omega,\Omega)$ with
$(\check X,\check{\bf B}+i\check\omega,\check\Omega)$, where
$\check X$ denotes a marked K3 surface with the following data:
$$\eqalign{\check\Omega&\equiv
(E.\re\Omega)^{-1}(\sigma_0+{\bf B}+i\omega)\quad\mod E\cr
\check{\bf B}&\equiv (E.\re\Omega)^{-1}\re\Omega-\sigma_0 \quad\mod E\cr
\check\omega&\equiv (E.\re\Omega)^{-1}\im\Omega\quad\mod E.\cr}$$
The actual classes of $\check\Omega$ and $\check\omega$ are determined 
completely by the relation $(\re\check\Omega)^2=(\im\check\Omega)^2
=\check\omega^2$, $\check\omega.(\re\check\Omega)=\check\omega.(\im\check
\Omega)=(\re\check\Omega).(\im\check\Omega)=0$.

We can now identify the large complex structure limit of $\check X$.
This limit is mirror to the large K\"ahler limit of $X$. In the latter
limit, we keep the complex structure on $X$ fixed but allow the K\"ahler
form to go to infinity. More precisely, if $\{ {\bf B}_l+it_l\omega\}$
is a sequence of complexified K\"ahler forms on $X$ with $t_l>0$,
$t_l\rightarrow\infty$, then we say $\{{\bf B}_l+it_l\omega\}$ are
approaching the large K\"ahler limit in the complexified K\"ahler moduli
space of $X$. 

We will now take, for our purposes,

\proclaim Definition 1.1. For each $l$, let $\check X_l$ be the K3 surface
given by the data
$(\check X_l,\check {\bf B}_l+i\check\omega_l,\check\Omega_l)$
mirror to $(X,{\bf B}_l+it_l\omega,t_l\Omega)$. The sequence of surfaces 
$\{\check X_l\}$ is said to approach a large complex structure limit
point.

We will take this as the starting point of our analysis, and will not prove 
here that this is equivalent to other reasonable definitions of large
complex structure limits found in the literature (but see discussions in
[10]). 

The reader will note that we are cheating to some extent here, by only
approaching the large K\"ahler limit along a ray. The more general
approach might be to allow a more general sequence of K\"ahler forms.
However, this is more difficult to deal with because the elliptic
fibration which arises below will be varying. We will ignore this
difficulty in this paper, as it obscures our main objectives. 

Note that
$$\check\Omega_l=(t_l E.\re\Omega)^{-1}(\sigma_0
+{\bf B}_l+it_l\omega)\quad\mod E$$
and
$$\check\omega_l=(E.\re\Omega)^{-1}\im\Omega\quad\mod E.$$
More precisely, if a representative ${\bf B}_l$ for
${\bf B}_l\mod E$ is chosen in $E^{\perp}\otimes{\bf R}$ with
the property that ${\bf B}_l\cdot\sigma_0=0$, then the requirement
that $\check\Omega_l^2=0$ yields
$$\check\Omega_l=(t_lE.\re\Omega)^{-1}
(\sigma_0+({\bf B}_l+it_l\omega)+((t_l^2\omega^2-{\bf B}_l^2)/2+1
-it_l\omega.(\sigma_0+{\bf B}_l))E).$$
Furthermore
the requirement that $\check\omega_l.\check\Omega_l=0$ yields
$$\check\omega_l=(E.\re\Omega)^{-1}(\im\Omega
-(\im\Omega.(\sigma_0+{\bf B}_l))E).$$

The K\"ahler class $\check\omega_l$ is represented by a Ricci-flat
metric $\check g_l$, and we would like to understand the behaviour
of this metric as $t_l\rightarrow\infty$. It is convenient to perform
a hyperk\"ahler rotation, i.e. $\check g_l$ is also a K\"ahler metric on
the K3 surface $\check X_{l,K}$ with
$$\eqalign{\check\Omega_{l,K}&=\im\check\Omega_l+i\check\omega_l\cr
\check\omega_{l,K}&=\re\check\Omega_l.\cr}$$
This equality holds on the level of forms. Explicitly, in cohomology,
$$\eqalign{
\check\Omega_{l,K}&=(E.\re\Omega)^{-1}(\omega
+i\im\Omega-((\omega+i\im\Omega).(\sigma_0+{\bf B}_l))E)
\cr
\check\omega_{l,K}&=(t_lE.\re\Omega)^{-1}(\sigma_0
+{\bf B}_l)\quad\mod E.\cr}$$

We will assume that, for all $l$,
$E$ represents the class of a fibre of an elliptic
fibration $ f_l : \check X _{l,K} \to \P ^1$.  This elliptic 
fibration coincides with a special Lagrangian $T^2$-fibration on 
$\check X _l$.  For general choice of data, such elliptic fibrations
with fibre class $E$ automatically exist, since then $\Pic \check X _{l,K} =
\Z E$ and $E^2 = 0$.  For any choice of data, there always exists an 
elliptic fibration on $\check X _{l,K}$, but the class of the fibre 
might only be the image
of $E$ under reflections by $-2$ curves in $\Pic \check X_{l,K}$. (See
[17], \S 1 for further details.)

Note that the area of the fibre of $f_l$ under the metric $\check g_l$ is
$(t_l E.\re\Omega)^{-1}$, which goes to zero as $t_l\rightarrow\infty$.

Now $\check\Omega_{l,K}$ depends on $l$, but these classes only
differ by the pull-back of a class from $\Pone$. This in fact
tells us the elliptic
K3 surfaces $\check X_{l,K}$ are closely related. Indeed, 
if $f:X\rightarrow\Pone$
is an elliptic K3 surface, with holomorphic 2-form $\Omega$, then
whenever $\alpha$ is a 2-form on $\Pone$, $\Omega'=\Omega+f^*\alpha$
satisfies $\Omega'\wedge\Omega'=0$ as forms, and thus $\Omega'$ induces
another complex structure on $X$ such that $f$ remains a holomorphic
elliptic fibration in this new complex structure. All the surfaces
$\check X_{l,K}$ are clearly related in this way.
In particular, all these
elliptic surfaces have the same jacobian $\check J_K$, which is the
unique elliptic K3 surface with a holomorphic section
with complex structure induced by
$\check\Omega_{l,K}+f_l^*\alpha$ for some $\alpha$.

This now leads us to the following question:

\proclaim Question 1.2. Let $j:J\rightarrow\Pone$ be an
elliptic K3 surface with a section, and let $f_l:X_l\rightarrow\Pone$
be a sequence of elliptic K3 surfaces with jacobian $j:J\rightarrow
\Pone$. Let $\omega_l$ be a Ricci-flat K\"ahler metric on
$X_l$ with $Vol(X_l)$ independent of $l$. Let $\epsilon_l=
Area_{\omega_l}(f_l^{-1}(y))$ for any point $y\in\Pone$, and suppose 
$\epsilon_l\rightarrow 0$ as $l\rightarrow\infty$. Describe the
behaviour of the metric $\omega_l$ as $l\rightarrow\infty$. 

We will solve this question in this paper in the case
that the map $j$ has 24 Kodaira type $I_1$ fibres. This is true
for the generic K3 elliptic fibration.

We end this section with a few additional important comments about
this setup.

First, it is often convenient to identify the underlying
differentiable manifold of an elliptic K3 surface $f:X\rightarrow
B$ with that of its jacobian. This can be done in a reasonably
canonical fashion by choosing a $C^{\infty}$ section $\sigma_0:
B\rightarrow X$ of $f$. If $\Omega_X$ is a holomorphic 2-form
on $X$, then $\Omega_J=\Omega_X-f^*\sigma_0^*\Omega_X$
defines a new complex structure on $X$, in which $\sigma_0(B)$
is a holomorphic section. This new complex structure yields
the jacobian.

Another important point is that once a $C^{\infty}$ zero-section $\sigma_0$ for
$f:X\rightarrow B$ is chosen, we obtain a group structure
on the non-singular part of each fibre of $f$. Let $X^0\subseteq X$
be obtained by taking the union of
the identity components of each fibre. Then given a holomorphic
2-form $\Omega$ on $X$, we can construct a map from the holomorphic
cotangent bundle $\T_B^*$ to $X^0$, taking
the zero section of $\T_B^*$ to $\sigma_0(B)$, and with the property
that the pull-back of $\Omega$ to $\T_B^*$ is a form $\Omega_{can}+\alpha$,
where $\alpha$ is a 2-form pulled back from the base and $\Omega_{can}$
is the canonical holomorphic symplectic 2-form on $\T_B^*$. (See [14],
\S\S 2 and 7 for further details of this map.)  The canonical holomorphic
symplectic 2-form can be defined
in local coordinates. If
$y$ is a local holomorphic coordinate on the base $B$, we can take $x$
to be the corresponding canonical coordinate on the fibres of 
$\T^*_{B}$, so that the coordinate $(x_0,y_0)$ represents the
1-form $x_0 dy$ at the point in $B$ with coordinate
$y_0$. The pair $x,y$ are called {\it holomorphic canonical coordinates}.
Then the canonical 2-form on $\T^*_{B}$ is $dx\wedge dy$ in these
coordinates. 

The map $\T_B^*\rightarrow X^0$ also gives an exact sequence
$$\exact{R^1f_*\boldz}{\T^*_{\Pone}}{ X^0}.$$
$R^1f_*\boldz$ gives a degenerating family of lattices
in the fibres of the complex line
bundle $\T^*_B$. Thus working on the cotangent bundle of $B$
gives useful coordinates for $X$ away from the singular fibres,
and these coordinates will be used repeatedly in later sections.

\bigskip

{\hd  \S 2. Equations for Ricci-flatness.}

In this section we will discuss equations for Ricci-flatness in
different coordinate systems. We are interested in the behaviour
of the metric on an elliptic K3 fibration, and this metric behaves
in radically different ways away from the singular fibres as opposed
to a neighbourhood of the singular fibres. In these two different
cases, it will be useful to have two different coordinate systems
to study the metrics.

For studying the metric away from the singular fibres,
we adopt the set-up from the previous section, with $\pi : 
\T^*_B \to B$ where $B$ is an open subset of $\C$.
We are actually working on $X=\T^*_B/\Lambda$, where
$\Lambda$ is a holomorphically varying family of lattices in $\T^*_B$.
We will assume in this section that the zero section is holomorphic, so
that
the holomorphic 2-form on $X$ is 
induced by $\Omega=dx\wedge dy$
on $\T^*_B$,
where $y=y_1+iy_2$ and $x=x_1+ix_2$ are holomorphic canonical coordinates on
$\T^*_B$.
The K\"ahler form in these coordinates takes the form
$$\eqalign{
\omega=&{i\over 2}W(dx\wedge d\bar x+\bar b dx\wedge d\bar y
+b dy\wedge d\bar x
+(W^{-2}+|b|^2)dy\wedge d\bar y)\cr
=&{i\over 2}(W(dx+bdy)\wedge \overline{(dx+bdy)}+W^{-1}dy\wedge d\bar y).\cr}$$
Here $W$ and $b$ are defined by the above expression,
and the coefficient of $dy\wedge d\bar y$ 
is chosen to ensure the normalisation $\omega^2=(\im\Omega)^2$.
The function $W$ is real-valued and the function $b$ is complex-valued.
The K\"ahler condition is now $d\omega=0$.
This equation can be written as
$$\eqalign{
\py W&=\px(Wb)\cr
\py(W\bar b)&=
\px(W(W^{-2}+|b|^2)).\cr
}$$

Note that expanding the second equation out gives
$$W\py\bar b+\bar b\py W
=-W^{-2}\px W+(\px W)|b|^2
+W(b\px\bar b+\bar b\px b).$$
Using the first equation to replace $\py W$
and simplifying gives the above two equations being equivalent to
$$(\py-b\px)\bar b
=-W^{-3}\px W\leqno{(2.1)}$$
$$
(\py-b\px)W=W\px b.
\leqno{(2.2)}$$

Define the vector fields
$$\eqalign{\pv&=W^{-1}\px\cr
\ph&=\py-b\px\cr}$$
and denote by $\pvc$ and $\phc$ the complex
conjugate vector fields. The subscripts v and h denote the vertical
and horizontal vector fields respectively.
Let $\tv$ and $\th$ denote the dual
frame of one-forms, i.e.
$$\eqalign{\tv &=W(dx+bdy)\cr
\th &=dy.\cr}$$
Then 
$$\omega={i\over 2}W^{-1}(\tv\wedge\tvc
+\th\wedge\thc).$$
In addition, equations (2.1) and (2.2) take the simpler form
$$\ph \bar b=\pv W^{-1}\leqno{(2.1')}$$
$$-\ph W^{-1}=\pv b\leqno{(2.2')}$$

\noindent {\bf Remark 2.1.} While we don't use this here,
one can calculate that the holomorphic curvature
$\Theta=(\Theta_{ij})_{1\le i,j\le 2}$ of this metric is given by
$$\eqalign{\Theta_{11}=-\Theta_{22}&=\partial W\wedge\bar\partial W^{-1}
+W\partial\bar\partial W^{-1}+W^2\partial\bar b\wedge\bar\partial b
\cr
\Theta_{21}=-\bar\Theta_{12}
&=-W^{-1}\partial(W^2\bar\partial b).\cr}$$\bigskip

\noindent{\bf Example 2.2.  The standard semi-flat metric.}\hfill

We call a metric {\it semi-flat} if it restricts to a flat metric on
each elliptic fibre. As above, let
$B\subseteq {\bf C}$ an open
subset, $y$ the coordinate on ${\bf C}$. Let $\tau_1(y),\tau_2(y)$
be two holomorphic functions on $B$ such that $\tau_1(y)dy,\tau_2(y)dy$
generate a lattice $\Lambda(y)\subseteq \T_{B,y}^*$ for each $y\in B$, giving
us the holomorphically varying family of lattices $\Lambda\subseteq
\T_B^*=B\times{\bf C}$. Typically, we may allow $\tau_1$ and $\tau_2$
to be multi-valued.
Assuming without loss of generality that 
$\Im(\bar\tau_1\tau_2) > 0$, 
then a Ricci-flat metric on $X=(B\times {\bf C})/\Lambda$ is given by  the data
$$\eqalign{W&={\epsilon\over  \Im(\bar\tau_1\tau_2)} \cr
b&=-{W\over\epsilon}[\Im(\tau_2\bar x)\partial_y\tau_1+\Im(\bar\tau_1 x)
\partial_y\tau_2]\cr}$$
It is easy to check that these satisfy the equations
(2.1) and (2.2).
This metric, a priori defined on $\T_B^*$, descends to a metric on $X$,
and the area of a fibre of $f:X\rightarrow B$ is $\epsilon$.
We call this metric on $X$ {\it the standard semi-flat metric},
with K\"ahler form $\omega_{SF}$. 

The reader may check explicitly that this metric is independent of
the particular choice of generators for $\Lambda$, so that multi-valuedness
of $\tau_1$ and $\tau_2$ do not cause a problem. Furthermore, the
metric is independent of the choice of the coordinate $y$ (keeping
in mind that a change of the coordinate $y$ necessitates a change
of the canonical coordinate $x$, and hence the functions $\tau_1$, $\tau_2$).
This may also be seen as follows: The inclusion $R^1f_*\Z\cong \Lambda
\subseteq \T_B^*$ allows one to identify $(R^1f_*{\bf R})\otimes 
C^{\infty}(B)$  with the underlying $C^{\infty}$ vector
bundle $\T_B^*$, along with the Gauss-Manin connection $\nabla_{GM}$
on $\T_B^*$, the flat connection whose flat sections are sections of
$R^1f_*{\bf R}$. The standard semi-flat metric is the unique
semi-flat 
Ricci-flat K\"ahler metric satisfying the conditions
\item{(1)} The area of each fibre is $\epsilon$;
\item{(2)} $\omega_{SF}^2=(\re\Omega)^2=(\im\Omega)^2$;
\item{(3)} The orthogonal complement of each vertical tangent space 
is the horizontal tangent space of $\nabla_{GM}$ at that point.

This metric was described in [12], and in the more
general context of special Lagrangian fibrations in
[19], as well as [14], Example 6.4.

The reader should be aware however that if $T_{\sigma}:X\rightarrow X$
denotes translation by a holomorphic section $\sigma$, 
then $T_{\sigma}^*\omega_{SF}$
may give rise to a different semi-flat metric, satisfying conditions 
(1) and (2) but not (3). However, if $\sigma$ is not only
holomorphic but a flat section with respect to the Gauss-Manin connection
(so that $\sigma(y)=a_1\tau_1(y)+a_2\tau_2(y)$ for constants
$a_1$, $a_2$) then  $T_{\sigma}$ is an isometry and $T_{\sigma}^*\omega_{SF}
=\omega_{SF}$, $T_{\sigma}^*\Omega=\Omega$.

It will also be useful to have the K\"ahler potential for the metric. 
This is a function $\varphi$ such that $\omega={i\over 2}\partial\bar\partial
\varphi$. Let $\phi_1$ and $\phi_2$ be anti-derivatives of $\tau_1$ and
$\tau_2$ respectively. Then we can take
$$\varphi={\epsilon\over \im(\bar\tau_1\tau_2)}\left(-{\bar x^2\over 2}
{\tau_1\over\bar\tau_1}+|x|^2-{x^2\over 2}{\bar\tau_1\over\tau_1}\right)
+{i\over 2\epsilon}(\phi_1\bar\phi_2-\bar\phi_1\phi_2).$$
This is well-defined on subsets $\T_U^*\subseteq\T_B^*$ for $U$ simply
connected, but not on $\T_B^*/\Lambda$.
\bigskip

\noindent {\bf Construction 2.3.  The Gibbons--Hawking Ansatz.}

We now describe the system of coordinates which is most suited to studying the 
hyperk\"ahler metric in a neighbourhood of a singular fibre of the elliptic fibration.
This system of coordinates goes under the name of the Gibbons--Hawking Ansatz, and the 
description in terms of a connection form on an $S^1$-bundle 
explained below is essentially the same as that given in [2], which
in turn is based on earlier work of Gibbons and Hawking, Hitchin, and others.

Let $U\subseteq {\bf R}^3$
be an  open set with the Euclidean metric, 
with coordinates $u_1,u_2,u_3$.
Let $\pi:X\rightarrow U$ be a principal $S^1$ bundle, with $S^1$ action
$S^1\times X\rightarrow X$ written as $(e^{it},x)\mapsto e^{it}\cdot x$.
Let $\theta$ be a connection 1-form on $X$, i.e. a $\u (1)=i{\bf R}$-valued
1-form invariant under the $S^1$-action and such that $\theta(\partial/\partial
t)=i$. The curvature of the connection $\theta$ is $d\theta=\pi^*\alpha$
for a 2-form $\alpha$ on $U$, and $i\alpha/2\pi$ represents the first
Chern class of the bundle (see [8], 
Appendix). Suppose $V$ is a positive real
function on $U$ satisfying $*dV=\alpha/2\pi i$. Let
$$\eqalign{
\omega_1&=du_1\wedge \theta/2\pi i+Vdu_2\wedge du_3\cr
\omega_2&=du_2\wedge \theta/2\pi i+Vdu_3\wedge du_1\cr
\omega_3&=du_3\wedge \theta/2\pi i+Vdu_1\wedge du_2.\cr}$$
Then $\omega_1^2=\omega_2^2=\omega_3^2$ is nowhere zero, and
$\omega_i\wedge\omega_j=0$, for $i\not=j$. Furthermore, 
$*dV=\alpha/2\pi i$ implies $d\omega_i=0$ for all $i$, since for instance 
$$d \omega _1 = -d u_1 \wedge d\theta /{2\pi i} + dV \wedge du_2 \wedge du_3 =
-d u_1 \wedge *dV + dV \wedge du_2 \wedge du_3 =0.$$
Therefore $\omega_1,\omega_2,\omega_3$ define a hyperk\"ahler metric on $X$.
Note $V$ is harmonic, since $d \alpha = 0$ implies that $*d*d V = 0$.  

Let $\theta _0$ denote the real 1-form $\theta/{2\pi i}$, and observe that 
$$-\omega _1 - i \omega _2 =  (\theta _0 - iVdu_3) \wedge (du_1+idu_2).$$
By taking this to be the (holomorphic) 2-form $\Omega$ on $X$, this determines an integrable almost complex
structure on
$X$, where $du_1 + i du_2$ and $\theta _0 -iVdu_3$ span the holomorphic cotangent space 
inside the complexified cotangent space. It follows that the (integrable) almost 
complex structure $J$ on the cotangent space is given by
$$ J(du_1 ) = -du_2, \quad J ( du_3 ) = -V^{-1} \theta _0 .$$  Thus, if we consider the K\"ahler form 
$\omega = \omega _3$ as an alternating tensor, and use the relation that if $g$ is the Riemannian metric, 
then $g(\zeta ,\xi ) = \omega (\zeta , J \xi)$, we obtain an expression for the metric 
$$ds^2 = V d{\bf u}\cdot d{\bf u} + V^{-1} \theta _0 ^2.$$
Usually, we shall in fact start from a positive harmonic function $V$ on $U$ such that $-* dV$ represents the
Chern class of the bundle.  Then we can always find a connection 1-form $\theta$ with $d\theta/{2\pi i} = *dV$, 
such a $\theta$ being uniquely determined up to pull-backs of closed 1-forms from $U$, and hence we 
obtain hyperk\"ahler metrics as above.\bigskip

\noindent \bf Remark 2.4.\quad \rm 
We will need to calculate some information about the curvature
of this metric. We can work
locally,  and therefore take the orthonormal moving coframe given by $V^{1/2} du_1 ,V^{1/2} du_2,
V^{1/2} du_3$ and  
$V^{-1/2} \theta_0$.  We can moreover write the connection form locally as 
$$\theta_0 = {dt\over {2\pi}} + A_1 du_1 + A_2 du_2 +A_3 du_3,$$ where $\nabla V = \nabla \times {\bf
A}$. To calculate the curvature, we may then apply Cartan's  method.
We obtain
$$  \| R\| ^2 = 12 V^{-6} |\nabla V |^4 +  V^{-4} \Delta (|\nabla V |^2)
- 6 V^{-5} (\nabla V)\cdot (\nabla (|\nabla V |^2)). $$
Using the fact that $V$ is harmonic, we then
recover the compact formula given in equation (32) of [28] 
that $$ \|R\|^2 =  {1\over 2} V^{-1} \Delta \Delta (V^{-1}). $$
\bigskip

\noindent {\bf Example 2.5.}\quad  If we consider the natural map $\C ^2 \setminus (0,0) \to \P
^1 (\C ) = S^2$, and restrict to $S^3 \subset \C ^2$, 
we easily check that the image of $(z_1 ,z_2)\in S^3 $ is  
$$  (2\, \Re (z_1 \bar z_2 ), 2\, \Im (z_1 \bar z_2 ), |z_1|^2 - |z_2|^2).$$
This is the standard Poincar\'e map.  The formula also
defines a map $X= \C ^2 \setminus (0,0) \to  \R ^3 \setminus (0,0,0)$; we compose this map with 
complex conjugation on $z_2$ to obtain a map $p : X= \C ^2 \setminus (0,0) \to  
\R ^3 \setminus (0,0,0)$, given by 
$$ p (z_1 , z_2) = (2\, \Re (z_1  z_2 ), 2\, \Im (z_1  z_2 ), |z_1|^2 - |z_2|^2).$$
This map exhibits 
 $X$ as an 
$S^1$-bundle over $R ^3 \setminus (0,0,0)$, with Chern class $\pm 1$.  The action of $S^1$ on 
$X$ is given by $e^{it}\cdot (z_1,z_2)=(e^{it}z_1, e^{-it}z_2)$.  Note also that if we compose $p$
with projection onto the first two factors, we obtain the map sending $(z_1 ,z_2 )$ to 
$2z_1 z_2 $, holomorphic with respect to the standard complex structures.

We now choose a positive harmonic function $V$ on 
${\bf R}^3 \setminus (0,0,0)$ such that $$- \int _{S^2} *dV = \int _{S^2} i\alpha/{2\pi } =  \pm 1,$$ 
i.e. the Chern number is correct.  The particular examples of such $V$ we consider are
$$V = e + {1\over 4\pi |{\bf u}|} = e + {1\over 4\pi\sqrt{u_1^2+u_2^2+u_3^2}},$$
where $e\ge 0$.
The integral  $$ \int _{S^2} *d \left({1\over 4\pi\sqrt{u_1^2+u_2^2+u_3^2}}\right)$$ is easily seen
to be
$\pm 1$ (depending on the orientation of the sphere).

Now we take as connection form
$$\theta=i\Im (\bar z_1dz_1-\bar z_2 dz_2)/(|z_1|^2+|z_2|^2).$$  Then
$$d \theta/2\pi i = {-(u_1du_2\wedge du_3+u_2 du_3\wedge du_1+u_3du_1\wedge
du_2)\over 4\pi (u_1^2+u_2^2+u_3^2)^{3/2}} = * dV$$
as required.  We therefore obtain hyperk\"ahler metrics on $X$, which, for all $ e \ge 0$,
extend to metrics on $\C ^2$.  In fact, such metrics are  ALF (asymptotically locally flat), 
approaching a flat metric when
$|{\bf u}| \to \infty$, whilst being periodic in $t$.  When 
$e = 1$, the metric obtained is the Taub-NUT metric, and when $e = 0$, it is just a flat
metric on $\C ^2$.  To prove the assertions for $e=0$, straightforward calculations show that,
with $z_j=x_j+iy_j$,
$$\eqalign{
\omega_1&={1\over\pi}(dx_2\wedge dy_1-dx_1\wedge dy_2)\cr
\omega_2&={1\over\pi}(dx_1\wedge dx_2-dy_1\wedge dy_2)\cr
\omega_3&={1\over\pi}(dx_1\wedge dy_1+dx_2\wedge dy_2)\cr}$$
So $\omega_1,\omega_2,\omega_3$ extend to ${\bf C}^2$,
and yield a flat metric, as claimed.
\bigskip 
\eject

\noindent {\bf Construction 2.6.  
Gibbons--Hawking versus holomorphic coordinates}\hfill

In the Gibbons--Hawking Ansatz, we 
consider the case when $U = B\times \R$, with $B$ a contractible open subset of $\R^2$  --- in particular,
the $S^1$-bundle $X$ over $U$ is topologically trivial.  Set $y_1 = u_1$, $y_2 = u_2$, so then $y = y_1 +
iy_2$ is a complex coordinate on $B$.  We will
see below how the hyperk\"ahler structure on $X$ gives rise to 
a complex structure on $X$ under which the function $y$ is holomorphic, i.e. the map $X\to B$ is holomorphic.
Moreover, if we pass to the universal cover $\tilde X$ on $X$, we can construct a holomorphic coordinate $x$ 
(depending on a choice of  holomorphic section of $\tilde X$ over $B$) such that the holomorphic 2-form is
just $dx\wedge dy$.  This in turn  enables us to identify $\tilde X$ with $\T ^* _B $ over $B$, with $x,y$ 
then being holomorphic canonical coordinates on $\T ^* _B$, where the identification depends on our choice of 
holomorphic section.  The $S^1$-action on $X$ yields an $\R$-action on
$\T ^* _B $, which we shall see is just translation on $x_1 = \Re \, x$, and so $X$ is isomorphic to $\T ^* _B
/\Z$.  The K\"ahler form  provided by the Gibbons--Hawking Ansatz yields a K\"ahler form $\omega$ on $\T ^*
_B$, corresponding of course to a Ricci-flat metric, and for which the functions $W$ and $b$ are independent
of $x_1$.  The  K\"ahler form therefore descends to $\T ^* _B /\Z$, and is invariant under the obvious
$S^1$-action.

Conversely, we shall see that any Ricci flat, $S^1$-invariant K\"ahler structure on $\T ^* _B /\Z$ of the
above type (i.e. we have $x,y$ holomorphic canonical coordinates on $\T ^* _B $ over $B$, for which $W$ and $b$
are independent of $x_1$) does in fact arise from the Gibbons--Hawking Ansatz in the way that has just been
described.  Moreover, Gibbons--Hawking coordinates $u_1 ,u_2 , u_3 $ and the connection form $\theta$ on
$X$ may be recovered  from the holomorphic canonical coordinates $x,y$ on $\T ^* _B $.  Here we
have $u_1 =y_1$, $u_2 =y_2$, and $u_3$ determined up to a constant. \medskip

We now give the details for the construction.  We have $U = B\times \R$, with $B$ a contractible open subset
of $\R^2$, and we set $y = y_1 + iy_2$, a complex coordinate on $B$.
Then $dy_1 + i dy_2 = d y$, and from the Gibbons--Hawking Ansatz equations we observe that
$dy \wedge d(\theta _0  - iVdu_3) = 0$.  By the theorem on integrability of almost complex structures,
$\Omega=(\theta_0-iVdu_3)\wedge dy$ is a holomorphic 2-form for an integrable
almost complex structure $J$ on $X$, and
locally  there exists a holomorphic coordinate $z$ such that $dz=(\theta _0-iVdu_3)\ \mod\ dy$.  Moreover it 
is then clear that $z$ is determined up to a holomorphic function of $y$, and that locally the holomorphic 
coordinates recover the (integrable) complex structure.  We now pass to the universal cover 
$\tilde X$ of $X$, topologically $B\times \R^2$, together with its integrable complex structure 
$\tilde J$ obtained from $J$ (from now on, we shall work on $\tilde
X$, but omit tildes from forms and functions pulled back from $X$).  We note that the complex structure is
invariant under the $\R$-action 
on $\tilde X$ induced from the given $S^1$-action on $X$.  The (global) form $\theta _0 - iVdu_3$ restricts
down to a holomorphic 1-form on each fibre, locally just $dz$.  Therefore, by integrating 
$\theta _0 - iVdu_3$ along paths in the
fibre from some fixed point, we obtain a holomorphic coordinate on the fibre, which locally (up to a constant
depending on the choice of base point) will coincide with $z$.

In order to get a global 
holomorphic coordinate $x$ on $\tilde X$, we choose a holomorphic section of $\tilde X$ over $B$ (such sections 
always exist), which will then be regarded as giving  
the required base point in each fibre for the path integration.  In this way, we obtain a global 
holomorphic function $x$ on  $\tilde X$ such that $x,y$ are holomorphic coordinates everywhere, and where 
$x$ is uniquely determined up to a holomorphic function of $y$ (corresponding to the choice of holomorphic 
section).  By construction, the global holomorphic coordinates $x,y$ on $\tilde X$ 
realize the almost complex structure, with $y$ a holomorphic coordinate on the base and $x$ a holomorphic 
coordinate on the fibres. 
Moreover $\Omega =-\omega _1 - i \omega _2 =  dx\wedge dy$,
and  so we can identify $\tilde X \to B$ with ${\cal T}^*_B \to B$ (with holomorphic canonical coordinates,
as described in \S 1), where the chosen holomorphic section of $\tilde X$ over $B$ is identified the zero 
section of the holomorphic cotangent bundle.  
Choosing a section of $\tilde X$ over $U= B \times {\bf R}$ 
enables us to consider the coordinate $t$ on $S^1$ as a coordinate
on the fibres; the above derivation of the holomorphic coordinate $x$ 
then shows that its real part $ x_1 = {t\over 2\pi } + g( y_1 , y_2 , u_3 )$,
for some function $g$, and that the action of ${\bf R}$ on $\tilde X$ 
is the obvious one given by translation on $x_1$. Explicitly
$X$ is obtained as a quotient of $\tilde X$ under the action of $\Z$ given by $x_1 \mapsto x_1 +1$.

Since $dx = \theta _0 - iVdu_3\ \mod\ dy$,  
there exists a complex-valued function $b$ on $\tilde X$ such that 
$dx + bdy = \theta _0 - iVdu_3$.  
Also
$$ (dx + bdy)\wedge \overline {(dx + bdy)} = 2iV \theta _0 \wedge du_3.$$
We now set $W= V^{-1}$ and calculate the K\"ahler form $\omega_3$ in terms of the holomorphic coordinates:
$$\omega _3 = du_3\wedge\theta_0+V du_1 \wedge du_2 = {i\over 2}( W (dx + bdy)\wedge \overline {(dx +
bdy)} + W^{-1}dy \wedge d\bar y),$$
which we observe has the same form as our original general
formula for $\omega$ in holomorphic canonical coordinates.   
Since we started with a Ricci-flat metric, 
the  previous equations for Ricci-flatness (2.1) and (2.2) which we derived are then automatically
satisfied. 

The next point is to observe that $W$ and $b$ are independent of $x_1$,
the real part of $x$. To see this,
recall now that $\tv = W(dx +bdy) = (V^{-1}\theta _0 
-idu_3)$ and $\th = dy$ is a globally defined coframe for the holomorphic cotangent bundle of $\tilde X$.  In
particular, since the imaginary part of $\tv$ is $-du_3$,
we will have that the imaginary part of $d\tv$ is zero.  We first calculate 
$$\partial \tv = W^{-1} \partial W \wedge \tv + W \partial b \wedge \th = (W^{-1} \ph W - W \pv b )\th \wedge \tv
. \leqno{(2.3)}$$  From this it is seen that equation ($2.2'$) is just the statement that $\partial \tv =
0$.  We next calculate 
$$\eqalign{  \pc \tv & =  W^{-1} \pc W \wedge \tv + W \pc b \wedge \th  \cr 
& =  W^{-1} \pvc W \ \tvc \wedge \tv + W^{-1} \phc W \ \thc \wedge \tv + W \pvc b \ \tvc \wedge \th 
+ W \phc b \ \thc \wedge \th.\cr}\leqno{(2.4)}$$
If then equation ($2.1'$) also holds, it is easily checked that the imaginary part of $\pc \tv$ is zero 
if and only if $\pv W =-\pvc W$ and $\pv b =-\pvc b$, that is $W$ and $b$ are independent of $x_1$, 
the real part of $x$. Thus $b$ is invariant under the ${\bf R}$-action,
that is $b$ is the pull-back of a function from $U$.

Conversely, if we start from a 
Ricci-flat, $S^1$-invariant K\"ahler metric on $\T ^* _B /\Z$ of the
above type (i.e. we have $x,y$ holomorphic canonical coordinates on $\T ^* _B $ over $B$, for which $W$ and $b$
are independent of $x_1$), we can pass to the universal cover $\tilde X =  \T ^* _B $ over $B$.  The above
construction then reverses.  We set $\phi$ to be the imaginary
part of $\tv = W(dx +bdy)$.  Clearly $\phi$ is invariant under the given
$\R$-action on $\tilde X$.
Reversing the derivation of the previous paragraph ensures that $d\phi =0$ on 
$\tilde X$, and so there is a global function $u_3$ with $\phi = - du_3$, where $u_3$ is  
invariant under the action of $\R$, and is determined up to a constant.  We set $V = W^{-1}$,  
$\theta _0 = (dx + bdy) - iV\phi$ and $\theta = 2\pi i \theta _0$; thus both $V$ and $\theta$ 
are also invariant under the action of $\R$.
It is straightforward now to verify that we get back the above form of 
the Gibbons--Hawking Ansatz, with $u_1 = y_1$ and $u_2 = y_2$, and where $U = B\times \R$ is the
quotient of $\tilde X$ by the $\R$-action. The periodicity of this 
$\R$-action  then yields an $S^1$-bundle $X$ over $U$ (to which $V$ and $\theta$ descend, and
on which the corresponding
$S^1$-action leaves $V$ and $\theta$ invariant).\medskip

Finally, we calculate (for use in \S 4) the differential $p_*$, where 
$p:\tilde X\rightarrow U= B\times
{\bf R}$ is the natural projection. Using the expression
$du_3={i\over 2}W(dx-d\bar x)+{i\over 2}W(bdy-\bar bd\bar y)$, we obtain
$$\eqalign{p_*\px&={iW\over 2}\partial_{u_3}\cr
p_*\pxc&={-iW\over 2}\partial_{u_3}\cr
p_*\py&=\py+{iW\over 2}b\partial_{u_3}\cr
p_*\pyc&=\pyc-{iW\over 2}\bar b\partial_{u_3}\cr}\leqno{(2.5)}$$
Thus
$p_*\ph=\py$, $p_*\pv={i\over 2}\partial_{u_3}$.

Also, as $W$ and $b$ can be thought of as functions on $B\times {\bf R}$,
being independent of $x_1$, the formula $(2.2')$ translates into
$$-\py V={i\over 2}\partial_{u_3} b. \leqno{(2.6)}$$
Thus $b$ can be calculated as
$$b(y,u_3)=\sigma(y)+\int 2i\py V du_3\leqno{(2.7)}$$
where $\sigma(y)$ is some constant of integration.

\bigskip

\noindent {\bf $S^1$-invariant Ricci flat metrics on elliptic fibrations}

We shall be most interested in the transformation described above when  $V$ and $\theta$ are
themselves periodic in $u=u_3$.  
The hyperk\"ahler metric descends to
one on the corresponding $S^1$-fibration over $Y=B\times S^1$ if and only if the three 2-forms 
$\omega _1 , \omega _2 , \omega _3$ are invariant under changing $u$ by a period, which in turn is
saying that the periodicity in $u$ is independent of $y$.  We shall 
now change notation and denote this $S^1 \times S^1$ fibration  over $B$ by $X$ (the universal cover
$\tilde X$ being the same as before).  
Since the restriction of the K\"ahler form $\omega _3$ to a fibre $X_y$ is just $du \wedge \theta _0 
=du \wedge dt/{ 2 \pi }$, the volume of any fibre is just the periodicity 
in $u$.
Changing coordinates to the holomorphic coordinates of Construction 2.6, 
we obtain 
a holomorphic map $f : X \to B$ to a contractible open subset $B$ of
$\C$, whose fibres are elliptic curves.  Having chosen a holomorphic section, we obtain 
holomorphic canonical coordinates $x,y$  on
the corresponding line bundle $\tilde X$ over $B$, where the holomorphic 2-form 
$\Omega = dx\wedge dy$, and where 
the K\"ahler form $\omega$ (as defined by the usual formula) determines a hyperk\"ahler metric on
$X$.  Moreover, both $W$ and $b$ are independent of $x_1$.

The periods of the above elliptic fibration have a basis $\{ 1, \tau (y) \}$, for some
holomorphic function $\tau$ of $y$.
If we wish to have an explicit formula for $\tau (y)$, we take a basis of homology $\{ \gamma _1 ,
\gamma _2 \}$, where $\gamma _1$ is an $S^1$ in a fibre $X_y$ of $X \to B$ given by the orbit of the 
$S^1$-action, and $\gamma _2$ is an $S^1$ in  $X_y$ mapping isomorphically to $\{ y\} \times S^1 
\subset Y$.  Restricted to the fibre $X_y$, we have $dx = \theta _0 - iVdu_3$;  one of the periods is then
$$\int _{\gamma _1} dx = \int _{\gamma _1} \theta _0 = 1,$$ as already observed, whilst the other period
$$ \tau (y) = \int _{\gamma _2} dx = \int _{\gamma _2} \theta _0 - i\int _{\gamma _2} Vdu_3.$$
By choosing the appropriate orientation for $\gamma _2$, we may also assume that $\Im \,\tau (y) >0$.\medskip

If we have  such a holomorphic elliptic fibration $f : X \to B$ and Ricci-flat metric
(independent of
$x_1$), we shall refer to it as an \it $S^1$-invariant Ricci-flat metric (on $X$) in canonical form.\medskip

\rm Conversely, if we are given such an $S^1$-invariant Ricci-flat metric on $X$, we saw above how
this does indeed arise from the Gibbons--Hawking Ansatz.  Moreover, in this case, we also have that 
$V$ and $\theta$ are periodic in $u$, with the period in $u$ being constant, namely the volume of the 
elliptic fibres of $f : X \to B$.
\medskip

\noindent {\bf Remark 2.7.} A particular case of an  
$S^1$-invariant Ricci-flat metric in canonical form
is a semi-flat metric: 
 Given, locally, two
periods $\tau_1$ and $\tau_2$, these should be interpreted as 1-forms
on $B$, i.e. are $\tau_1 dy$, $\tau_2dy$. We can then locally
replace $y$ with a holomorphic function $g$ on an open set $U$ such
that $dg=\tau_1 dy$, and thus can assume $\tau_1=1$. Then in these
coordinates, the semi-flat metric coincides with the Gibbons--Hawking
metric obtained by taking $V=\im\tau_2/\epsilon$ on 
$U\times {\bf R}/\epsilon\boldz$.  
We can then use the formula of Remark 2.4 to compute $\| R\|^2$ for a semi-flat metric (which will coincide
with the value calculated via Remark 2.1). Thus
$$\|R\|^2={1\over 2}V^{-1}\Delta\Delta V^{-1}
={\epsilon^2\over 2} (\im\tau_2)^{-1}\Delta\Delta (\im\tau_2)^{-1}.$$
In particular, $\|R\|^2\rightarrow 0$ as $\epsilon\rightarrow 0$.\bigskip

Returning now to the set-up in Question 1.2;
away from the singular fibres, we expect that, as the volume $\epsilon$ of the fibres tends to zero,
the metric (suitably normalized) will approach a semi-flat one.  
This expectation is motivated by the following result, which proves a slightly weaker version of the 
expected convergence for the $S^1$-invariant Ricci-flat case, purely by local considerations, as 
a consequence of Harnack's inequality for harmonic functions.  Whilst we don't expect a purely local 
proof of convergence in general (i.e. not assuming the $S^1$-invariance of the metrics), the main result
of this paper (Theorem 5.6) will prove a very strong form of the expected convergence to a semi-flat 
metric (locally over the base) by means of global methods.

\proclaim Proposition 2.8.  Let $\pi : X \to B$ be an elliptic fibration with 
periods $\{ 1,\tau (y)
\}$,  over the open disc $B$ of radius $R$ in
$\C$, with $\Im\, \tau (y) >0$, and let 
$B_0 \subset B$ denote a smaller disc of radius $R_0 < R$.  
Suppose we have a sequence on $X$ of 
$S^1$-invariant Ricci-flat metrics $g_i$ in canonical form 
(and with constant volume form), 
for which the volume $\epsilon _i := \epsilon (g_i )$ of the fibres tends to
zero as $i\to \infty$. 
Then on $\pi ^{-1} (B_0)$ we have 
$W_i := W(g_i ) \to
0$ uniformly as $i\to \infty$.  On a fixed fibre, with periods $\{ 1,\tau \}$,
we have the stronger statement that $ \epsilon _i ^{-1}  W_i\ \Im\, \tau \to 1$
uniformly as $i\to \infty$.\par

\noindent \it Proof. \rm   Our assumption that the volume form is constant 
ensures that we can fix the holomorphic canonical coordinates coordinates $x,y$ and the holomorphic 
2-form $\Omega = dx\wedge dy$ independent of $i$.  
We now transform the coordinates to Gibbons--Hawking coordinates;
 the first claim is equivalent to $ W_i \to 0$ uniformly on $B_0 \times \R $.  If the volume of the fibres
is $\epsilon_i$, then the periodicity in $u$ is $\epsilon_i$.  
Fix $R_1$ with $R_0 < R_1 < R$; then for $\epsilon_i <<1$,
the ball $\tilde B_1$ in $U = B\times \R$ with centre the origin and radius $R_1$ will contain the set 
$B_0 \times [0, \epsilon_i ]$, so it will suffice to show that $ W_i \to 0$ uniformly on $\tilde B_1$ 
as $i\to \infty$. Fix $i$ for the moment so that $\epsilon_i$ is
sufficiently small as above, and drop the subscript for convenience. 
Let $\tilde B$ denote the ball of radius $R$, with centre the origin.  Recall now that $W= V^{-1}$.
Given that $V$ is harmonic on $\tilde B$, this is precisely the situation in which we can 
apply the strong form of Harnack's inequality, as stated in Problem 2.6 on page 29 of [11], namely
that for any point $P\in \tilde B_1$,
$${(1 - R_1/R) \over (1+ R_1/R )^2 } \le {V(P) \over V(0)} \le {(1 + R_1/R) \over (1- R_1/R )^2 }.$$
Thus, for $P \in \pi ^{-1} (B_0)$, the ratio $W(P)/W(0)$ is bounded above and below by 
appropriate positive constants.
For each $y \in B_0$, we can calculate the volume of the fibre $X_y$ as
$$\epsilon = \int _{X_y} W dx_1 \wedge dx_2 = \int _0 ^ {\Im\, \tau (y)} W dx_2.$$  For $R_0$ fixed and 
for $y\in B_0$, we also have that 
$\Im\, \tau (y)$ is bounded above 
and below by appropriate positive constants.  On any fibre $X_y$ with $y\in B_0$, 
we can find a point at which $W$ takes the average value on the fibre, namely $\epsilon / \Im\, \tau (y)$.
Putting all these facts together yields the claim that $W_i \to 0$ uniformly on $\tilde B_1$ as $i\to \infty$.

For the stronger statement on a fixed fibre, we can assume that the fibre is $X_0$, and that $W$ 
takes the average value 
$\epsilon / \Im\, \tau $ at the centre $0$ of the ball $\tilde B$.  If we take a concentric ball $\tilde
B(r)$ of small radius $r$, it will still contain all of $\{ 0 \} \times [0, \epsilon ]$, provided 
$\epsilon < r$.  Harnack's inequality then yields 
$${(1 - r/R) \over (1+ r/R )^2 } \le {W(0) \over W(P)} \le {(1 + r/R) \over (1- r/R )^2 }$$
for all $P \in X_0$.  By taking $r$ arbitrarily small, these upper and lower bounds are arbitrarily 
close to 1, and hence $ \epsilon _i^{-1}  W_i \ \Im\, \tau \to 1$ 
uniformly on $X_0$. $\bullet$

\bigskip
{\hd \S 3. The Ooguri--Vafa metric.}

The aim of this section is to describe a certain hyperk\"ahler metric on
a neighbourhood of each singular fibre in our elliptically fibred K3 surface, and to derive
various estimates associated with this metric.  If the fibres are assumed to have volume
$\epsilon$, then away from the singular fibre, this metric decays
very rapidly, for small $\epsilon$,
to a semi-flat metric.  We shall assume throughout that we only have
singular fibres of Kodaira type 
$I_1$, and so locally around the singular fibre, one of the periods is invariant under monodromy (and
in fact, by an appropriate choice of holomorphic coordinate $y$ on the base, 
may be taken to be constant, value
$1$), whilst the other period will  be multivalued and tend to infinity.  The metric we define will
be an $S^1$-invariant  
metric (as described in the previous section) on the smooth part of the fibration, and will
be most conveniently described in the Gibbons--Hawking coordinates.

The metric we describe was first written down (in a slightly different form) by Ooguri and Vafa [29],
and so will be referred to as the \it Ooguri--Vafa \rm metric.   In \S 4, we shall start 
with an Ooguri--Vafa metric in a neighbourhood 
of each singular fibre; by appropriately twisting these metrics, we'll show that they may be 
glued with a semi--flat metric away from the singular fibres, hence obtaining a global metric, which
is Ricci-flat away from the gluing regions, and which represents the correct K\"ahler class.  For 
small $\epsilon$, it is these metrics which approximate very accurately the global Ricci flat metric
with the given K\"ahler class.

Before launching into the technical details, we shall briefly describe the basic idea behind the
construction of the Ooguri--Vafa metric, which, given the description of the Gibbons--Hawking Ansatz
in \S 2, should strike the reader as very natural.
The harmonic function $V$ we use will be periodic in $u$ of period $\epsilon$ (the 
volume of the fibres), but have Taub--NUT type singularities on the fibre $y=0$ at
the points $u\in \epsilon \Z$.

We take $U = D \times \R \setminus 
\{ 0\} \times \epsilon\Z $, or more precisely its quotient by $\epsilon \Z$, where $D \subset \C$
is an open disc centred at the origin.  
We denote by $y_1 ,y_2$ the
coordinates on $D \subset \C$, and by $u$ the coordinate on
$\R$.   We want to write down $V$ harmonic on $U$, periodic in $u$ and with singularities 
of the correct type at the points $\{ 0\} \times \epsilon\Z $.  For instance, around zero, $V$
should behave like a  harmonic function
plus a term ${1\over 4\pi |{\bf x}|}$, from which it will follow that 
the total space $X$ of the
$S^1$-fibration over $U$ extends (by adding a single point) to a manifold $\bar X$ mapping onto 
$\bar U = D \times \R /\epsilon \Z$. In addition,
the hyperk\"ahler metric extends to $\bar X$.
We are led therefore to take $V = V_0 + f(y_1 ,y_2 )$,
where $f$ is a harmonic function in $y_1 ,y_2$ on $D$, and 
$$ V_0 = {1\over 4\pi} \sum _{-\infty} ^{\infty} \left( {1\over \sqrt{(u+n\epsilon )^2 + y_1 ^2 + y_2 ^2}} 
- a_{|n|}\right) ,$$
where $a_n = {1\over n\epsilon }$ ($n >0$), thus ensuring appropriate convergence, and $a_0$ is
chosen appropriately to ensure that the periods do not change as we change $\epsilon$ --- that is, 
we are defining metrics on a fixed elliptic fibration.  This choice of $a_0$ also ensures, 
on a fixed annulus in
$D$, that $\epsilon V_0 \sim - {1\over 2\pi } \log r$ as $\epsilon \to 0$, where 
$r^2 = y_1^2 + y_2 ^2$.  In general, the periods around an $I_1$ fibre may be 
assumed to be $1$ and 
$\tau (y) =   {1\over 2\pi i } \log y + i h(y) $, where $h$ is holomorphic in $y = y_1
+ iy_2$, and these may be achieved in our construction by taking $V = V_0 + f(y_1 ,y_2)$, where $f$
denotes the real part of $h$.

We now give the technical details.  

\proclaim Lemma 3.1. Let
$$T_j={1\over 4\pi} \sum_{n=-j}^{j} \left(
{1\over \sqrt{(u+n\epsilon)^2+y_1^2+y_2^2}}-a_{|n|}\right)$$
where
$$a_n=\cases{1/n\epsilon& $n\not=0$\cr
2(-\gamma+log (2\epsilon))/\epsilon&$n=0$\cr}$$
and $\gamma$ is Euler's constant.
Then 
\item{(a)} the sequence $\{T_j\}$
converges uniformly on 
compact sets in $D\times {\bf R}-\{0\}\times\epsilon\Z$
to a harmonic function $V_0$. Here $D\subseteq {\bf C}$ is the unit
disc centred at the origin.
\item{(b)} $V_0$ has an expansion, valid when $ |y| \not =0$, 
$$V_0=-{1\over 4\pi \epsilon}\log |y|^2+\sum_{m=-\infty\atop m\not=0}^{m=
\infty} {1\over 2\pi\epsilon}e^{2\pi imu/\epsilon} K_0(2\pi|my|/\epsilon)$$
where $y=y_1+iy_2$ and 
$K_0$ is the modified Bessel function. (See [3], pg. 374.)
\item{(c)} There exists a constant $C$ such that for any
$0< r_0<1$, there exists an $\epsilon_0>0$ such that for all $\epsilon
<\epsilon_0$, $|y|>r_0$,
$$\left| V_0+{1\over 4\pi\epsilon}\log|y|^2\right|\le
{C\over\epsilon} e^{-2\pi|y|/\epsilon}.$$
\item{(d)} If $r\le 1$, and $f$ is a harmonic function on the disc
$D_r$ of radius $r$ such that $f(y)-{1\over 4\pi}\log|y|^2
>0$ for $|y|\le r$, then there exists an $\epsilon_0$ such that for
all $\epsilon<\epsilon_0$,
$$V_0+f(y)/\epsilon>0$$
in $D_r\times {\bf R}$.

\noindent {\it Proof.} (a) Let $p$ be the smallest integer greater than $\epsilon^{-1}\sqrt{
1+\epsilon^2}$. Then for $0\le u\le\epsilon$, $y_1^2+y_2^2\le 1$,
we have
$${1\over \sqrt{(u+n\epsilon)^2+y_1^2+y_2^2}}>a_{|n|+p}$$
for all $n$.
Let 
$$R_j=
{1\over 4\pi} \sum_{n=-j}^{j} \left(
{1\over \sqrt{(u+n\epsilon)^2+y_1^2+y_2^2}}-a_{|n|+p}\right).$$
Then for $j>2p$, 
$$T_j-R_j={1\over 4\pi}\left(-a_p-a_0-2\sum_{n=1}^{p-1} a_n
+2\sum_{n=j-p+1}^j a_{n+p}\right).$$
Put $C(\epsilon)={1\over 4\pi}(-a_p-a_0-2\sum_{n=1}^{p-1}a_n)$. Note that
$\sum_{n=j-p+1}^j a_{n+p}\rightarrow 0$ as $j\rightarrow\infty$, 
so if $R_j$ converges uniformly on compact sets to a harmonic
function $R$, then $T_j$ converges to a harmonic function
$R+C(\epsilon)$. Now for $0<u<\epsilon$, $y_1^2+y_2^2<1$, $R_j$ is
a monotonically increasing sequence of harmonic functions (since all terms
are positive). Furthermore, it is easy to check that, say, the 
sequence $R_j$ is bounded at $u=\epsilon/2,y_1=y_2=0$. Thus by the Harnack
convergence theorem, (Theorem 2.9, [11]), the $R_j$ converge
uniformly on compact subsets to a harmonic function $R$, and
$V_0=R+C(\epsilon)$. Since $R$ is positive, we see $V_0>C(\epsilon)$.
For $u=0,\epsilon$, we  merely omit the term which blows up and then repeat
the previous argument.

(b) The part which requires care is the constant term of the Fourier expansion,
i.e computing ${1\over \epsilon}\int_0^{\epsilon} V_0 du$. To do so, consider
the following variant on the $T_j$:
$$S_j={1\over 4\pi} \sum_{n=-j}^{j} \left(
{1\over \sqrt{(u+n\epsilon)^2+y_1^2+y_2^2}}-b_{|n|}\right)$$
where
$$b_n=\cases{(\log (n+1)-\log n)/\epsilon& $n\not=0$\cr
0&$n=0$\cr}$$
Then 
$$T_j-S_j={1\over 4\pi}\left({2\over\epsilon}
\log(j+1)-a_0-{2\over\epsilon}\sum_{n=1}^j {1\over n}\right).$$
As $j\rightarrow\infty$, this converges to
$${1\over 4\pi} (-{2\gamma\over\epsilon}-a_0)= 
-{1\over 2\pi\epsilon}\log(2\epsilon).$$
Now we calculate
$$\eqalign{4\pi \int_0^{\epsilon} S_jdu
=&\sum_{n=-j}^j \int_0^{\epsilon} \left( {1\over \sqrt{
(u+n\epsilon)^2+|y|^2}}-b_{|n|}\right) du\cr
=&
\sum_{n=-j\atop n\not=0}^j(-\log(|n|+1)+\log|n|)+
\sum_{n=-j}^j
\int_{n\epsilon}^{(n+1)\epsilon} {1\over\sqrt{u^2+|y|^2}}du \cr
=&
\int_{\epsilon}^{(j+1)\epsilon}\left( {1\over \sqrt{u^2+|y|^2}}-{1\over u}\right)du
+\int_{-j\epsilon}^{0}\left( {1\over \sqrt{u^2+|y|^2}}+{1\over u-\epsilon}\right)du\cr
&+\int_0^{\epsilon} {1\over \sqrt{u^2+|y|^2}}du\cr
=&\log\left( u^{-1}\left( u+\sqrt{u^2+|y|^2}\right) \right)
\bigg|_{\epsilon}^{(j+1)\epsilon}
+\log\left( |u-\epsilon|\left( u+\sqrt{u^2+|y|^2}\right)\right)
\bigg|_{-j\epsilon}^{0}\cr
&+\log\left( u+\sqrt{u^2+|y|^2}\right)
\bigg|_0^{\epsilon}.\cr
}$$
Evaluating this and letting $j\rightarrow\infty$, one obtains
$${1\over\epsilon}\lim_{j\rightarrow\infty}\int_0^{\epsilon} S_jdu
={1\over 4\pi\epsilon}(\log 2+2\log\epsilon-\log{|y|^2\over 2})$$
from which we conclude that
$${1\over\epsilon}\int_0^{\epsilon} V_0du=-{1\over 4\pi\epsilon}
\log|y|^2.$$

To compute the other terms in the Fourier expansion, we just need
to calculate 
$$\eqalign{
{1\over \epsilon}\int_0^{\epsilon} V_0e^{2\pi imu/\epsilon} du
&={1\over 4\pi\epsilon}\int_{-\infty}^{\infty} {e^{2\pi imu/\epsilon}
\over\sqrt{u^2+|y|^2}} du\cr
&={1\over 2\pi\epsilon}\int_0^{\infty} {\cos(2\pi mu/\epsilon)\over
\sqrt{u^2+|y|^2}}du\cr
&={1\over 2\pi\epsilon}\int_0^{\infty} {\cos(2\pi|my|v/\epsilon)
\over \sqrt{v^2+1}} dv\cr
&={1\over 2\pi\epsilon}K_0(2\pi|my|/\epsilon).\cr}$$
The last equality follows from [3], page 376, formula 9.6.21.

(c) By [3], 9.8.6, there exists a constant $C_1$ such that
$\sqrt{x}e^xK_0(x)\le C_1$ for $x\ge 2$. (In fact $C_1\le 2$). In
particular, $K_0(x)\le C_1e^{-x}$ for $x\ge 2$.
Thus
$$\eqalign{\left|\sum_{m=-\infty\atop m\not=0}^{\infty}
{1\over 2\pi\epsilon}e^{2\pi i m u/\epsilon}K_0(2\pi|my|/\epsilon)\right|
&\le {C_1\over \pi\epsilon}\sum_{m=1}^{\infty} e^{-2\pi|my|/\epsilon}\cr
&={C_1\over\pi\epsilon} {e^{-2\pi|y|/\epsilon}\over 
1-e^{-2\pi|y|/\epsilon}}\cr}$$
for $2\pi|y|/\epsilon\ge 2$. From this follows (c).

(d) By the maximum principal, the minimum value $M$ of $f$ occurs
on the boundary of $D_r$. On the other hand, for
fixed $u$, it is clear $V_0$ is monotonically decreasing in $|y|$.
Thus the minimum value of $V_0+f/\epsilon$ must occur on
$(\partial D_r)\times {\bf R}$. But taking $r_0<r$, by (c)
there exists an $\epsilon_0$ such that for all $\epsilon<\epsilon_0$,
$$\left|V_0+{1\over 4\pi\epsilon}\log|y|^2\right|<-{1\over 4\pi\epsilon}
\log r^2+M/\epsilon$$ 
whenever $|y|=r$. Thus $V_0+f/\epsilon$ is positive on
$\partial D_r\times{\bf R}$ for $\epsilon<\epsilon_0$, hence
$V_0+f/\epsilon$ is positive on $D_r\times {\bf R}$. $\bullet$
\bigskip

With this rather technical lemma out of the way, we may now proceed to 
the construction of our metric, using the Gibbons--Hawking Ansatz formalism, as
developed in \S 2.  Suppose $D_r \subset \C$ is the disc of radius 
$r < 1$, centre the origin, and $f : \bar X \to D_r$ an elliptic fibration,
with singular fibre over the origin of type $I_1$.  Let
$\bar Y=D_r\times{\bf R}/\epsilon\Z$ and 
$Y=(D_r\times {\bf R}-\{0\}\times
\epsilon\Z)/
\epsilon\Z$.  It is straightforward to check that there is an 
induced map $\bar \pi : \bar X \to \bar Y$ of $C^\infty$ manifolds, which
restricts to an $S^1$-bundle $\pi : X \to Y$ with Chern class $\pm 1$, the sign 
dependent on the choice of orientation for the fibre.  For further justification 
of these statements, the reader is referred to [15], Example 2.6 (1).  The plan 
now is to define a hyperk\"ahler metric on $X$ via the Gibbons--Hawking Ansatz 
applied to $\pi : X \to Y$, and then check that it extends to a hyperk\"ahler 
metric on $\bar X$.

\bigskip

\proclaim Proposition 3.2. With the notation as above, let $h (y)=f(y_1,y_2)
+ig(y_1,y_2)$ be a holomorphic function on $D_r$, so that $-{1\over 4\pi}
\log |y|^2+f(y_1,y_2)>0$ on $D_r$. 
Let $V_0$ be the harmonic function on $Y$ defined in Lemma 3.1, and $V=V_0+f(y_1,y_2)/\epsilon $,
with $\epsilon$ chosen small enough so that $V>0$ on $Y$. Then
there exists a connection 1-form $\theta$ on $X$ such that
$d\theta/2\pi i=*dV$, and this defines 
a hyperk\"ahler metric on $X$ with
$$\eqalign{
-\re\Omega&=dy_1\wedge \theta/2\pi i+Vdy_2\wedge du\cr
-\im\Omega&=dy_2\wedge \theta/2\pi i+Vdu\wedge dy_1\cr
\omega&=du\wedge \theta/2\pi i+Vdy_1\wedge dy_2.\cr}$$
These forms extend to $\bar X$, giving a hyperk\"ahler
metric on $\bar X$, and a holomorphic elliptic fibration 
$\bar X\rightarrow D_r$  with periods $1$
and ${1\over 2\pi i}\log y +i h(y) + C$, for some real constant $C$.
By appropriate choice of $\theta$, this constant $C$ may be taken to be zero.\par

\noindent {\it Proof.} Since $V$ is harmonic, recall that $*dV$ is closed.
Taking a sphere $S^2$ of radius $<\epsilon$ centred at $0\in D_r\times
{\bf R}$, we have 
$$\int_{S^2} *dV=\int_{S^2} *d\left({1\over 4\pi\sqrt{u^2+y_1^2+y_2^2}}\right)$$
since all other terms in $*dV$ are defined at 0, and hence are exact
on an $\epsilon$-ball around $0$, and therefore do not contribute to the
integral. In Example 2.5, it was however observed that this latter integral is 
$\pm 1$ (depending on the orientation of the sphere). Thus,
since a connection form $\theta$ can be found such that $id\theta/2\pi$
is any desired representative of $c_1$, we can find a connection form
$\theta$ such that $d\theta/2\pi i=*dV$.  Applying now the 
Gibbons--Hawking Ansatz construction described in \S 2,  
we obtain a hyperk\"ahler metric on $X$, with the forms $\re\Omega$, $\im\Omega$ and 
$\omega$ as described in the Proposition.

To see that these forms extend to $\bar X$, focus on an $\epsilon/2$-ball
$B$ around $0$ in $\bar Y$. Then $\bar\pi^{-1}(B)\rightarrow B$ can
be identified with the map given in Example 2.5, restricted to the inverse
image of the $\epsilon/2$-ball in ${\bf C}^2$. Let $\theta '$ be the
connection form given in that example. Now $d(\theta)-d(\theta ')$
is the pull-back of an {\it exact} form on $B$, since all other terms
of $V$ besides the $n=0$ term are defined on $B$. Thus on $\bar\pi^{-1}
(B-\{0\})$
we can write
$\theta=\theta ' +\bar\pi^*\beta$ for a form $\beta$ defined on all of $B$.
Now consider for example on $\bar\pi^{-1}(B-\{0\})$
$$\eqalign{\omega&=du\wedge \theta/2\pi i+V dy_1\wedge dy_2\cr
&=du\wedge (\theta ' +\bar\pi^*\beta)/2\pi i+(1/4\pi\sqrt{u^2+y_1^2+y_2^2}+V')
dy_1\wedge dy_2\cr}$$
where $V'$ is a function defined everywhere on $B$. Thus we obtain
$$(du\wedge\theta '/2\pi i+(1/4\pi\sqrt{u^2+y_1^2+y_2^2})dy_1\wedge dy_2)
+du\wedge\beta+V' dy_1\wedge dy_2.$$
The first two terms were seen to extend to all of $\bar\pi^{-1}(B)$
in Example 2.5, and the last two terms are defined everywhere on 
$B$, so $\re\Omega$ extends to $\bar X$. Note that $(\re\Omega)^2\not=0$
at the singular point of the singular fibre, because $du\wedge\beta
+V' dy_1\wedge dy_2=0$ at that point.

Finally, we compute the periods.  Referring back to our discussion of 
$S^1$-invariant Ricci-flat metrics
in \S 2, one of the periods is constant, value 1.  The other period $\tau (y)$ is locally  
holomorphic in $y$, and given by 
$$\int _{\gamma _2} dx = \int _{\gamma _2} \theta _0  - i\int _{\gamma _2} Vdu,$$
where $\gamma _2$ is an $S^1$ in the fibre $X_y$ mapping isomorphically to $\{ y\} \times S^1 
\subset Y$.   Calculating the imaginary part of this,
$$\int _{\gamma _2} dx_2 =  - \int _{\gamma _2} Vdu = \pm({1\over 4\pi}\log |y|^2 -f(y_1,y_2)),$$
using the Fourier expansion for $V_0$ proved in Proposition 3.1 (b).
We choose the orientation of $\gamma_2$ to obtain the choice of sign to
be minus.
Then $\int _{\gamma _2} dx_1$ is necessarily locally a harmonic conjugate of 
$-{1\over 4\pi}\log |y|^2 +f(y_1,y_2)$, and so the period of 
$\gamma_2$ is $${1\over 2\pi i}\log y + i h(y) + C$$ for
some real constant $C$. 
Now $\theta _0$ may be modified by adding a term $adu$ ($a\in \R$) without changing the 
fact that $d\theta _0 = *dV$.  If $\theta _0$ is changed in this way, we have  
$$\eqalign{\int_{\gamma_2} dx_1 
&=\int_{\gamma_2} \theta _0 + a du\cr
&=a \epsilon +\int_{\gamma_2} \theta _0 .\cr}$$
We can therefore choose $a$ suitably to obtain $C=0$, and hence the periods as claimed. $\bullet$
\bigskip

\noindent \bf Remark 3.3. \hfill\break
\rm (1)\quad There is still some remaining flexibility over choosing $\theta$, as we can change
$\theta$ by the pull-back of a closed form from $\bar Y$.  This however need not worry us, since
in order to perform the gluing in \S 4, we will in any case need to twist the Ooguri--Vafa
metrics, the twist given as translation by an appropriate local section.\hfill\smallskip
\noindent (2)\quad  Recall that in the holomorphic canonical
coordinates $x,y$, the holomorphic 2-form 
$\Omega$ on $X$ is just $dx\wedge dy$, 
and so the complex structure on $X$ will be the one desired.
This 2-form extends uniquely to give the correct complex structure on $\bar X$.

\noindent {\bf Remark 3.4. A useful transformation.}

As $\epsilon \to 0$, the behaviour near the singular fibre of the Ooguri--Vafa metric is understood
best by making a change of variables.  The periods may be assumed to be 
$1, {1\over 2\pi i}\log y +ih(y)$, as in Proposition 3.2, and we take $V = V_0 + f(y)/\epsilon $.
We make the change of variables 
$s = u/\epsilon , v_1 = y_1/\epsilon , v_2 = y_2/\epsilon$.  Thus the disc of radius $\epsilon$ in
the complex $y$-plane corresponds to the unit disc in the complex $v$-plane.  If we now consider 
$V_0$ as a function of these new variables, we observe that
$$\epsilon 
 V_0 = {1\over 4\pi} \sum _{n=-\infty} ^{\infty} \left( {1\over \sqrt{(s +n)^2 + v_1 ^2 + v_2 ^2}} 
- c_{|n|}\right) ,$$
where $c_n = {1\over n}$ ($n >0$), and $$c_0 = 2(-\gamma + \log (2\epsilon )) = 
2(-\gamma + \log 2) + 2\log \epsilon .$$
So, if $\tilde V_0$ is the standard function $V_0$ in variables $s, v_1 , v_2 $ for $\epsilon = 1$,
we deduce that $$ \epsilon V = \tilde V_0 - 
{1\over 2\pi }\log \epsilon +f.$$  
Thus, if
we start with an Ooguri--Vafa metric with fibres of volume $\epsilon$ over the disc of radius $\epsilon$, 
make the change of variables described above, and then rescale the metric by $\epsilon ^{-1}$, we
obtain the Ooguri--Vafa metric over the unit disc, with fibres of volume one, corresponding to the 
harmonic function $\tilde V_0  + f - 
{1\over 2\pi }\log \epsilon $.  Thus the periods of the corresponding 
elliptic fibration are seen to be $1, {1\over 2\pi i}\log v +{1\over 2\pi i}\log \epsilon +ih(y)$.
This transformation lies behind the various estimates for diameters and curvature we 
derive below. We note here that in fact the formula given in [29] was
for $\epsilon V_0$, rather than $V_0$, except that the constant $a_0$
was not specified. The exact value for $a_0$ greatly influences the
behaviour of the metric as $\epsilon\rightarrow 0$, so this is quite
important.

To understand the metric for $|y| >\epsilon $, we can use the Fourier expansion for $V_0$ from 
Lemma 3.1 (b), and use the same change of variables as above.  Thus 
$$\epsilon V_0= -{1\over 4\pi }\log (v_1^2+v_2^2) -{1\over 2\pi }\log \epsilon +
\sum_{m=-\infty\atop m\not=0}^{m=
\infty} {1\over 2\pi }e^{2\pi ims} K_0(2\pi|mv|)$$
where $v=v_1+iv_2$.
\bigskip

\noindent {\bf Estimates for diameter and curvature}

We now consider a fixed elliptic fibration 
 $f : X' \to D'$  over a disc $D'$ of radius $a'<1$,
with singular fibre of Type $I_1$ over the origin, and which we assume extends to an elliptic
fibration over some larger disc.  We assume that the periods are of the form
$1, \tau (y)$, where $\tau (y) ={1\over 2\pi i}\log y +i h(y)$ as in Proposition 3.2.  We then 
wish to study 
sequences of Ooguri--Vafa metrics yielding the correct holomorphic 2-form $\Omega$, but with the
volume 
$\epsilon$ of the fibres tending to zero --- such metrics exist on $X'$ for small enough 
$\epsilon$ by Proposition 3.2.   We first ask about the diameters of the fibres.

\proclaim Proposition 3.5.  There exists a positive constant $C_1$ 
(independent of $\epsilon$) such that, 
 for metrics as above with fibre volume $\epsilon$, the diameters of the fibres over $D'$ are
bounded above by $C_1 (\epsilon\, \log {\epsilon ^{-1}})^{1/2}$.  Moreover, there exists a second 
constant $C_2$ such that the diameter $d(\epsilon )$ of the singular fibre is at least 
$C_2 (\epsilon\, \log {\epsilon ^{-1}})^{1/2}$.  

\noindent{\bf Remark 3.6.}
In particular, it follows that
$d(\epsilon ) \to 0$ as $\epsilon \to 0$. If however
we rescale the metric by
$\epsilon  ^{-1}$ as in Remark 3.4
to obtain fibres of volume one, then the diameter of the singular fibre is of
order $(\log {\epsilon ^{-1}})^{1/2}$, and therefore becomes arbitrarily large as $\epsilon \to 0$.
This then contrasts with the situation for a non-singular fibre,
where for sufficiently small $\epsilon$, the Ooguri--Vafa
metric near this non-singular
fibre is close to being semi-flat. Thus the diameter of the fibre in the
rescaled metric remains bounded.
\par

\noindent\it Proof. \rm To calculate the diameter of a fibre, we recall from \S 2 the formula 
for $ds^2$ in the Gibbons--Hawking Ansatz, namely 
$$ds^2 = V d{\bf u}\cdot d{\bf u} + V^{-1} \theta _0 ^2.$$  From this, it is clear that the diameter
of a fibre is at least $\int_0^{\epsilon/2} V^{1/2}du
={1\over 2}\int _0 ^{\epsilon} V^{1/2} du$.  Recall however that for all $y\ne 0$, there
exists a point on the fibre over $y$ at which $V = \im \tau (y)/\epsilon$, where now 
$\im \tau (y)$ is bounded below by a positive constant 
for $y\in D'$.  For some constant $C$ therefore, we have 
on each fibre $0\ne y \in D'$, a point at which $V^{-1/2} \le C \epsilon ^{1/2}$; by
continuity, this is also true for the singular fibre.  Thus each fibre over $D'$ contains 
an $S^1$ in the $S^1$-bundle (where $u$ is constant) of length at most $ C \epsilon ^{1/2}$, 
and hence the diameter of the fibre is  at most  ${1\over 2}\int _0 ^{\epsilon} V^{1/2} du 
+ C \epsilon ^{1/2}$.  Since $V = V_0 + f(y_1 ,y_2)/\epsilon $, and $|f|$ is bounded on $D'$ by some 
constant $A>0$, we 
have $$ \int _0 ^{\epsilon} V_0 ^{1/2} du - A^{1/2}\epsilon ^{1/2} \le \int _0 ^{\epsilon} V^{1/2} du \le 
\int _0 ^{\epsilon} V_0 ^{1/2} du + A^{1/2}\epsilon ^{1/2} .$$  Since $\int _0 ^{\epsilon} V_0 ^{1/2} du$ 
clearly takes its maximum when $y=0$, we are reduced to estimating $\int _0 ^{\epsilon} V_0 ^{1/2}
du$ on the singular fibre only, and showing that it is of order 
$(\epsilon \, \log {\epsilon ^{-1}})^{1/2}$.

We now let $\bar V_0$ denote the restriction of $V_0$ to the singular fibre, that is we take $y=0$.  
Making the substitution $s = u/\epsilon$ as above, we observe that, for $0< s <1$, 
$$4 \pi \epsilon \bar V_0 = \sum _{n=1} ^{\infty} \left( {1\over {s+n}} - {1\over n} \right) 
+ \sum _{n=1} ^{\infty} \left( {1\over {-s+n}} - {1\over n} \right) +  {1\over s} + 2\gamma -
2\log 2\epsilon .$$  We now quote formula 6.3.16 from [3], for the fact that, for $0 < s <1$, 
$$- \sum _{n=1} ^{\infty} \left( {1\over {s+n}} - {1\over n} \right) - \gamma = \psi (1+s) ,$$
where $\psi$ denotes the psi function.
Thus, for $0< s <1$, 
$$ 4 \pi \epsilon \bar V_0 = - \psi (1+s) - \psi (1-s)  + {1\over s} - 2\log 2\epsilon.$$
Using formula 6.3.15 from [3], we know that 
$$ - \left( \psi (1+s) + \psi (1-s) \right) = 2(1-s^2)^{-1} + 2\gamma - 2 
+ \sum _{n=1} ^{\infty} 2(\zeta (2n+1) -1)s^{2n},$$ where $\zeta$ denotes the usual zeta function.
Hence, for $0< s <1$, 
$$ 4 \pi \epsilon \bar V_0 = 2(1-s^2)^{-1} + {1\over s} - 2\log \epsilon + G + 2g(s) ,$$
with $G = - 2\log 2 + 2\gamma - 2$, and 
where $$g(s) = \sum _{n=1} ^{\infty} (\zeta (2n+1) -1)s^{2n}$$ has radius of 
covergence at least 2 (by inspection of the coefficients), and so defines a continuous 
(non-negative) function on $[0,1]$.  Now observe 
that $\int _0 ^{\epsilon} \bar V_0 ^{1/2} du = \epsilon ^{1/2} \int _0 ^1 (\epsilon \bar V_0 )^{1/2}
ds$.  The lower  bound now follows immediately by ignoring the first two terms in the expression for 
$ 4 \pi \epsilon \bar V_0$.  The upper bound follows by using the 
elementary fact that for $\alpha ,\beta $ non-negative real numbers, 
$(\alpha +\beta )^{1/2} \le \alpha^{1/2} + \beta^{1/2}$,
along with the fact that the integrals $\int _0 ^1 s^{-1/2} ds$ and 
$\int _0 ^1 (1-s^2)^{-1/2} ds$ are finite.  $\bullet$ \bigskip

\proclaim Corollary 3.7.  With notation as in Proposition 3.5, we suppose $D \subset D'$ 
is a disc centred on the origin of radius $a\le a'<1$, and 
let $Diam (\epsilon )$ denote the
diameter of the total space of the elliptic fibration over $D$, under an Ooguri--Vafa metric on 
$X'$ with fibre volume $\epsilon$.  There exists a constant $C_3$ (independent of both
$\epsilon$ and $a$) such that,
if $\epsilon \le a$, then
 $$Diam (\epsilon ) < C_3 a^{1/2} \epsilon ^{-1/2}.$$

\noindent \it Proof. \rm  Consider the slice $u= \epsilon /2$ of $Y$, and a radial curve 
$\gamma$ from $y=0$ to $y = a e^{i\theta }$ within this slice.  
There is a horizontal lift $\tilde \gamma$ of $\gamma$ to $X$; recalling that 
$$ds^2 = V d{\bf u}\cdot d{\bf u} + V^{-1} \theta _0 ^2,$$
we deduce that the length of $\tilde \gamma$ is just
$$ \int _{\gamma} V(y, \epsilon /2 )^{1/2} |dy| = \int _0 ^a V(re^{i\theta }, \epsilon /2 )
^{1/2} dr .$$

Since by Proposition 3.5, the diameters of the fibres are bounded above by 
$$ C_1 \epsilon ^{1/2} (\log \epsilon ^{-1} )^{1/2} < C_1 a^{1/2} \epsilon ^{-1/2},$$
if we can show that the latter integral is bounded above by 
$C a^{1/2} \epsilon ^{-1/2}$, for some constant $C$ independent of both
$\epsilon$ and $a$, then the desired 
bound for $Diam (\epsilon )$ will follow (to go between any two fibres, we can always take 
the route via the central fibre).

We estimate the above integral in two parts, from $0$ to $\epsilon$, and from $\epsilon$ to $a$.
We can estimate the first of these integrals most easily by performing the useful transformation 
described in Remark 3.4.  Recall that 
$$ \epsilon V = \tilde V_0  +  f - {1\over 2\pi }\log \epsilon .$$
Now $\tilde V_0 ( |v|, 1/2 )$ is bounded above for $0\le |v| \le 1 $ by 
$\tilde V_0 ( 0, 1/2 )$, and so 
$ \epsilon V( r e^{i\theta } , \epsilon /2 ) \le A' - {1\over 2\pi }\log \epsilon$ 
for $0\le r \le \epsilon$, where $A'$ is some positive constant.  Thus 
$$ \int _0 ^\epsilon V(re^{i\theta }, \epsilon /2 )^{1/2} dr \le \epsilon ^{-1/2} 
 \int _0 ^\epsilon (A'  - {1\over 2\pi }\log \epsilon )^{1/2} dr \le C' 
\epsilon ^{1/2} (\log \epsilon ^{-1} )^{1/2},$$ for some positive constant $C'$ independent of 
$\epsilon$ (and of course $a$).

We therefore now need to demonstrate that 
 $$ \int _\epsilon ^a V(re^{i\theta }, \epsilon /2 )^{1/2} dr$$ 
has a bound of the desired type.
To do this, we use the expression for
$V_0$ given in Lemma 3.1(b).  From the proof of  Lemma 3.1(c), we deduce that, for $|y| \ge \epsilon
/\pi$, we have 
$$ 2 \pi \epsilon V_0 < - \log |y| + 
2C'_1  {e^{-2\pi|y|/\epsilon}\over 
1-e^{-2\pi|y|/\epsilon}} .$$
In particular, since the second term is decreasing in the range, we have, 
for $|y| \ge \epsilon$, that
$$ 2 \pi \epsilon V_0 < -\log |y| + 2C'_1 {e^{-2\pi}\over 
1-e^{-2\pi}},$$
and hence that 
$$ 2 \pi \epsilon V < -\log |y| + C'_2, $$ for some constant $C'_2$ independent of 
$\epsilon$ and $a$.  Using the assumption that $a \le a' <1$, we have 
$$\eqalign{ \int _\epsilon ^a V(re^{i\theta }, \epsilon /2 )^{1/2} dr & < 
(2\pi \epsilon )^{-1/2} \int _\epsilon ^a (C'_2 -  \log r )^{1/2} dr \cr
& < \epsilon ^{-1/2} C'_3 \int _\epsilon ^a (\log {r^{-1}})^{1/2} dr \cr
& < \epsilon ^{-1/2} C'_3 \int _\epsilon ^a r^{-1/2} dr \cr
& < 2 C'_3 \epsilon ^{-1/2} a^{1/2},\cr}$$
for an appropriate constant $C'_3$, depending on $a'$ but independent of $\epsilon$ and $a$.  
The result then follows immediately. $\bullet$\bigskip

\proclaim Proposition 3.8.  With notation as in Proposition 3.5, let $R(\epsilon )$ denote the
curvature  tensor of the total space $X'$ of the elliptic fibration over $D'$, under an Ooguri--Vafa
metric  on $X'$ with
fibre volume $\epsilon$.  Then there exists positive constants $C_4 , C'_4$ (independent of $\epsilon$) such
that, for all sufficiently small $\epsilon$, 
$$ C'_4 \epsilon ^{-1} \log (\epsilon ^{-1})^{-2} < \| R(\epsilon )\| _{C^0} <  
 C_4 \epsilon ^{-1} \log (\epsilon ^{-1}),$$  where $\|\ .\ \| _{C^0}$ denotes the usual $C^0$-norm
on $X'$.

\noindent \it Proof. \rm  Recall first from Remark 2.4 that 
$$ \| R\| ^2 = 12 V^{-6} |\nabla V |^4 +  V^{-4} \Delta (|\nabla V |^2) 
- 6 V^{-5} (\nabla V)\cdot (\nabla (|\nabla V |^2)). $$

We now perform our change of coordinates $s= u/\epsilon$, $v= y/\epsilon$.  We recall that 
$V = V_0 + f(y_1 ,y_2 )/\epsilon $ for some bounded harmonic function $f$ defined over $D'$, 
and that 
$$ \epsilon V = \tilde V_0  + f - {1\over 2\pi }\log \epsilon .$$ 
Also observe that $\nabla _{u, y_1 ,y_2} = \epsilon ^{-1} \nabla _{s, v_1 ,v_2}$;  from 
now on $\nabla$ will denote $\nabla _{s, v_1 ,v_2}$, and $\Delta$ will denote 
$\Delta _{s, v_1 ,v_2}$.  We set $V_1 = \epsilon V = \tilde V_0  +  f - 
{1\over 2\pi }\log \epsilon $, considered as a function of $s, v_1 ,v_2$.  Thus 
$$\epsilon ^2 \| R \| ^2  =  12 V_1 ^{-6} |\nabla V_1 |^4 +  V_1 ^{-4} \Delta (|\nabla V_1 |^2) 
- 6 V_1 ^{-5} (\nabla V_1 )\cdot (\nabla (|\nabla V_1 |^2)).$$  

We first prove the upper bound for $\| R(\epsilon )\| _{C^0}$, namely that 
$$\| R(\epsilon )\|  < C_4 \epsilon ^{-1} \log (\epsilon ^{-1})$$ at all points of $X'$.  
The easy part of this is to deal with the points in the range $ 1/2 \le |v| < a/\epsilon $ 
(where $a$ now denotes the radius of $D'$), corresponding to $|y| \ge \epsilon /2$ 
in the disc $D'$.  Here we use the Fourier expansion for $V_1$, namely 
$$ V_1 = -{1\over 2\pi }\log |v| + f -{1\over 2\pi }\log \epsilon +
\sum_{m=-\infty\atop m\not=0}^{m=
\infty} {1\over 2\pi }e^{2\pi ims} K_0(2\pi|mv|).$$ 
Recalling that $K_0 (x)$ and its derivatives decay at least as fast as $e^{-x}$ for large $x$, it is clear that 
$|\nabla V_1 |^4$, $\Delta (|\nabla V_1 |^2)$ and $(\nabla V_1 )\cdot (\nabla (|\nabla V_1 |^2))$ 
are bounded (independent of $\epsilon$) for $1/2 \le |v| < a/\epsilon$.  Moreover,  
for $\epsilon$ sufficiently small,  
$$ \epsilon V = 
-{1\over 2\pi}\log |y| + f(y) +\sum_{m=-\infty\atop m\not=0}^{m=
\infty} {1\over 2\pi}e^{2\pi imu/\epsilon} K_0(2\pi|my|/\epsilon)$$ is bounded below, over $D'$, by some
positive constant (independent of $\epsilon$).  Thus,
$V_1$  is bounded below on $1/2 \le |v| < a/\epsilon$, 
and hence $\epsilon \| R \|$ is bounded
above on the given range by some constant, again independent of $\epsilon$.

The trickier argument is of course for the range $0\le |v| \le 1/2$, corresponding to $0\le |y| \le 
\epsilon /2$.  We assume that $\epsilon$ is small enough that $3\epsilon /4  \le a$.
 We make our usual change of variables, so that 
$$ V_1 = \epsilon V = \tilde V_0  + f + {1\over 2\pi }\log (\epsilon ^{-1})$$ defines an Ooguri--Vafa metric 
over the disc $ |v| < 3/4 $, fibres of volume one, and periods $\{1, (2\pi i)^{-1}\log v +  
(2\pi i)^{-1}\log \epsilon + ih \}$.  Now choose $A \ge 0$ such that $f + A >0$ whenever $|v| < 3/4$, and set
$V_2 = \tilde V_0  + f +A$; $V_2$ then determines an Ooguri--Vafa metric 
over the disc $ |v| < 3/4 $, fibres of volume one, and periods $\{ 1, (2\pi i)^{-1}\log v +  ih + iA \}$.
We may obviously assume that $A < (2\pi )^{-1}\log (\epsilon ^{-1})$. 
Let $R_1$,  respectively $R_2$, denote the curvature tensors of the metrics determined by $V_1$, 
respectively $V_2$.  Our aim now is to show that $\| R_1 \|^2 < C\, \log (\epsilon ^{-1})^2$ over 
the disc $ |v| \le 1/2 $; 
if this is true, it follows from the above that 
$\| R \| < C_4 \epsilon ^{-1} \log (\epsilon ^{-1}) $ at all points over $D'$, for some positive
constant 
$C_4$.

Since the metric determined by $V_2$ is independent of $\epsilon$, it is clear that $\| R_2 \|$ is bounded 
over $ |v| \le 1/2 $, with the bound independent of $\epsilon$.  Hence 
$$ 12 V_2 ^{-6} |\nabla V_2 |^4 +  V_2 ^{-4} \Delta (|\nabla V_2 |^2) 
- 6 V_2 ^{-5} (\nabla V_2 )\cdot (\nabla (|\nabla V_2 |^2)) \leqno{(3.1)}$$ is bounded independent of
$\epsilon$ over the  range in question, $|v| \le 1/2$, which from now on will be taken as
understood.  We wish to show that 
$$  \| R_1 \| ^2 
= 12 V_1 ^{-6} |\nabla V_2 |^4 +  V_1 ^{-4} \Delta (|\nabla V_2 |^2) 
- 6 V_1 ^{-5} (\nabla V_2 )\cdot (\nabla (|\nabla V_2 |^2)) \le C (\log (\epsilon ^{-1}))^2. $$
Since $V_1 = V_2 + (2\pi )^{-1}\log (\epsilon ^{-1}) - A \ge V_2$, it will be enough to prove the same bound for 
$$ 12 V_1 ^{-2} V_2 ^{-4} |\nabla V_2 |^4 +  V_2 ^{-4} \Delta (|\nabla V_2 |^2) 
- 6 V_1 ^{-1} V_2 ^{-4}(\nabla V_2 )\cdot (\nabla (|\nabla V_2 |^2)). \leqno{(3.2)}$$ 
By subtracting our previously bounded expression (3.1), we need then only show boundedness for  
$$ 12 V_1 ^{-2} V_2 ^{-6} (V_2 ^2 - V_1 ^2)|\nabla V_2 |^4 
- 6 V_1 ^{-1} V_2 ^{-5}(V_2 - V_1) (\nabla V_2 )\cdot (\nabla (|\nabla V_2 |^2)).\leqno{(3.3)}$$
Expanding this latter expression out, we get
$$ 6 \left( (2\pi )^{-1} \log (\epsilon ^{-1}) - A \right) V_1 ^{-1} \left( V_2 ^{-5}(\nabla V_2
)\cdot (\nabla (|\nabla V_2 |^2))  - 2 V_2 ^{-6} (1+ {V_2 \over V_1}) |\nabla V_2 |^4 \right).
\leqno{(3.4)}$$

We now claim that $$ \left| V_1 ^{-1} \left( V_2 ^{-5}(\nabla V_2 )\cdot (\nabla (|\nabla V_2
|^2))  - 4 V_2 ^{-6} |\nabla V_2 |^4 \right)\right| $$
and $$ V_1 ^{-2} V_2 ^{-6} |\nabla V_2 |^4$$ are bounded independent of $\epsilon$.  If this 
is true, then the latter bound will
imply that $$V_1 ^{-1} (1 - {V_2 \over V_1}) V_2 ^{-6} |\nabla V_2 |^4 \le 
C' (\log (\epsilon ^{-1}) - 2\pi A)$$ for 
some positive $C'$, and then the former bound implies that the expression we are interested in has a 
bound of the form $ B_1 ( \log (\epsilon ^{-1}) - 2\pi A)^2 + B_2 ( \log (\epsilon ^{-1}) - 2\pi A)$, for suitable
positive  constants $B_1 , B_2$.  This then gives the required result.

To show boundedness for the two remaining quantities, it is sufficient to bound the functions  
$$ V_2 ^{-6} \left| 4 V_2 ^{-1} |\nabla V_2 |^4 - (\nabla V_2 )\cdot (\nabla (|\nabla V_2 |^2))\right| $$ and 
$$V_2 ^{-8} |\nabla V_2 |^4.$$  Both these functions are defined away from $\{ 0 \} \times \Z$ and are
 periodic in
$s$; moreover, they plainly do not depend on $\epsilon$.  If we show that they are in fact both
continuous at  the origin ($v=0, s =0$), the existence of the required bounds will follow
automatically.

We now write $4\pi V_2 = \rho ^{-1} + w$, where $\rho = (s^2 + v_1 ^2 + v_2 ^2)^{1/2}$ and $w$ is a
harmonic  function on a neighbourhood of the origin.  Then we see that   
$ (4\pi )^4 | \nabla V_2 |^4 = \rho ^{-8} + O(\rho ^{-6}) $.  Since $ (4\pi V_2 ) ^{-8} = \rho ^8 (1
+ w\rho )^{-8}$, we deduce that  
$  V_2 ^{-8} |\nabla V_2 |^4$ is regular at the origin, taking the value $(4\pi )^4$ there.
Moreover, it is easily checked that $$(4\pi )^3 (\nabla V_2 )\cdot (\nabla (|\nabla V_2 |^2)) = 4 \rho
^{-7} +  O(\rho ^{-5}) ,$$ and so in particular 
$$  4 V_2 ^{-7} |\nabla V_2 |^4 -  V_2 ^{-6}(\nabla V_2 )\cdot (\nabla (|\nabla V_2 |^2) $$ is also 
regular at the origin, and vanishes there. 

We now turn to the lower bound for $\| R(\epsilon )\| _{C^0}$.  We work on the transformed elliptic 
fibration over the disc $| v | \le 1/2$, and let $M$ denote the $C^0$-norm of the function given 
by (3.1).

\rm From the above calculations, at all points $P$ with sufficiently small value of $\rho$, we have 
$$ 
\left| V_2 ^{-6} (\nabla V_2 )\cdot (\nabla (|\nabla V_2 |^2)
- 4 V_2 ^{-7} |\nabla V_2 |^4 \right|  < M \quad , \quad
V_2 ^{-7} |\nabla V_2 |^4 > 2M . $$  We now fix such a point $P$; 
note that the coordinates $s,v_1 ,v_2$ 
are then taken to be fixed, and so this does not correspond to taking a 
fixed point (independent of $\epsilon$) on our original family $X'$.

Observe now that 
$$ V_1 /V_2 = 1 + {(2\pi )^{-1} \log (\epsilon ^{-1}) - A \over V_2};$$
so for $P$ fixed, $V_1 (P)/V_2 (P) > 2$ for $\epsilon$ sufficiently small.  From this it 
follows that, when evaluated at $P$,
$$ \left|  V_2 ^{-6}(\nabla V_2
)\cdot (\nabla (|\nabla V_2 |^2))  - 2 V_2 ^{-7} (1+ {V_2 \over V_1}) |\nabla V_2 |^4 \right| 
 > M ,$$ for $\epsilon$ sufficiently small.  Hence, for $\epsilon$ sufficiently small, the 
modulus of (3.4) evaluated at $P$ is at least $3M$ say, and thus the same is 
true of (3.3).  From this, and our original choice for $M$, it follows that the modulus of (3.2) evaluated at $P$ 
is at least $2M$.  Therefore, for $\epsilon$ sufficiently small, 
$$ \| R_1 (P)\|^2 > B (\log (\epsilon ^{-1}))^{-4}$$
for some constant $B$ independent of $\epsilon$.  Thus 
$$ \| R(\epsilon ) \|_{C^0} > B^{1/2} \epsilon ^{-1} (\log (\epsilon ^{-1}))^{-2},$$
as required.   $\bullet$
\bigskip

{\hd \S 4. Almost Ricci-flat metrics on Elliptic K3 Surfaces.}

Our goal in this section is to construct K\"ahler metrics on elliptic
K3 surfaces which are very close to being Ricci-flat by gluing
the Ooguri--Vafa metric in neighbourhoods of singular fibres
to the semi-flat metric away from the singular fibres.

We begin by producing one such metric on a Jacobian elliptic fibration.
Fix a K3 surface $X$ with a fixed holomorphic 2-form $\Omega$ and
an elliptic fibration $f:X\rightarrow B=\Pone$, which we will take
to have a holomorphic section $\sigma_0$. Furthermore, assume all
singular fibres of $f$ are of Kodaira type $I_1$; there will then
be 24 such fibres. Let $p_1,\ldots,p_{24}\in B$ be those
points for which $X_{p_i}=f^{-1}(p_i)$ is singular, 
$\Delta=\{p_1,\ldots,p_{24}\}$, $B_0=B-\Delta$, $X_0=f^{-1}(B_0)$, and
$X^{\#}=X-Sing(f^{-1}(\Delta))$. There is an exact sequence, already mentioned
in \S 1,
$$0\mapright{} R^1f_*\boldz\mapright{}\T_B^*\mapright{\phi} X^{\#}\mapright{} 0,$$
with the property that $\phi$ maps the zero section of $\T_B^*$ to
$\sigma_0$ and $\phi^*\Omega$ is the canonical holomorphic 2-form on
$\T_B^*$, which is $dx\wedge dy$ if $y$ is a coordinate on $B$ and
$x$ a canonical fibre coordinate. (See [14], Proposition 7.2). Here
$x=0$ defines the zero section.

Given this data, by Example 2.2, for each $\epsilon$, there exists 
a well-defined Ricci-flat metric on $X_0$, the standard semi-flat metric
$\omega_{SF}$, with the area of each fibre being $\epsilon$. The
reader should keep in mind the dependence of $\omega_{SF}$ on $\epsilon$.

Now let $y$ be a holomorphic coordinate on $B$ defined in a neighbourhood
$U$ of $p\in\Delta$, $U$ contractible with $U\cap\Delta=\{p\}$, and
$y=0$ at the point $p$. Let $x$ be the corresponding canonical fibre
coordinate. Let $U^*=U-\{p\}$, $X_{U^*}=f^{-1}(U^*)$. We can then
choose over $U^*$ holomorphic periods $\tau_1(y),\tau_2(y)$, representing
possibly
multi-valued holomorphic sections of $\T^*_{U^*}$ generating the period lattice.
Because the monodromy about an $I_1$ fibre in a suitable basis is
$\pmatrix{1&1\cr 0&1\cr}$, we can take
one of these, say $\tau_1$, to be 
single-valued, though $\tau_2$ will be multi-valued around the $I_1$ fibre.
We will always choose $\tau_1$ and $\tau_2$ so that $\im(\bar\tau_1\tau_2)>0$.
Set
$$\eqalign{
W_0(y)&=1/\Im(\bar\tau_1\tau_2)\cr
b_0(x,y)&=
-{\Im(\tau_2\bar x)\partial_y\tau_1+\Im(\bar\tau_1 x)
\partial_y\tau_2\over \Im(\bar\tau_1\tau_2)}\cr}$$
and
$$\eqalign{\pv&=W_0^{-1}\px\cr
\ph&=\py-b_0\px\cr
\tv&=W_0(dx+b_0dy)\cr
\th&=dy\cr}$$
as in \S 2. The latter two 1-forms are well-defined on $X_U^*$,
so form a basis for $(1,0)$ forms. We denote by $\pvc$ et cetera
the complex conjugates of the above as usual.

\proclaim Lemma 4.1. Let $\omega$ be a real closed $(1,1)$ form
on $X_{U^*}$, with
$$\omega={i\over 2}(\alpha\tv\wedge\tvc+\beta \th\wedge\tvc
+\bar\beta\tv\wedge\thc+\gamma \th\wedge\thc).$$
There exists a function $\varphi$ on $X_{U^*}$ such that $\omega
=i\partial\bar\partial\varphi$ if and only if $\omega$ represents
the zero cohomology class on $X_{U^*}$ and 
$$\int_{X_b}\beta dx_1\wedge dx_2=0$$
for all $b\in U^*$. Furthermore, for $0<r_1<r_2$,
let $U_{r_1,r_2}=\{y\in U\,|\, r_1< |y| < r_2\}$. If $r_1<r_1'<r_2'<r_2$ and
$\overline{U}_{r_1,r_2}\subseteq U^*$,
then there exists a constant
$C$ depending only on $r_1,r_2,r_1',r_2'$ and the periods
of $f$ over $U_{r_1,r_2}$ such that $\varphi$ can be chosen with
$$\|\varphi\|'_{C^{k+2,\alpha}}\le C(\|\alpha\|_{C^{k,\alpha}}
+\|\beta\|_{C^{k,\alpha}}+\|\gamma\|_{C^{k,\alpha}}).$$
Here, we compute the $C^{k,\alpha}$ norm of a function on 
$f^{-1}(U_{r_1,r_2})$ by thinking of them as functions on $\T^*_{U_{r_1,r_2}}$,
which we embed in ${\bf C}^2$ by the coordinates $x$ and $y$. We can then
use the standard $C^{k,\alpha}$ norms on a bounded open set of 
$\T_{U_{r_1,r_2}}^*$ which contains a fundamental domain of each fibre.
The norm $\|\cdot\|'_{C^{k,\alpha}}$ denotes the similar norm of
a function over $U_{r_1',r_2'}$.

\noindent {\bf Remark 4.2.} We note that the definition of the $C^{k,\alpha}$
norm given above depends on the choice of holomorphic coordinate $y$
and the bounded open set, but any two such norms will be equivalent.

\noindent {\it Proof.} Before beginning the proof, we observe from (2.4) and (2.3) that 
$$\eqalign{
\partial\tvc=&((\pv W_0)\tv+(\ph W_0)\th)\wedge W_0^{-1}\tvc\cr
&+W_0((\pv\bar b_0)\tv+(\ph\bar b_0)\th)\wedge\thc\cr
\bar\partial \tvc=&0.\cr}\leqno{(4.1)}$$

Furthermore, locally for the base, a function on $X_U^*$ can be expanded
in a Fourier series on the fibres, yielding
$$f(x,y)=\sum_{n,m\in\boldz} a_{n,m}(y)e^{2\pi i(n\Im(\tau_2\bar x)
+m\Im(\bar\tau_1 x))/\Im(\bar\tau_1\tau_2)}.$$
A direct calculation shows that
$$\ph f=
\sum_{n,m\in\boldz}\py(a_{n,m})e^{2\pi i(n\Im(\tau_2\bar x)
+m\Im(\bar\tau_1 x))/\Im(\bar\tau_1\tau_2)}.$$

If a $\varphi$ exists, then of course $\omega$ represents
the zero cohomology class on $X_{U^*}$. Also,
$$i\partial\bar\partial\varphi=i\partial((\pvc\varphi)\tvc+
(\phc\varphi)\thc).$$
 From (4.1) it then follows that if $\omega=i\partial\bar\partial\varphi$, then
$$\beta=2(\ph\pvc\varphi+(\pvc\varphi)W_0^{-1}\ph W_0),$$
and then by looking at the constant term $a_{0,0}$ of the Fourier
expansion of $\beta$, it is clear that $\int_{X_b} \beta dx_1\wedge dx_2=0$.

Conversely, first suppose $\omega$ is cohomologically trivial.
Then there exists
a one-form $\xi$ of type $(1,0)$ such that ${i\over 2}d(\bar\xi-\xi)=\omega$
(since $\omega$ is real). Necessarily $\partial\xi=\bar\partial\bar\xi
=0$. Thus $\bar\xi$ represents a class in $H^{0,1}(X_{U^*})$. If this class is 
zero, then there exists a function $\varphi$ such that $\bar\partial\varphi
=\bar\xi$, and then $\partial\bar\varphi=\xi$, so
$$\omega={i\over 2}(\partial\bar\partial\varphi-\bar\partial\partial\bar\varphi)
={i\over 2}\partial\bar\partial(\varphi+\bar\varphi)
=i\partial\bar\partial\re\varphi,$$
as desired. Thus we need to understand when $\bar\xi$ represents the zero class.

Now $H^{0,1}(X_{U^*})=H^1(X_{U^*},\O_{X_{U^*}})$, which by the Leray spectral
sequence for $f$ is isomorphic to $H^0(U^*,R^1f_*\O_{X_{U^*}})$, as
$H^i(U^*,f_*\O_{X_{U^*}})=H^i(U^*,\O_{U^*})=0$ for $i\ge 1$. Thus
$\bar\xi$ represents zero in $H^{0,1}(X_{U^*})$ if and only if
$\bar\xi|_{X_b}$ represents the zero class in $H^{0,1}(X_b)$ for all
$b\in U^*$. If we write $\bar\xi=\bar g\tvc+\bar h\thc$,
this is equivalent to the 
constant term in the Fourier expansion of $\bar g$ on the fibre being zero.
Denote this constant term by $\bar g_0(y)$.

What kind of function is $\bar g_0$? Well, by (4.1),
$$0=\bar\partial\bar\xi=(\phc\bar g-\pvc\bar h)\thc\wedge\tvc.$$
By looking at the constant term of the Fourier expansion of this
coefficient, we see $\pyc\bar g_0=0$, so $\bar g_0$ is a
holomorphic function on $U^*$. This function
gives the section of $R^1f_*\O_{X_{U^*}}$ defined by $\bar\xi$.

Now let us compute the coefficient $\beta$ of $\th\wedge\tvc$
in $\omega$ in terms of $g$ and $h$. 
 From $\omega={i\over 2}(\partial\bar\xi-\bar
\partial\xi)$ and (4.1), it follows that
$$\beta=\ph \bar g+\pvc h+\bar g W_0^{-1}\ph W_0+gW_0\pvc b_0.$$
If $\beta_0$ is the constant term in the Fourier expansion of $\beta$,
then we get, using $(2.2')$ for the second line,
$$\eqalign{\beta_0&=\py\bar g_0+\bar g_0W_0^{-1}\ph W_0+g_0 W_0
\pvc b_0\cr
&=\py\bar g_0+\bar g_0\px b_0+g_0\pxc b_0\cr
&=\py\bar g_0+{\bar g_0(\bar\tau_2\py\tau_1-\bar\tau_1\py\tau_2)
-g_0(\tau_2\py\tau_1-\tau_1\py\tau_2)\over \bar\tau_1\tau_2-\tau_1\bar\tau_2}
\cr
&=\py\bar g_0+b_0(\bar g_0,y).\cr}\leqno{(4.2)}$$
If we now assume in addition that $\int_{X_b}\beta dx_1\wedge dx_2=0$
for all $b\in U^*$, then $\beta_0=0$, so
$$\py \bar g_0+b_0(\bar g_0,y)=0.$$
Now write $\bar g_0(y)=a_1(y)\tau_1(y)+a_2(y)\tau_2(y),$
where $a_1,a_2$ are real functions of $y$. Then $b_0(\bar g_0,y)
=-a_1\py\tau_1
-a_2\py\tau_2$, so 
$$0=\py \bar g_0+b_0(\bar g_0,y)=(\py a_1)\tau_1+
(\py a_2)\tau_2.$$
But 
$$0=\pyc \bar g_0=(\pyc a_1)\tau_1+(\pyc a_2)\tau_2.$$
Thus combining these two equations gives
$$\eqalign{(\partial_{y_1}a_1)\tau_1+(\partial_{y_1}a_2)\tau_2&=0\cr
(\partial_{y_2}a_1)\tau_1+(\partial_{y_2}a_2)\tau_2&=0,\cr}$$
and by linear independence of $\tau_1$ and $\tau_2$ we see $a_1$ and
$a_2$ are constant. Since $\bar g_0$ is well-defined and we are assuming
$\tau_1$ is the monodromy invariant period, we have $\bar g_0=a\tau_1$, 
$a$ a constant.
Now a calculation shows that
$${i\over 2}d(\tau_1\tvc-\bar\tau_1\tv)=0,$$
so we can subtract $a\tau_1\tvc$
from $\bar\xi$ without affecting ${i\over 2}d(\bar\xi-\xi)=\omega$. Thus
we can assume $\bar g_0=0$, 
and then $\bar\xi$ represents the zero class in $H^{0,1}
(X_{U^*})$, allowing us to complete the proof of the existence of
$\varphi$. 

Now we need to control the norm of $\varphi$. First note that
$$W_0^2\alpha=2\partial_x\partial_{\bar x}\varphi={1\over 2}\Delta_x\varphi,$$
where $\Delta_x=\partial_{x_1}^2+\partial_{x_2}^2$ denotes the standard
Laplacian on fibres. Writing $\varphi=\varphi_0+\varphi_v$ where
$\varphi_0$ is the pull-back of a function on $U^*$ and 
$\int_{X_b}\varphi_v dx_1 dx_2=0$ for all $b\in U^*$, we have
$W_0^2\alpha=\Delta_x\varphi_v/2$. It then follows that $|\varphi_v|$ is bounded
on each fibre (being a torus) with the bound proportional to a bound
for $|\alpha|$ on that fibre, with the constant of proportionality
depending on the periods at that point. (To see this, one can just
work with Fourier series). Thus $$\|\varphi_v\|_{C^0}\le C_1\|\alpha\|_{C^0}$$
on $U_{r_1,r_2}$, where $C_1$ depends on the periods over $U_{r_1,r_2}$.

Next restrict $\varphi$ and $\omega$ to the zero section of 
$f:X_U\rightarrow U$. On this zero-section, 
$$\eqalign{\omega&={i\over 2}\gamma dy\wedge d\bar y\cr
&=i\partial\bar\partial\varphi\cr
&=i\partial_y\partial_{\bar y}\varphi dy\wedge d\bar y\cr
&={i\over 4}\Delta_y\varphi dy\wedge d\bar y.\cr}$$
so $\Delta_y\varphi=2\gamma$ on the zero section, where $\Delta_y$
is the the standard Laplacian $\partial^2_{y_1}+\partial^2_{y_2}$
on $U^*$.

Now let $\psi$ be a harmonic function on $U_{r_1,r_2}$ such that
$\psi|_{\partial\overline{U}_{r_1,r_2}}=\varphi|_{\partial\overline{
U}_{r_1,r_2}}$. (Here we are identifying $U_{r_1,r_2}$ with its image
under the zero section.) This function exists and is unique.
Then
$$\Delta_y(\varphi-\psi+{\|2\gamma\|_{C^0}\over 4}(y_1^2+y_2^2))
=2\gamma+\|2\gamma\|_{C^0}\ge 0.$$ Thus by the maximum principal,
$\varphi-\psi+\|2\gamma\|_{C^0}(y_1^2+y_2^2)/4$ achieves its maximum
when either $|y|=r_1$ or $|y|=r_2$, and since $\varphi-\psi=0$ on the boundary
of $U_{r_1,r_2}$, we have
$$\varphi-\psi\le \|2\gamma\|_{C^0}r_2^2/4.$$
Similarly, $\psi-\varphi\le \|2\gamma\|_{C^0}r_2^2/4$, so
$$\|\varphi-\psi\|_{C^0}\le \|\gamma\|_{C^0}r_2^2.$$
This estimate holds on $U_{r_1,r_2}$, but from $\|\varphi_v\|_{C^0}
\le C_1\|\alpha\|_{C^0}$, it is clear that the oscillation of $\varphi$
along the fibres is bounded by $C_1\|\alpha\|_{C^0}$, and thus
on $f^{-1}(U_{r_1,r_2})$,
$$\|\varphi-\psi\|_{C^0}\le C_2(\|\alpha\|_{C^0}+\|\gamma\|_{C^0})$$
for some constant $C_2$ depending on the periods, $r_1$, and $r_2$.
Noting that $\partial\bar\partial\psi=0$, we can replace $\varphi$
by $\varphi-\psi$. Then the $C^{k+2,\alpha}$ estimates follow from
the standard interior Schauder estimates for the Laplacian 
(see [11], problem 6.1.) This is because the ordinary Laplacian
(in the coordinates $x_1,x_2,y_1,y_2$) of $\varphi$ can be expressed
in terms of $\alpha,\beta$ and $\gamma$.
$\bullet$

\proclaim Lemma 4.3. Let $\omega$ be a K\"ahler form on $X_U$, $\omega_{SF}$
the semi-flat K\"ahler form on $X_0$, such that 
$$\int_{X_b}\omega=\int_{X_b}\omega_{SF}=\epsilon.$$
Then $[\omega_{SF}-\omega]=0$ in $H^2(X_{U^*},{\bf R})$, and
furthermore, there exists a holomorphic section $\sigma$ of
$f:X_U\rightarrow U$ and a function $\varphi$ on $X_{U^*}$ such that
$$\omega_{SF}-T_{\sigma}^*\omega=i\partial\bar\partial\varphi,$$
where $T_{\sigma}$ is translation by the section $\sigma$.

\noindent {\it Proof.} To show the first part, we first observe that $H_2(X_{U^*},\boldz)$
is generated by the homology classes of two submanifolds: $X_b$ for
some $b\in U^*$, and $T$, where $T$ is a torus fibred in circles over
a simple closed loop $\gamma:[0,1]\rightarrow U^*$ 
generating $\pi_1(U^*)$, with the class of the
fibre being the monodromy invariant cycle. To show $[\omega_{SF}-\omega]=0$,
we just need $\int_{X_b} \omega_{SF}-\omega=0$, which is obvious, and
$\int_T \omega_{SF}-\omega=0$. Now $\int_T\omega=0$ since $\omega$ is defined
on $X_U$, where $T$ is homologous to zero. On the other hand,
if we describe $T$ explicitly, parametrised by coordinates $s,t$
with $\mu:[0,1]^2\rightarrow X_{U^*}$ given by
$$\mu(s,t)=(x(s,t),y(s,t))=(s\tau_1(\gamma(t)),\gamma(t)),$$
then a calculation
shows that $\mu^*\omega_{SF}=0$, and hence $\int_T\omega_{SF}=0$.
Thus $[\omega_{SF}-\omega]=0$.

As in Lemma 4.1, 
write, for each section $\sigma$ of $f:X_{U^*}\rightarrow U^*$,
$$\omega_{SF}-T_{\sigma}^*\omega={i\over 2}
(\alpha_{\sigma}\tv\wedge\tvc+
\beta_{\sigma} \th\wedge
\tvc+\cdots).$$
Let $\sigma_0$ be the zero section, so that $T_{\sigma_0}$ is the identity.
We showed in (4.2) that the function $\beta_0$, the constant term in the
Fourier expansion of $\beta_{\sigma_0}$,
was of the form
$$\beta_0=\py k+b_0(k,y)$$
where $k(y)$ is a holomorphic function on $U^*$.

Now write
$$\omega={i\over 2}(WW_0^{-2}\tv\wedge\tvc+\beta_{\omega}\th\wedge\tvc
+\cdots)$$
where necessarily the constant term of $W$ is
$$(\Im\bar\tau_1\tau_2)^{-1}\int_{X_b} W dx_1\wedge dx_2=
(\Im\bar\tau_1\tau_2)^{-1}\int_{X_b}{i\over 2}WW_0^{-2}\tv\wedge\tvc=
\epsilon/\Im(\bar\tau_1\tau_2).$$
We calculate $T_{\sigma}^*\omega$. First note that
$$\eqalign{T^*_{\sigma}(\tv)&=
W_0(d(x+\sigma(y))+b_0(x+\sigma(y),y) dy)\cr
&=W_0(dx+b_0 dy)+W_0(\py\sigma(y)+b_0(\sigma(y),y))dy\cr
&=\tv+W_0(\py\sigma+b_0(\sigma(y),y))\th.\cr}$$
Thus the coefficient of $\th\wedge\tvc$ in $T_{\sigma}^*(\omega)$ is
$${i\over 2}(\beta_{\omega}\circ T_{\sigma}+{W\over W_0}(\py\sigma+b_0(\sigma(y),
y))).$$
On the other hand, $\omega_{SF}={i\over 2}W_0^{-1}(\epsilon\tv\wedge\tvc
+\epsilon^{-1}\th\wedge\thc)$. Thus $\beta_{\sigma_0}=-\beta_{\omega}$,
and $$\beta_{\sigma}=\beta_{\sigma_0}\circ T_{\sigma}-{W\over W_0}(\py\sigma
+b_0(\sigma(y),y)).$$
So
$$\eqalign{
(\Im\bar\tau_1\tau_2)^{-1}
\int_{X_b}\beta_{\sigma} dx_1\wedge dx_2
=&(\Im\bar\tau_1\tau_2)^{-1}\bigg(
\int_{X_b}\beta_{\sigma_0}\circ T_{\sigma} dx_1\wedge dx_2\cr
&-W_0^{-1}(\py \sigma+b_0(\sigma(y),y))\int_{X_b} W dx_1\wedge dx_2\bigg)\cr
=&\beta_0
-\epsilon(\py \sigma+b_0(\sigma(y),y)).\cr}$$
If we take $\sigma(y)=k(y)/\epsilon$, this will yield zero.
So for this choice of $\sigma$,
$\omega_{SF}-T_{\sigma}^*\omega=i\partial\bar\partial\varphi$ for
some function $\varphi$ on $X_{U^*}$.

Note that a holomorphic section of $f$ over $U^*$ always extends to a
holomorphic section of $f$ on $U$. $\bullet$

\proclaim Theorem 4.4. Let $f:X\rightarrow\Pone$ be an elliptically
fibred K3 surface with a holomorphic
section and 24 singular fibres over $\Delta=\{p_1,\ldots,
p_{24}\}$ as above. Then there exists open sets $U_1^i\subseteq U_2^i\subseteq
\Pone$,
$i=1,\ldots,24$, each diffeomorphic to a disc, $U^i_j\cap \Delta=\{p_i\}$,
positive constants $D_1,\ldots,D_6$ and $\epsilon_0$ such that, for all
$\epsilon<\epsilon_0$, there exists a K\"ahler metric $\omega_{\epsilon}$
on $X$ with the following properties:
\item{(1)} $$\int_X\omega_{\epsilon}^2=\int_X(\re\Omega)^2=\int_X(\im\Omega)^2$$
\item{(2)} $$\int_{X_b}\omega_{\epsilon}=\epsilon$$
\item{(3)} $$\omega_{\epsilon}|_{f^{-1}(\Pone\setminus
\bigcup_i U_2^i)}=\omega_{SF}$$
\item{(4)} $\omega_{\epsilon}|_{f^{-1}(U_1^i)}=T_{\sigma_i}^*\omega_{OV},$
where $\omega_{OV}$ is an Ooguri--Vafa metric and
$T_{\sigma_i}$ denotes translation by some holomorphic section $\sigma_i$.
\item{(5)} If $F_{\epsilon}=\log\left({\Omega\wedge\bar\Omega/2\over 
\omega_{\epsilon}^2}\right)$, then
$$\|F_{\epsilon}\|_{C^0}\le D_1e^{-D_2/\epsilon}$$
and
$$\|\Delta F_{\epsilon}\|_{C^0}\le D_1e^{-D_2/\epsilon}$$
where $\Delta$ denotes the Laplacian with respect to the metric
$\omega_{\epsilon}$.
\item{(6)} $$inf_v\{Ric(v,v)\, |\ |v|_{\omega_{\epsilon}}=1\}\ge -D_3
e^{-D_4/\epsilon}.$$
\item{(7)} With the Riemannian metric induced by $\omega_{\epsilon}$,
$Diam(X)\le D_5\epsilon^{-1/2}$.
\item{(8)}  If $R$ denotes the Riemann curvature tensor, then
$$\|R\|_{C^0}\le D_6\epsilon^{-1}\log\epsilon^{-1},$$
$$\hbox{$\|R\|_{C^0}\rightarrow\infty$ as $\epsilon\rightarrow 0$,}$$
and on any non-singular fibre, there exists a constant $C$ depending on
the fibre such that 
$$\|R\|\le C\epsilon.$$

\noindent {\it Proof.} Let $p\in\Delta$; we fix our attention near this point.
Choosing a holomorphic coordinate $y$ in a neighbourhood of $p$,
we can express the holomorphic periods of $f$ as $\tau_1(y),\tau_2(y)$,
where $\tau_1$ is taken to be single valued. In $\T_B^*$, this
coincides with the holomorphic differential $\tau_1(y)dy$. Locally,
there exists a function $g(y)$ with $dg=\tau_1(y)dy$;
since $\tau_1(p)\not=0$, we can use $g$ as a local holomorphic coordinate
in a neighbourhood of $p$. Replacing $y$ by $g$, we can
then assume that $\tau_1(y)=1$ and also that $y=0$ at $p$. By
results of \S 3, we can then construct for all $\epsilon$ less than some
$\epsilon_0$, a metric $\omega_{OV}$ on
$f^{-1}(U)$, for some $U=\{y\, | \ |y|<r\}$, for some $r$ which only
depends on the period $\tau_2$ and $\epsilon_0$, but not
$\epsilon$. Fix $r_1<r_2<r$, and let
$U_i=\{y\, |\ |y|<r_i\}$. If $p=p_j$, we set $U_i^j=U_i$.

Remaining focused near $p$,
let $\psi:(0,(r_2+\delta)^2)\rightarrow [0,1]$
be a fixed $C^{\infty}$ 
cut-off function, with $\psi(r^2)=1$ for $r\le r_1$, $\psi(r^2)=0$
for $r\ge r_2$. Now apply Lemma 4.3 with $\omega=\omega_{OV}$.
Then there exists a holomorphic section $\sigma$ of $f$ over $U$, such that
$$\omega_{SF}-T_{\sigma}^*\omega_{OV}=i\partial\bar\partial\varphi$$
for some function $\varphi$ on $X_{U^*}$.
We can then glue $T_{\sigma}^*\omega_{OV}$ and $\omega_{SF}$ by
$$\omega_{new}=\omega_{SF}-i\partial\bar\partial(\psi(|y|^2)\varphi).$$
For $|y|\ge r_2$, $\omega_{new}$ coincides with $\omega_{SF}$;
for $|y|\le r_1$, $\omega_{new}$ coincides with $T_{\sigma}^*\omega_{OV}$.
This can be done at each singular fibre, obtaining a global closed
real $(1,1)$ form $\omega_{new}$.

We still need to check $\omega_{new}$ is positive. One
calculates that on $X_{U^*}$, 
$$\eqalign{\omega_{new}=&(1-\psi(|y|^2))\omega_{SF}+\psi(|y|^2)T_{\sigma}^*
\omega_{OV}\cr
&-i(\psi'(|y|^2)\bar y dy\wedge\bar\partial\varphi+\psi'(|y|^2)y
\partial\varphi\wedge d\bar y+\psi''(|y|^2)|y|^2\varphi dy\wedge d\bar y).\cr}$$
The sum of the first two terms is positive, so we need to make sure
the last three terms are small. Thus we need to control the size
of $\varphi$. To do so, we need to show $\omega_{SF}-T_{\sigma}^*\omega_{OV}$
is small. Now
$$\omega_{SF}={i\over 2}W_0^{-1}(\epsilon\tv\wedge\tvc+\epsilon^{-1}
\th\wedge\thc).$$
On the other hand, we can assume $\sigma$ is the zero section by having
chosen the right holomorphic section in Construction 2.6 to perform
the transformation between coordinates, and write, with $W=V^{-1}$,
$$\eqalign{\omega_{OV}=&{i\over 2}(W(dx+bdy)\wedge\overline{(dx+bdy)}
+W^{-1}dy\wedge d\bar y)\cr
=&{i\over 2}(W(dx+b_0dy)\wedge\overline{(dx+b_0dy)}+W(b-b_0)dy\wedge
\overline{(dx+b_0dy)}\cr&+W(\bar b-\bar b_0)(dx+b_0dy)\wedge d\bar y
+(W|b-b_0|^2+W^{-1})dy\wedge d\bar y)\cr
=&{i\over 2}(WW_0^{-2}\tv\wedge\tvc +{W\over W_0}(b-b_0)\th\wedge\tvc
\cr&+{W\over W_0}(\bar b-\bar b_0)(\tv\wedge\thc)+(W|b-b_0|^2+W^{-1})\th\wedge
\thc).\cr}$$
Thus we are applying Lemma 4.1 with $\alpha=\epsilon W_0^{-1}-WW_0^{-2}$,
$\beta={W\over W_0}(b_0-b)$, and $\gamma=\epsilon^{-1}W_0^{-1}-W^{-1}
-W|b-b_0|^2$. Now we work in Gibbons--Hawking coordinates using the
fact that $\alpha,\beta$ and $\gamma$ are invariant under the action
$x_1\mapsto x_1+t$. So if we bound the $C^k$ norm of $\alpha$, $\beta$
and $\gamma$ as functions on $U_{r_1,r_2}\times {\bf R}/\epsilon\boldz$,
with coordinates $y$ and $u$, we can apply (2.5) to bound the $C^k$ norms
of $\alpha,\beta,\gamma$ with respect to the coordinates $x$ and $y$.
The interpolation inequalities then give $C^{k',\alpha'}$ bounds for any
$k'<k$.

First look at $\alpha$. Now 
$$\eqalign{V&=V_0+{1\over 4\pi\epsilon}\log|y|^2+\epsilon^{-1}\im(\tau_2)\cr
&=\epsilon^{-1}\im(\tau_2)+g(u,y)\cr
&=\epsilon^{-1}W_0^{-1}+g(u,y),\cr}$$
where $g(u,y)$ is, by Lemma 3.1 (c), a harmonic function on
$U_{r_1,r_2}\times {\bf R}/\epsilon\boldz$ with $\|g\|_{C^0}$ being $O(e^{-C/\epsilon})$.
It then follows from [11], Theorem 2.10, that for each $k$, $\|g\|_{C^k}$
is also $O(e^{-C/\epsilon})$. Thus
$$\eqalign{\alpha&=\epsilon W_0^{-1}-WW_0^{-2}\cr
&={\epsilon\over W_0}-{\epsilon\over W_0+\epsilon g W_0^2}\cr
&={\epsilon^2 g W_0^2\over W_0(W_0+\epsilon g W_0^2)}.\cr}$$
Now using the fact that the denominator is bounded above and below, and
observing that any derivative of $\alpha$ will have, in the numerator,
only terms which include factors of $g$ or its derivatives, we
see that for each $k$, $\|\alpha\|_{C^k}$ is $O(e^{-C/\epsilon})$.

Next look at $\beta$. By construction,
$$0=\int_{X_b}\beta dx_1\wedge dx_2=\int_{X_b}\beta\theta_0\wedge Vdu
=\int_0^{\epsilon} W_0^{-1}(b-b_0)du.$$
Thus $b-b_0$, which is a function on $U_{r_1,r_2}\times S^1(\epsilon)$
(even though $b$ and $b_0$ are not) has no constant term in its
Fourier expansion. Both $b$ and $b_0$, however, are quasi-periodic in $u$,
i.e. consist of a linear plus periodic term. Let $\tilde b$ and $\tilde b_0$
denote the periodic part (not including the constant term) of $b$ and
$b_0$ respectively. Then $b-b_0= \tilde b -\tilde b_0$, and we can bound the
$C^k$ norm of $\tilde b$ and $\tilde b_0$ separately. For example, by (2.6),
$\partial_u \tilde b =2i\py g(u,y)$, which is $O(e^{-C/\epsilon})$, and then
the Poincar\'e inequality implies $\|\tilde b\|_{C^0}$ is $O(e^{-C/\epsilon})$.
Similar arguments apply to $\|\tilde b_0\|_{C^0}$, using the explicit form
for $b_0$, and from this one obtains $O(e^{-C/\epsilon})$ bounds on
$\|\beta\|_{C^k}$ for each $k$. Similar arguments apply for $\gamma
=-g(u,y)-W|b-b_0|^2$.

Thus the last three terms of $\omega_{new}$ are $O(e^{-C/\epsilon})$, and
since the sum of the first two terms have eigenvalues $O(\epsilon^{-1})$
and $O(\epsilon)$, it is clear that for sufficiently small $\epsilon$,
$\omega_{new}$ is positive definite. On the other hand, it is not clear
that $\int_X\omega_{new}^2=\int_X(\re\Omega)^2$. However, by construction
$\omega_{new}^2=(\re\Omega)^2$ outside of $f^{-1}(U_{r_1,r_2})$, and
$\omega_{new}^2$ and $(\re\Omega)^2$ differ only by $O(e^{-C/\epsilon})$
on $f^{-1}(U_{r_1,r_2})$, so
$\int_X\omega_{new}^2-\int_X(\re\Omega)^2=O(e^{-C/\epsilon})$. Now noting
that
$$([\omega_{new}]+aE)^2=[\omega_{new}]^2+2a\epsilon,$$ we can find a two-form
$\alpha$ on $B$ supported on $U_{r_1,r_2}$ with $C^k$ norm $O(e^{-C/\epsilon})$
such that $\int_X(\omega_{new}+f^*\alpha)^2
=\int_X(\re\Omega)^2$.
Set $\omega_{\epsilon}=\omega_{new}+f^*\alpha$. Because $\alpha$ is still small,
$\omega_{\epsilon}$ is still positive and defines the desired K\"ahler metric.
Properties (1)-(4) are then satisfied by construction. Note that $F_{\epsilon}
=\log(\Omega\wedge\bar\Omega/2\omega_{\epsilon}^2)$
is zero outside of $f^{-1}(U_{r_1,r_2})$, and $\omega_{\epsilon}$
is within $O(e^{-C/\epsilon})$ of $\omega_{SF}$ on $f^{-1}(U_{r_1,r_2})$.
Thus $\|F_{\epsilon}\|_{C^0}$ is $O(e^{-C/\epsilon})$. 
The same is true of $\|F_{\epsilon}\|_{C^2}$, and since the coefficients
of the metric are at worst $O(\epsilon)$ or $O(\epsilon^{-1})$ in
$f^{-1}(U_{r_1,r_2})$, $\|\Delta F_{\epsilon}\|_{C^0}$ is also
$O(e^{-C/\epsilon})$.
Furthermore, the Ricci form is $i\partial\bar\partial
F_{\epsilon}$, which is $O(e^{-C/\epsilon})$. This gives (5) and (6).

To bound the diameter of $X$ with the metric $\omega_{\epsilon}$,
first restrict the metric to the zero section $\sigma_0$ of $f$.
Identifying $\sigma_0$ with the base $B$, we note that on $B\setminus
\bigcup U_2^i$, the K\"ahler form of this restricted metric is
${i\over 2}(\epsilon W_0)((\epsilon W_0)^{-2}+|b_0|^2) dy\wedge d\bar y$.
But $b_0=0$ on $\sigma_0$, so this is just ${i\over 2}\epsilon^{-1}
W_0^{-1}dy\wedge d\bar y$. Let $D$ be the diameter of
$B\setminus \bigcup U_2^i$ under the metric ${i\over 2} W_0^{-1}dy
\wedge d\bar y$; this is independent of $\epsilon$. Thus $Diam(B
\setminus \bigcup U_2^i)=D\epsilon^{-1/2}$ under the metric induced by
$\omega_{\epsilon}$. On the other hand, the diameter of each fibre over
$B\setminus \bigcup U_2^i$ is bounded by some constant times $\epsilon^{1/2}$,
so $Diam(f^{-1}(B\setminus \bigcup U_2^i))\le D'\epsilon^{-1/2}$ for
sufficiently small $\epsilon$. Then applying Corollary 3.7 to each
$f^{-1}(U_2^i)$, (keeping in mind that the changes to the metric in 
the gluing area are negligible for small $\epsilon$), we see in
fact that
$$Diam(X)\le D'\epsilon^{-1/2}+D''\epsilon^{-1/2},$$
which we can always bound by $D_5\epsilon^{-1/2}$
for some constant $D_5$. 

Finally, (8) follows immediately from Proposition 3.8, Remark 2.7, and
the fact that any non-singular fibre has a neighbourhood in which
$\omega_{\epsilon}$ is arbitrarily close to the semi-flat metric for
$\epsilon$ sufficiently small. $\bullet$

\proclaim Theorem 4.5. Let $j:J\rightarrow\Pone$ be an elliptically
fibred K3 surface with section and 24 singular fibres over $\Delta=\{p_1,\ldots,
p_{24}\}$ as above. Then there exists open sets $U_1^i\subseteq U_2^i\subseteq
\Pone$,
$i=1,\ldots,24$, each diffeomorphic to a disc, $U^i_j\cap \Delta=\{p_i\}$,
positive constants $D_1,D_2,D_3,D_4,D_5$ and $\epsilon_0$ such that, for all
$\epsilon<\epsilon_0$, for any elliptic fibration $f:X\rightarrow\Pone$
with Jacobian $j:J\rightarrow\Pone$ with holomorphic 2-form $\Omega$ with
$[\re\Omega]^2=[\re\Omega_J]^2$,
and for any K\"ahler class $[\omega_{\epsilon}]$
on $X$ with $[\omega_{\epsilon}].X_b=\epsilon$ and $[\omega_{\epsilon}]^2
=[\re\Omega]^2=[\im\Omega]^2$,
there exists a K\"ahler metric $\omega_{\epsilon}$ representing 
$[\omega_{\epsilon}]$ on $X$ 
with the following properties:
\item{(1)} $\omega_{\epsilon}|_{f^{-1}(\Pone\setminus
\bigcup_i U_2^i)}$ is a semi-flat
metric (not necessarily the standard one).
\item{(2)} $\omega_{\epsilon}|_{f^{1}(U_1^i)}=T_{\sigma _i}^*\omega_{OV},$
where $\omega_{OV}$ is an Ooguri--Vafa metric and
$T_{\sigma _i}$ denotes translation by a (not necessarily holomorphic) section.
\item{(3)} If $F_{\epsilon}=\log\left({\Omega\wedge\bar\Omega/2\over 
\omega_{\epsilon}^2}\right)$, then
$$\|F_{\epsilon}\|_{C^0}\le D_1e^{-D_2/\epsilon}$$
and 
$$\|\Delta F_{\epsilon}\|_{C^0}\le D_1e^{-D_2/\epsilon},$$
where $\Delta$ denotes the Laplacian with respect to $\omega_{\epsilon}$.
\item{(4)} $$inf_v\{Ric(v,v)\, |\ |v|_{\omega_{\epsilon}}=1\}\ge -D_3
e^{-D_4/\epsilon}.$$
\item{(5)} With the Riemannian metric induced by $\omega_{\epsilon}$,
$Diam(X)\le D_5\epsilon^{-1/2}$.
\item{(6)}  If $R$ denotes the Riemann curvature tensor, then
$$\|R\|_{C^0}\le D_6\epsilon^{-1}\log\epsilon^{-1},$$
$$\hbox{$\|R\|_{C^0}\rightarrow\infty$ as $\epsilon\rightarrow 0$,}$$
and on any non-singular fibre, there exists a constant $C$ depending on
the fibre such that 
$$\|R\|\le C\epsilon.$$

\noindent {\it Proof.} First note that as in \S 1, we think of $X$ as a K3 surface
obtained from $J$ simply by altering the holomorphic 2-form
$\Omega_J$ on $J$ to $\Omega_J+j^*\alpha$, for some 2-form $\alpha$
on $\P ^1$. Thus it is natural to identify the underlying
manifolds $X$ and $J$, and we are only changing the complex structure. So we can
think of $[\omega_{\epsilon}]\in H^2(J,{\bf R})$, and in particular,
in the notation of \S 1, we can write 
$$[\omega_{\epsilon}]=\epsilon(\sigma_0+{\bf B})\mod E$$
for some ${\bf B}\in E^{\perp}/E\otimes {\bf R}$. Furthermore,
given the values of the classes $[\omega_{\epsilon}]$ and $[\Omega]$
modulo $E$, and given that $[\omega_{\epsilon}],[\re\Omega],[\im\Omega]$
form a hyperk\"ahler triple, the classes $[\omega_{\epsilon}]$ and
$[\Omega]$ are completely determined. Thus the choice of ${\bf B}$
uniquely determines the K\"ahler class and complex structure.

We modify the role of ${\bf B}$ slightly. Let $\omega_{\epsilon}^0$
be the K\"ahler form on $J$ provided by Theorem 4.4. Then in fact we
can write
$$[\omega_{\epsilon}]-[\omega^0_{\epsilon}]=\epsilon{\bf B}\mod E$$
for some ${\bf B}\in E^{\perp}/E\otimes {\bf R}$. The class
${\bf B}$ still determines all data. So fix this class in $E^{\perp}/E\otimes
{\bf R}$. This latter vector space is naturally identified with
$H^1(\Pone,R^1j_*{\bf R})$. Consider the exact sequence
$$\exact{R^1j_*{\bf R}}{C^{\infty}(\T_B^*)}{\F}.$$
Here $C^{\infty}(\T_B^*)$ denotes the sheaf of $C^{\infty}$ sections
of $\T_B^*$, and the first map is induced by tensoring the inclusion
$R^1j_*\boldz\hookrightarrow \T_B^*$ with ${\bf R}$. This gives a surjection 
$H^0(\Pone,\F)\rightarrow H^1(\Pone,R^1j_*{\bf R})$. Now a section
of $\F$ is given by an open covering $\{U_i\}$ of $\Pone$ and sections
$\sigma_i\in\Gamma(U_i,C^{\infty}(\T_B^*))$ with $\sigma_i-\sigma_j
\in \Gamma(U_i,R^1j_*{\bf R})$. This open covering $\{U_i\}$
can always be chosen with the following properties:
\item{(1)} Each $U_i$ contains at most one point of $\Delta$,
and if $p_j\in U_i$, then $\overline{U_2^j}\subseteq U_i$.
\item{(2)} Each $U_i$ is convex with respect to some metric on
$\Pone$, so that all multiple intersections of the $U_i$'s are contractible.
\item{(3)} If $U_i\cap \Delta=\phi$, then $U_i\cap \bigcup_j U_2^j=\phi$.

In fact, fixing one such open covering, all sections of $\F$ can be represented
over this open covering.

Now represent ${\bf B}$ by $(U_i,\sigma_i)$, and let $T_{\sigma_i}:
f^{-1}(U_i)\rightarrow f^{-1}(U_i)$ denote translation by the section
$\sigma_i$. Now consider the forms $T_{\sigma_i}^*\omega^0_{\epsilon}$
and $T_{\sigma_i}^*\Omega_J$. On $f^{-1}(U_i\cap U_j)$, $\omega^0_{\epsilon}$
is the standard semi-flat metric, 
by condition (3) on the open covering above, and since 
$\sigma_j-\sigma_i$ is a flat section with respect to the Gauss-Manin
connection, $T_{\sigma_j-\sigma_i}$ is an isometry, i.e. 
$T_{\sigma_j-\sigma_i}^*\omega^0_{\epsilon}
=\omega^0_{\epsilon}$, $T^*_{\sigma_j-\sigma_i}
\Omega_J=\Omega_J$ (see Example 2.2). Thus
$$T^*_{\sigma_i}\omega^0_{\epsilon}=T^*_{\sigma_i}T^*_{\sigma_j-\sigma_i}
\omega^0_{\epsilon}=T^*_{\sigma_j}\omega^0_{\epsilon},$$
and similarly $T^*_{\sigma_i}\Omega_J=T^*_{\sigma_j}\Omega_J$. Thus these
forms glue, to give global forms $\omega_{\epsilon}$, $\Omega$ on the
manifold $J$. The 2-form $\Omega$ satisfies $\Omega\wedge\Omega=0$, and thus
induces a new complex structure on $J$.  An easy local calculation shows that
$\Omega = \Omega_J+j^*\alpha '$, for some 2-form $\alpha '$
on $\P ^1$.  Furthermore,
it is clear that
$$\hbox{$\int_{J_b}\omega_{\epsilon}=\epsilon$, $\int_J\omega_{\epsilon}
\wedge \Omega=0$, and $\int_J \omega_{\epsilon}^2=\int_J (\re\Omega)^2
=\int_J (\im\Omega)^2$.}$$ 
Thus the cohomology classes $[\omega_{\epsilon}]$, $[\re\Omega]$
and $[\im\Omega]$ form a hyperk\"ahler triple. If we show that
$$[\omega_{\epsilon}]-[\omega^0_{\epsilon}]=\epsilon{\bf B}\mod E,$$
then we have constructed a K\"ahler form in the desired class (deducing 
moreover that the new complex structure is just that obtained from $X$).

To see the required identity, observe that we have
an exact sequence
$$H^2_{f^{-1}(\Delta)}(X,{\bf R})\mapright{\varphi} H^2(X,{\bf R})
\mapright{} H^2(X_0,{\bf R})\mapright{} H^3_{f^{-1}(\Delta)}(X,{\bf R}).$$
Now $H^2_{f^{-1}(\Delta)}(X,{\bf R})=H^0(f^{-1}(\Delta),{\bf R})={\bf R}^{24},$
and the image of $\varphi$ is just the one-dimensional subspace of
$H^2(X,{\bf R})$ spanned by $[E]$, the class of a fibre.
Thus it is enough to show that
$$[\omega_{\epsilon}|_{X_0}]-[\omega^0_{\epsilon}|_{X_0}]
=\epsilon{\bf B}\in H^2(X,{\bf R})/E \subseteq H^2(X_0,{\bf R}).$$
Now on $X_0$, $\omega^0_{\epsilon}$ is cohomologous to $\omega_{SF}$
by construction, and $\omega_{\epsilon}$ is cohomologous to a K\"ahler
form $\omega'_{SF}$ obtained in the same way as $\omega_{\epsilon}$ via
translation and gluing, but starting from $\omega_{SF}$ rather than 
$\omega^0_{\epsilon}$. Thus it is enough to show that on $X_0$
$$[\omega'_{SF}]-[\omega_{SF}]=\epsilon{\bf B}.$$

Over an open set $U_i$, write 
$$\omega_{SF}={i\over 2}W_0^{-1}(\epsilon\tv\wedge\tvc+\epsilon^{-1}
\th\wedge\thc).$$
Now
$$T_{\sigma_i}^*(\tv)=\tv+W_0(\py\sigma_i+b_0(\sigma_i(y),y))\th
+W_0(\pyc\sigma)\thc,$$
so over $U_i$ 
$$\eqalign{\omega_{SF}'-\omega_{SF}=&
T_{\sigma_i}^*\omega_{SF}-\omega_{SF}\cr
=&{\epsilon i\over 2}((\py\sigma_i+b_0(\sigma_i(y),y))\th\wedge\tvc
+\pyc\sigma_i\thc\wedge\tvc
\cr&+\py\bar\sigma_i\tv\wedge \th+(\pyc\bar\sigma_i+\bar b_0(\sigma_i(y),y))
\tv\wedge\thc))\mod \th\wedge\thc\cr
=&{\epsilon i\over 2}d(\sigma_i\tvc-\bar\sigma_i\tv)\mod\th\wedge \thc\cr}$$
as can be easily seen using (4.1) and a calculation similar to
that of (4.2).

How does the two-form $\omega_{SF}'-\omega_{SF}$ determine an element of
$H^1(B_0,R^1f_{0*}{\bf R})$? Given the open covering $\{U_i\}$ of
$B_0$, if $\omega_{SF}'-\omega_{SF}$ is an exact form on each
$f^{-1}(U_i)$, we can write $\omega_{SF}'-\omega_{SF}=d\alpha_i$
for some 1-form $\alpha_i$ on $f^{-1}(U_i)$. Then on
$f^{-1}(U_i\cap U_j)$, $\alpha_i-\alpha_j$ is closed, and hence
determines an element of $H^1(f^{-1}(U_i\cap U_j),{\bf R})
=\Gamma(U_i\cap U_j,R^1f_{0*}{\bf R})$ for our choice of open covering.
Now we have found such $\alpha_i$ modulo $\th\wedge\thc$, so
$\omega_{SF}'-\omega_{SF}$ is represented by a \v Cech cocycle
for $R^1f_{0*}{\bf R}$ given by
$$(U_i\cap U_j, {\epsilon i\over 2}((\sigma_i-\sigma_j)\tvc-(\bar\sigma_i
-\bar\sigma_j)\tv)).$$
By integrating this one-form over the periods, one sees this is precisely
the section of $\Gamma(U_i\cap U_j,R^1f_{0*}{\bf R})$ given by 
$\epsilon(\sigma_i-\sigma_j)$. Thus $\omega_{SF}'-\omega_{SF}$ represents the
class $\epsilon{\bf B}$.

Finally, properties (1)-(4) and (6) follow immediately from 
Theorem 4.4, (3)-(6) and (8).
On the other hand, the diameter of $f^{-1}(U_i)$ with respect to 
$\omega_{\epsilon}$ is the same as the diameter of $f^{-1}(U_i)$ 
with respect to the metric $\omega_{\epsilon}^0$. Since there are a fixed
number of $U_i$'s, the estimate on the diameter continues to hold
from Theorem 4.4, (7).  $\bullet$

\noindent {\bf Remark 4.6.} In the construction of the proof of
Theorem 4.5, we may sometimes want to be able to control the sections
$\sigma_i$ we use to represent the class ${\bf B}$. This can be
done as follows. The class of ${\bf B}$ depends on the choice of the
zero section $\sigma_0$. Changing the class of the zero section
changes ${\bf B}$ by an element of $E^{\perp}/E$. Thus ${\bf B}$
really should be thought of as living in
$E^{\perp}/E\otimes {\bf R}/\boldz$. (See [13] or [14], \S 7.)
Thus in some cases we might want to choose a compact set $F$ in
$E^{\perp}/E\otimes {\bf R}$ containing a fundamental domain for
$E^{\perp}/E$. We can then choose for each ${\bf B}\in F$ a representative
$(\sigma_i)$ of ${\bf B}$ with various norms, as required, bounded
by constants independent of ${\bf B}\in F$. We will say we are choosing
the $B$-field ${\bf B}$ in a fundamental domain for the $B$-field.
\bigskip

{\hd \S 5. Ricci-flat metrics.}

We will continue with the setting of Theorem 4.5. In other words, we have 
a fixed Jacobian elliptic fibration $j:J\rightarrow\Pone$. Our goal is to
show that there exists an $\epsilon_0$ such that for any 
$f:X\rightarrow\Pone$ with Jacobian $j:J\rightarrow\Pone$,
and any $\epsilon<\epsilon_0$, and any metric $\omega_{\epsilon}$ given
by Theorem 4.5, there exists a function $u_{\epsilon}$ such that
$\omega_{\epsilon}+i\partial\bar\partial u_{\epsilon}$ is a Ricci-flat
metric, and furthermore that $u_{\epsilon}$ is very small in
the $C^{k,\alpha}$ sense. Of course, that such a $u_{\epsilon}$ exists
is Yau's proof of the Calabi conjecture. Here we apply standard
techniques, following [20], to obtain control of
$u_{\epsilon}$. As mentioned in the introduction, the only subtle
difference is that as $\epsilon\rightarrow 0$, $Diam(X)\rightarrow\infty$,
and this requires us to be a bit more careful in estimating constants.
However, we follow [20] closely.

More precisely, we wish to solve the equations
$$\eqalign{
(\omega_{\epsilon}+i\partial\bar\partial u_{\epsilon})^2
&=e^{F_{\epsilon}}\omega_{\epsilon}^2\cr
\int_X u_{\epsilon}\omega_{\epsilon}^2&=0.\cr}\leqno{(5.1)}$$
Here $F_{\epsilon}=
\log({\Omega\wedge\bar\Omega/2\over \omega_{\epsilon}^2})$. 
By [35], we know such a
$u_{\epsilon}$ exists. 

We begin with some standard lemmas. For convenience, we will assume
$Vol(X)=1$. This can be achieved since we are holding the volume of
$X$ constant anyway, so we just scale the original $\Omega$ so
that $\int_J (\re\Omega)^2=1$.

\proclaim Lemma 5.1. Let $X$, $\omega_{\epsilon}$ be as in Theorem 4.5.
Assume $Vol(X)=1$. Then there exists a function
$I(\epsilon)$ depending only on $\epsilon$ and $J$ with $I(\epsilon)
\ge C\epsilon^{5}$, $C$ depending only on $J$, such that
\item{(1)} For any function $f$ on $X$ such that
$\int_X f\omega_{\epsilon}^2=0$,
$$\|df\|^2_2\ge I(\epsilon)\|f\|^2_4.$$
\item{(2)} For any function $f$ on $X$,
$$\|df\|^2_2\ge I(\epsilon)(\|f\|^2_4-\|f\|^2_2).$$

\noindent {\it Proof.} These are the standard Sobolev inequalities, but we just need to be
careful about the constants. We have, by [23], Lemmas 1 and 2,
for a function $f$
such that $\int_X f\omega_{\epsilon}^2=0$,
$$\|df\|_2^2\ge C_2\|f\|^2_4$$
while for an arbitrary function, we have
$$\|df\|^2_2 \ge D(4)C_2(\|f\|^2_4-\|f\|^2_2).$$
Here, we are using Li's notation for the constants $C_0, C_1, C_2, D(n)$
and the fact that the volume is 1 and the dimension is 4. Again by [23],
$D(4)$ is an absolute constant, $C_2=D(4)C_0^{1/2}$, and
$2C_1\ge C_0\ge C_1$, where $C_1$ is the constant in the isoperimetric
inequality
$$C_1(\min\{V(M_1),V(M_2)\})^3\le V(N)^4$$
where $V$ denotes volume, and $N$ is any codimension one submanifold
of $X$ dividing it into $M_1$ and $M_2$. In [9], Croke calls this constant
$\Phi(M)$.

Theorem 13 from [9] says that 
$$C_1\ge C_4\left(\int_0^{Diam(X)}  ((\sqrt{1/K})\sinh (\sqrt{K}r))^3 dr
\right)^{-5},$$
where $C_4$ again is an absolute constant, and $Ric(X)\ge -3K$, where
$3K\le D_3e^{-D_4/\epsilon}$ by
Theorem 4.5, (4). 
Now the integral is bounded above by 
$$Diam(X)(\sqrt{1/K}\sinh (\sqrt{K}Diam(X)))^3.$$
Now by Theorem 4.5, (5),
$\sqrt{K}Diam(X)\rightarrow 0$ as $\epsilon\rightarrow 0$,
so for sufficiently small $\epsilon$, using the first term of the Taylor
series expansion of $\sinh$, this is bounded by
$$C_5Diam(X)^4\le C_6\epsilon^{-2}$$
so $C_1\ge C_7\epsilon^{10}$, hence $C_0\ge C_8\epsilon^{10}$ and we can take
$$I(\epsilon)=\min(D(4),1)C_2\ge C_9\epsilon^{5}.\quad\bullet$$
\bigskip

\proclaim Lemma 5.2. (The $C^0$ estimate.) Let $u_{\epsilon}$ be the 
solution to equations (5.1). There exists a constant $C$
depending only on $J$, such that for all $\epsilon<\epsilon_0$,
($\epsilon_0, D_2$ as in Theorem 4.5)
$$\|u_{\epsilon}\|_{\infty}\le C\epsilon^{-5}e^{-D_2/\epsilon}.$$

\noindent {\it Proof.}
The starting point is the inequality
(23) of [20]:
$$\int_X|d|u_{\epsilon}|^{p/2}|^2\le Ap\int_X|F_{\epsilon}||u_{\epsilon}|^{p-1}.
$$
All integrals are with measure $\omega_{\epsilon}^2$. 
See also the expanded derivation of this inequality in [24]. 
One can check the constant
$A$ is independent of $p$ and $\epsilon$. 

We apply this first with $p=2$. The left-hand side is
$\|du_{\epsilon}\|^2_2\ge A_1\epsilon^{5}\|u_{\epsilon}\|_4^2$ by
Lemma 5.1, (1), so by applying
H\"older's inequality to the right-hand side,  we get
$$\left(\int_X |u_{\epsilon}|^4\right)^{1/2}\le
A_2\epsilon^{-5}\left(\int_X |F_{\epsilon}|^{4/3}\right)^{3/4}
\left(\int_X |u_{\epsilon}|^4\right)^{1/4}$$
or
$$\eqalign{\|u_{\epsilon}\|_4&\le A_3\epsilon^{-5}
\left(\int_X |F_{\epsilon}|^{4/3}\right)^{3/4}\cr
&\le C_1\epsilon^{-5} e^{-D_2/\epsilon}.\cr}\leqno{(5.2)}$$

Now for arbitrary $p$, using Lemma 5.1, (2)
$$\eqalign{
\|u_{\epsilon}\|^p_{2p}&=
\left(\int_X |u_{\epsilon}^{2p}|\right)^{1/2}=\| \, |u_{\epsilon}|^{p/2}\|^2_4\cr
&\le A_4\epsilon^{-5}\|\, d|u_{\epsilon}|^{p/2}\|^2_2+\| \, 
|u_{\epsilon}|^{p/2}\|^2_2 \cr
&\le A_5p\epsilon^{-5}\left(\int_X |F_{\epsilon}|\, |u_{\epsilon}|^{p-1}\right)
+\|\, |u_{\epsilon}|^{p/2}\|^2_2.\cr}$$
Applying H\"older's inequality to the first term, we have, with $q=p$,
$q'=1/(1-1/p)=p/(p-1)$,
$$\eqalign{
\int_X |F_{\epsilon}|\, |u_{\epsilon}|^{p-1}&\le \|F_{\epsilon}\|_p \, 
\|\, |u_{\epsilon}|^{p-1}\|_{p/(p-1)}\cr
&=\|F_{\epsilon}\|_p \, \| u_{\epsilon}\|^{p-1}_p,\cr}$$
so
$$\eqalign{
\|u_{\epsilon}\|^p_{2p}&\le A_5p\epsilon^{-5}\|F_{\epsilon}\|_p\, 
\|u_{\epsilon}\|_p^{p-1}+\|u_{\epsilon}\|_p^p\cr
&=\left( A_5p\epsilon^{-5}\|F_{\epsilon}\|_p+\|u_{\epsilon}\|_p \right) \, \|u_{\epsilon}
\|_p^{p-1}.\cr}\leqno{(5.3)}$$

Now we claim that if we set $p_n=2^{n+1}$, there exists constants 
$C_n$ such that
$$\|u_{\epsilon}\|_{p_n}\le C_n\epsilon^{-5}e^{-D_2/\epsilon}$$
for all $\epsilon<\epsilon_0$.
This holds for $n=1$ by (5.2). Suppose it holds for
a given $n$. Then by (5.3),
$$\eqalign{\|u_{\epsilon}\|_{p_{n+1}}^{p_n}&
\le (A_5p_n\epsilon^{-5}D_1e^{-D_2/\epsilon}+
C_n\epsilon^{-5}e^{-D_2/\epsilon})
(C_n\epsilon^{-5}e^{-D_2/\epsilon})^{p_n-1}\cr
&\le\cases{(A_5D_12^{n+1}+1)(C_n\epsilon^{-5}e^{-D_2/\epsilon})^{p_n}&
if $C_n\ge 1$;\cr
(A_5D_12^{n+1}+1)(\epsilon^{-5}e^{-D_2/\epsilon})^{p_n}&
if $C_n\le 1$.\cr}\cr}$$
Thus we can take
$$C_{n+1}\le\cases{ (A_5D_12^{n+1}+1)^{2^{-(n+1)}}C_n& if $C_n\ge 1$;\cr
(A_5D_12^{n+1}+1)^{2^{-(n+1)}}&if $C_n\le 1$.\cr}$$
It then follows as in [20], page 299, that $C_n\le A_6$ for some
constant $A_6$ independent of $n$ and $\epsilon$. Thus we conclude that
$$\|u_{\epsilon}\|_{\infty}\le A_6\epsilon^{-5}e^{-D_2/\epsilon}$$
for all $\epsilon<\epsilon_0$. $\bullet$
\bigskip

\proclaim Lemma 5.3. (The $C^2$ estimate.) Let $u_{\epsilon}$ be the 
solution to equations (5.1). There are constants
$C$ and $\epsilon_0$ depending only on $J$ (possibly smaller than
the $\epsilon_0$ of Theorem 4.5) such that for all $\epsilon<\epsilon_0$,
$$C^{-1}\omega_{\epsilon}\le\tilde\omega_{\epsilon}\le C\omega_{\epsilon}$$
where $\tilde \omega_{\epsilon}=\omega_{\epsilon}+i\partial\bar\partial 
u_{\epsilon}$.

\noindent {\it Proof.} 
Let $R_{\epsilon}=\sup_{i\not= j}|R_{i\bar i j\bar j}|$, where 
$R_{i\bar  i j\bar j}$ is the holomorphic bisectional curvature
of the metric $\omega_{\epsilon}$, and the supremum is over all
points of $X$ and unitary bases at each point. Since the holomorphic
bisectional curvature determines the curvature, ([6], pg. 76)
and $\sup\|R\|\rightarrow\infty$ as $\epsilon\rightarrow 0$ by
Theorem 4.5 (6), we must have
$R_{\epsilon}>1$ for small $\epsilon$. So if $c_{\epsilon}
=2R_{\epsilon}$, then $c_{\epsilon}+\inf R_{i\bar ij\bar j}\ge R_{\epsilon}
>1$. Here the infinum is as before over all unitary frames and points on $X$.
Then [35], (2.22), reads
$$e^{c_{\epsilon}u_{\epsilon}}\Delta'(e^{-c_{\epsilon}u_{\epsilon}}
(2+\Delta u_{\epsilon}))
\ge (\Delta F_{\epsilon}-4\inf_{i\not=j} R_{i\bar ij\bar j}(x))
-2c_{\epsilon}(2+\Delta u_{\epsilon})
+(c_{\epsilon}+\inf_{i\not=j} R_{i\bar ij\bar j}(x))e^{-F_{\epsilon}}
(2+\Delta u_{\epsilon})^2$$
where the infina are now only at the given point (but still
over all unitary bases).
Here $\Delta'$ is the Laplacian with respect to the metric
$\omega_{\epsilon}+i\partial\bar\partial u_{\epsilon}$, and $\Delta$
is the Laplacian with respect to $\omega_{\epsilon}$. Let
$$k(x)=-\inf_{i\not=j} R_{i\bar ij\bar j}(x)/R_{\epsilon},$$
so that $|k(x)|\le 1$.

Now suppose $e^{-c_{\epsilon}u_{\epsilon}}(2+\Delta u_{\epsilon})$
assumes its maximum at $x\in X$. Then by the maximum principal,
the Laplacian must be non-positive there, so at the point $x$
$$\eqalign{
0\ge& e^{c_{\epsilon}u_{\epsilon}}\Delta'(e^{-c_{\epsilon}u_{\epsilon}}
(2+\Delta u_{\epsilon}))\cr
\ge& (\Delta F_{\epsilon}+4 k(x) R_{\epsilon})-2c_{\epsilon}(2+\Delta 
u_{\epsilon})
+(c_{\epsilon}-k(x)R_{\epsilon})e^{-F_{\epsilon}} (2+\Delta u_{\epsilon})^2\cr
=&(\Delta F_{\epsilon} +4 k(x)R_{\epsilon})-4R_{\epsilon} 
(2+\Delta u_{\epsilon})
+(2-k(x))R_{\epsilon} e^{-F_{\epsilon}}(2+\Delta u_{\epsilon})^2\cr
=&e^{-F_{\epsilon}}(2-k(x))R_{\epsilon}
\left[ \left( (2+\Delta u_{\epsilon})-{2e^{F_{\epsilon}}\over
2-k(x)}\right)^2-\left({2e^{F_{\epsilon}}\over 2-k(x)}\right)^2
+{e^{F_{\epsilon}}(\Delta F_{\epsilon}+4 R_{\epsilon} k(x))\over
(2-k(x))R_{\epsilon}}\right]\cr}$$
and since $|k(x)|\le 1$, we get
$$\left | (2+\Delta u_{\epsilon})-{2e^{F_{\epsilon}}\over 2-k(x)}\right|
\le \left | \left({2e^{F_{\epsilon}}\over 2-k(x)}\right)^2 
- {e^{F_{\epsilon}}(\Delta F_{\epsilon}+4 R_{\epsilon} k(x))\over
(2-k(x)) R_{\epsilon}}\right|^{1/2}$$
If we are outside of the region where the gluing is taking place, then
$F_{\epsilon}=0$, so we get
$$\left | (2+\Delta u_{\epsilon})-{2\over 2-k(x)}\right |
\le \left | \left({2\over 2-k(x)}\right)^2-{4 k(x)\over 2-k(x)}\right |^{1/2},$$
or
$$\eqalign{2+\Delta u_{\epsilon}
&\le {2\over 2-k(x)}+\left |\left ({2\over 2-k(x)}\right)^2 -
{4k(x)\over 2-k(x)}\right|^{1/2}\cr
&=2.\cr}$$

In the gluing region, by Theorem 4.5 (6), there is a constant $C_1$
$$|k(x)|\le C_1\epsilon/R_{\epsilon}.$$

Also in the gluing region, we can use the bounds of Theorem 4.5, (3)
on $F_{\epsilon}$ and $\Delta F_{\epsilon}$, to get, for
a constant $C_2$ bounding $e^{F_{\epsilon}}$,
$$\eqalign{
2+\Delta u_{\epsilon}&\le
{2 e^{F_{\epsilon}}\over 2-k(x)}+\left |\left ( {2e^{F_{\epsilon}}
\over 2-k(x)}\right)^2- {e^{F_{\epsilon}}(\Delta F_{\epsilon}+4 R_{\epsilon}
k(x)))\over (2-k(x))R_{\epsilon}}\right|^{1/2}\cr
&\le {2C_2\over 2-C_1\epsilon/R_{\epsilon}}+
\left | \left ( {2C_2\over 2-C_1\epsilon/R_{\epsilon}}\right )^2
+{C_2(D_1e^{-D_2/\epsilon}+4C_1\epsilon)\over
2-C_1\epsilon/R_{\epsilon}}\right |^{1/2}.\cr}$$

Now as $\epsilon\rightarrow 0$, $\epsilon/R_{\epsilon}\rightarrow 0$ by
Theorem 4.5, (6).  
So what we get is
$$(2+\Delta u_{\epsilon})(x)\le C_3$$
for sufficiently small $\epsilon$, and $C_3$ independent of $\epsilon$.

Now
$$e^{-c_{\epsilon}u_{\epsilon}(y)}(2+\Delta u_{\epsilon})(y)
\le e^{-c_{\epsilon}u_{\epsilon}(x)}(2+\Delta u_{\epsilon})(x)$$
for all points $y$, so
$$\eqalign{2+\Delta u_{\epsilon}
&\le e^{c_{\epsilon}(u_{\epsilon}(y)-u_{\epsilon}(x))} C_3\cr
&\le e^{c_{\epsilon}(\sup u_{\epsilon}-\inf u_{\epsilon})} C_3\cr
&\le e^{R_{\epsilon}C_4\epsilon^{-5}e^{-D_2/\epsilon}}C_3.\cr}$$
By Theorem 4.5, (6),
$$R_{\epsilon}\epsilon^{-5}e^{-D_2/\epsilon}\rightarrow 0,$$
so we get
$$2+\Delta u_{\epsilon}\le C_5$$ 
for sufficiently small $\epsilon$.

Now working in a choice of coordinates $z_1,z_2$ at a point so that 
$\partial_{z_1},\partial_{z_2}$ are unitary at the point with respect
to $\omega_{\epsilon}$ and which also diagonalizes $\tilde \omega_{\epsilon}
=\omega_{\epsilon}+i\partial\bar\partial u_{\epsilon}$, 
then 
$$(\tilde\omega_{\epsilon})_{i\bar j}=\delta_{ij}(1+(u_{\epsilon})_{i\bar i}),$$
and each $1+u_{i\bar i}$ is positive, so $1+(u_{\epsilon})_{i\bar i}
\le C_5$, so $\tilde\omega_{\epsilon}\le C_5\omega_{\epsilon}$.
Also,
$$\tilde\omega^2_{\epsilon}=\prod (1+(u_{\epsilon})_{i\bar i})
\omega^2_{\epsilon}=e^{F_{\epsilon}}\omega^2_{\epsilon}.$$
Since $1+(u_{\epsilon})_{i\bar i}$ is bounded above by $C_5$,
it must be bounded below by something close to $C_5^{-1}$, so changing
$C_5$ slightly if necessary, we get
$$C_5^{-1}\omega_{\epsilon}\le\tilde\omega_{\epsilon}\le C_5\omega_{\epsilon}.
\quad\bullet$$

We note here that for some purposes, Lemma 5.3 is already sufficient.
For example, if we wish to know that the fibres collapse to points
as $\epsilon\rightarrow 0$, Lemma 5.3 along with Proposition 3.5 tells
us the diameter of each fibre under the Ricci-flat metric
goes to zero as $\epsilon\rightarrow 0$. However, if we wish
to get a clearer picture of the asymptotic behaviour of
the metric, we need stronger results.

\proclaim Lemma 5.4. (The $C^{2,\alpha}$ estimate.) Let $u_{\epsilon}$
be the solution to equations (5.1). If $U\subseteq B$
is a simply connected open set with $\overline{U}\subseteq B_0
=B\setminus \Delta$, then there exists constants $\alpha$ and $\epsilon_0$
and a polynomial $P$,
depending on $J$ and $U$, such that 
$$\|u_{\epsilon}\|_{C^{2,\alpha}}\le P(\epsilon^{-1})$$
in $f^{-1}(U)$ for all $\epsilon<\epsilon_0$ and ${\bf B}$ in a fundamental
domain for the B-field (see Remark 4.6). 
Here the $C^{2,\alpha}$ norm is on $f^{-1} (U)$ as defined in Lemma 4.1,
and so $\alpha$, $\epsilon_0$ and
$P$ also depend on the choice of holomorphic coordinate $y$
and fixed bounded domain $T'$, as specified in the proof below.

\noindent {\it Proof.} We need to apply the basic result of [11], Theorem 17.14.
However, we must be careful about the constants. Let $\pi:\T_B^*\rightarrow
B$ be the projection, and let $T'\subseteq
\pi^{-1}(U)\subseteq\T_{B}^*$ be a fixed bounded domain which contains
a fundamental domain of each fibre of $f$ over $U$. We will be
computing norms in the domain $T'$. To do so,
we choose a holomorphic coordinate $y$ in the base, yielding holomorphic canonical 
coordinates $x,y$ on $\T_U^*$. Now take a bigger open set
$T(\epsilon)$ containing $T'$. This open set will also be bounded,
but will depend on $\epsilon$. We choose it as follows. First let
$V\subseteq B_0$ be an open set with $\overline{U}\subseteq V$, $\overline{V}
\subseteq B_0$, 
and the holomorphic coordinate $y$ extending to $V$. Let
$T\subseteq \pi^{-1}(V)$ be a domain containing $T'$ and containing
a fundamental domain of each fibre over $V$. Let
$$T(\epsilon)=\{(x,y)\in\T_{V}^*\quad|\quad\hbox{there exists
$(\tilde x,y)\in T$ with $|x-\tilde x|<\epsilon^{-1/2}$}\}.$$
The point of this choice is as follows. 
Consider
the change of variable $y'=\epsilon^{-1/2}y$, $x'=\epsilon^{1/2}x$. 
Then using $x',y'$ to identify $\T_U^*$ with a subset of ${\bf C}^2$,
we get $Dist(\partial\overline{T(\epsilon)},T')\ge 1$, for sufficiently
small $\epsilon$, in
the euclidean distance in ${\bf C}^2$. 

Pulling back $\omega_{\epsilon}$ to $T(\epsilon)$, we can write
$$\omega_{\epsilon}=i\partial\bar\partial(\varphi_1+\varphi_2)$$
where $i\partial\bar\partial\varphi_1$ is a semi-flat metric
and $i\partial\bar\partial\varphi_2$ is the correction to this metric
resulting from the gluing process. By applying Lemma 4.1, we can choose
$\epsilon$ sufficiently small so that the $C^{2,\alpha}$ norm of $\varphi_2$
on $T'$ is as small as we like (and the $C^2$ norm of $\varphi_2$ on
$T(\epsilon)$).
On the other hand, $\varphi_1$ can be taken to be
a translation of the K\"ahler potential given for the standard
semi-flat metric in Example 2.2. Since we have chosen ${\bf B}$ in
a fundamental domain, we can then bound the $C^{2,\alpha}$ norm
of $\varphi_1$ on $T'$ independently of ${\bf B}$ as a 
polynomial in $\epsilon^{-1}$. The same is true of the $C^2$ norm of
$\varphi_1$ on $T'(\epsilon)$.

Now the equation that $u_{\epsilon}$ satisfies is
$$(i\partial\bar\partial\psi_{\epsilon})^2=\Omega\wedge\bar\Omega/2$$
where $\psi_{\epsilon}=\varphi_1+\varphi_2+u_{\epsilon}$. Thus a
$C^{2,\alpha}$ bound on $\psi_{\epsilon}$ polynomial in $\epsilon^{-1}$
yields a $C^{2,\alpha}$ bound on $u_{\epsilon}$ polynomial in $\epsilon^{-1}$.
Now changing coordinates between $x,y$ and $x',y'$ also only
affects the $C^{2,\alpha}$ norm of a function by a factor polynomial
in $\epsilon^{-1}$, so we can work with respect to the coordinates
$x',y'$. Now in these coordinates,
$$i\partial\bar\partial(\varphi_1)={i\over 2}
(W_0(dx'+\epsilon bdy')\wedge\overline{(dx'+\epsilon bdy')}+W_0^{-1} dy'
\wedge d\bar y').$$
By looking at the explicit form of $b$ for the semi-flat metric,
we see $\epsilon b$ in fact goes to zero as $\epsilon\rightarrow 0$ on
$T(\epsilon)$. Thus the eigenvalues of $i\partial\bar\partial\varphi_1$
on $T(\epsilon)$, i.e. the eigenvalues of the matrix
$$\pmatrix{W_0&\epsilon bW_0\cr \epsilon\bar bW_0& W_0^{-1}+\epsilon^2|b|^2
\cr},$$ 
can be bounded below and above by some constants $\lambda$ and $\Lambda$
independently of $\epsilon$. Since $\varphi_2$ is small, the same is true of
$\omega_{\epsilon}$ on $T(\epsilon)$. Finally, by Lemma 5.3,
the eigenvalues of $i\partial\bar\partial\psi_{\epsilon}$ are bounded
below and above by $C^{-1}\lambda$ and $C\Lambda$, independently of
$\epsilon$. Furthermore, Lemmas 5.2 and 5.3 imply the $C^2$ norm
of $\psi_{\epsilon}$ on $T(\epsilon)$ is bounded
by a polynomial in $\epsilon^{-1}$.
We can now apply [11], Theorem 17.14 to the domains
$T'\subseteq T(\epsilon)$, to obtain the desired result. $\bullet$

\bigskip
We shall now follow the standard method of continuity from [35], and,
for $t\in [0,1]$, look 
at the solution $u_{\epsilon,t}$ to the equation
$$(\omega_{\epsilon}+i\partial\bar\partial u_{\epsilon,t})^2
=(1+t(e^{F_\epsilon }-1))\omega_{\epsilon}^2,\leqno{(5.4)}$$ 
$$\int_X u_{\epsilon,t}\ \omega_{\epsilon}^2=0.\leqno{(5.5)}$$

We set $\omega _{\epsilon ,t} = \omega_{\epsilon}+i\partial\bar\partial u_{\epsilon,t}$,
the K\"ahler form of a metric on the given complex manifold $X$.  For $t=0$, we just get 
back our original (glued) metric, whilst $t=1$ is the case we have just looked at, 
yielding the Ricci flat metric with K\"ahler form $\tilde \omega _{\epsilon }$.  
Since $\log (1+t(e^{F_{\epsilon}}-1))$ has the same properties as
$F_{\epsilon}$ for $t\in [0,1]$, all the above estimates of Lemmas 5.2--5.4
work
equally well for $u_{\epsilon,t}$.  In particular,
$$C^{-1}\omega_{\epsilon}\le \omega_{\epsilon ,t}\le C\omega_{\epsilon} $$
for some constant $C$ independent of $t\in [0,1]$ and $\epsilon$, and 
$$\|u_{\epsilon,t}\|_{C^{2,\alpha}}\le P(\epsilon ^{-1}),$$
with the polynomial $P$ independent of $t\in [0,1]$ and $\epsilon$.

Moreover, the Ricci form of the metric $\omega_{\epsilon ,t}$ is given by
$$ {i\over 2\pi} \partial \bar\partial (F_{\epsilon} - \log 
(1 + t(e^{F_\epsilon} -1))),$$ and so the Ricci curvature 
$Ric_{\omega_{\epsilon ,t}}$ has a similar lower bound (independent of $t$) 
as $Ric _{\omega _\epsilon }$.
\bigskip

\proclaim Lemma 5.5.  Let $G_{\epsilon ,t} (x,y)$ denote Green's function for the Laplacian 
$\Delta _{\epsilon ,t}$ associated to the metric $\omega_{\epsilon ,t}$, normalised so that 
$\int _X G_{\epsilon ,t} (x,y) \,\omega^2_{\epsilon,t}(x) = 0$.  Then, for $\epsilon$ sufficiently small and 
any $t \in [0,1]$,
$$G_{\epsilon ,t} (x,y) \ge - A \epsilon ^{-11},$$ for some constant $A$ independent
of $\epsilon$ and $t$.

\noindent {\it Proof.}  For ease of notation, we drop the suffices $\epsilon ,t$.
We follow the proof of Lemma 3.3 from [24], which is due to Peter Li.  
The volume of $X$ is 1, and we 
set $K(x,y,s) = H(x,y,s) -1$, where $H$ is the heat kernel on $X$.  
As in [24], we need to find a lower bound for the integral of $K(x,y,s)$ 
over $1\le s \le \infty$, of the same form as that claimed for $G(x,y)$.  Lu observes that 
$$ K(x,y,s) \ge - K^{1/2} (x,x,s) K^{1/2} (y,y,s) ,$$ and that furthermore, for any $x\in X$, 
$$ K(x,x,s) \le K(x,x,1) e^{-\lambda (s-1)},$$
for all $s\ge 1$, where $\lambda$ denotes the first (positive) eigenvalue of the Laplacian.  If now
we can suitably bound $\lambda$ from below, and $K(x,x,1)$ from above, we'll be able to integrate 
the resulting function which bounds $K^{1/2} (x,x,s) K^{1/2} (y,y,s)$ from above, obtaining a 
lower bound for $\int _1 ^\infty K(x,y,s) ds$.

The bound from below for $\lambda$ comes from Theorem 4 on page 116 of [31].  Since the metric is 
within a fixed constant factor of our original metric, all the quantities in the given formula are known,
and so using Theorem 4.5, and we deduce that
$$\lambda \ge A_1 Diam (X) ^{-2} \ge  A_2 \epsilon ,$$  for appropriate absolute constants
$A_1 , A_2$.  The proof of Lemma 5.1 may be applied to the metric $\omega_{\epsilon ,t}$ to 
obtain a similar bound on the Sobolev constant, and then 
the bound from above for
$K(x,x,1)$  is implied by equation (3.12) of [36], where the argument given there has been run for 
the function $K(x,y,s)= H(x,y,s) -1$ (so in particular $\int _X K(x,z,s)\, 
\omega^2_{\epsilon,t}(z)=0$).  
For an appropriate constant $A_3$ independent of $t$, we have 
$$ K(x,x,1) \le A_3 \epsilon ^{-10}.$$  Thus for all $s\ge 1$
$$ K(x,x,s) \le  A_3 \epsilon ^{-10} e^{- A_2 \epsilon  (s-1)},$$
which then implies that $$K(x,y,s) \ge -A_4 \epsilon ^{-10} e^{- A_2
\epsilon  (s-1)},$$ for some constant $A_4$ independent of $\epsilon$ and $t$.  
On 
integrating, we obtain the claimed bound in the form stated (a rather more involved argument 
in fact gives a bound $K(x,x,1) \le A'_3 \epsilon ^{-3}$, and hence 
$G_{\epsilon ,t} (x,y) \ge - A' \epsilon ^{-4}$, but this extra accuracy is not required). $\bullet$
\bigskip

We are now ready for our main theorem.

\proclaim Theorem 5.6. For any simply connected open set
$U\subseteq B_0$ with $\overline{U}\subseteq B_0$, and any
$k\ge 2$, $0<\alpha<1$, there exists constants $C, C'$,
and $\epsilon_0$ such that for all choices of ${\bf B}$ in a fundamental
domain for the B-field and any $\epsilon<\epsilon_0$ giving $\omega_{\epsilon}$
as in Theorem 4.5, and $u_{\epsilon}$ satisfying the equations
$$\eqalign{
(\omega_{\epsilon}+i\partial\bar\partial u_{\epsilon})^2
&=e^{F_{\epsilon}}\omega_{\epsilon}^2\cr
\int_X u_{\epsilon}\omega_{\epsilon}^2&=0\cr}$$
with $F_{\epsilon}=
\log({\Omega\wedge\bar\Omega/2\over \omega_{\epsilon}^2})$, 
we have
$$\|u_{\epsilon}\|_{C^{k,\alpha}}\le Ce^{-C'/\epsilon}.$$
Here, the norm is as in Lemma 4.1 on the region $f^{-1} (U)$, and
the constants $C,C'$ are independent of $\epsilon$.

\noindent {\it Proof.} This is now completely standard. Following [20] and [24],
we differentiate (5.4) with respect to $t$, getting
$$\Delta_{\epsilon,t}{d u_{\epsilon,t}\over dt}
={e^{F_{\epsilon}}-1\over 1+t(e^{F_{\epsilon}}-1)}.$$
The right hand side is very small, which along with the estimate
of Lemma 5.5, allows us to bound $du_{\epsilon,t}/dt$. Indeed, by
Green's formula and (5.5), we have
$${d u_{\epsilon,t}(x)\over dt}
=-\int_X\left(\Delta_{\epsilon,t}{du_{\epsilon,t}\over dt}\right)
\tilde G_{\epsilon,t}(x,y)\omega^2_{\epsilon,t}(y).$$
Here $\tilde G_{\epsilon,t}$ is the Green's function for the Laplacian
for $\omega_{\epsilon,t}$, normalized so that $\inf_X \tilde G_{\epsilon,t}
=0$. Lemma 5.5 tells us that
$\int_X \tilde G_{\epsilon,t}(x,y)\omega_{\epsilon,t}^2(y)\le A\epsilon^{-11}$
for some constant $A$ independent of $\epsilon$ and $t$, so bounds on 
$F_{\epsilon}$ imply
$$\|du_{\epsilon,t}/dt\|_{C^0}\le C_1e^{-C_2/\epsilon}\leqno{(5.6)}$$
for some constants $C_1$ and $C_2$ independent of $t$ and $\epsilon$,
for sufficiently small $\epsilon$.

We can now apply the interior Schauder estimates (see [11] Theorem 6.2)
to obtain
$$\|du_{\epsilon,t}/dt\|_{C^{2,\alpha}}\le C_3e^{-C_4/\epsilon}\leqno{(5.7)}$$
for sufficiently small $\epsilon$. This holds for $\alpha$ as given by Lemma
5.4.
We note that a certain amount of care must be taken in applying
these estimates: first, we need to use the estimate of (5.6)
and Lemma 5.4 on a larger open set $U'$ with
$\overline{U}\subseteq U'\subseteq B_0$. Second we note that by
Lemma 5.4, the $C^{0,\alpha}$ estimates for the coefficients of
the second order operator $\Delta_{\epsilon,t}$ depend only polynomially
on $\epsilon^{-1}$, and the same is true, much as in the proof
of Lemma 5.4, for the constants $\lambda$ and $\Lambda$ needed in applying
[11], Theorem 6.2. The constant arising in the Schauder estimate can
be verified to depend only polynomially on $\lambda$ and $\Lambda$. Taking
these things into account, one obtains (5.7). 

Now integrating (5.7) with respect to $t$ we obtain
$$\|u_{\epsilon}\|_{C^{2,\alpha}}\le C_3 e^{-C_4/\epsilon}.$$ 
Using Schauder estimates again repeatedly  in the standard way (see
[35], Formula (4.5) and following text), one can then find for each
$k$, constants $C$ and $C'$ such that 
$$\|u_{\epsilon}\|_{C^{k,\alpha}}\le Ce^{-C'/\epsilon}.$$
To get this inequality for any $\alpha$, one uses the interpolation
inequalities. $\bullet$
\bigskip
\noindent {\bf Remark 5.7.}  The construction of the Ooguri--Vafa metric in
\S 3 clearly works also for singular fibres of type $I_n$, simply by quotienting
at the appropriate stage by $\epsilon n \Z$ instead of $\epsilon \Z$, and the above 
proofs go through unchanged in this case.  Thus all the results of this section 
remain valid for elliptic K3 surfaces with semi-stable fibres.

\bigskip
{\hd \S 6. Gromov--Hausdorff convergence.}

We now return to the notion of convergence alluded to in the
introduction. We wish to show that with the proper normalization,
the results of \S 5 imply that in the large complex structure
limit, K3 surfaces
in fact converge to 2-spheres. To make this precise, we
first recall the notion of Gromov--Hausdorff distance. The definition
given below can be easily seen to be equivalent to a
definition in terms of $\epsilon$-dense
subsets, c.f. [30] pg. 276.

\proclaim Definition 6.1. Let $(X,d_X)$, $(Y,d_Y)$ be two compact
metric spaces. Suppose there exists maps $f:X\rightarrow Y$
and $g:Y\rightarrow X$ (not necessarily continuous) such that
for all $x_1,x_2\in X$,
$$|d_X(x_1,x_2)-d_Y(f(x_1),f(x_2))|<\epsilon$$
and for all $x\in X$,
$$d_X(x,g\circ f(x))<\epsilon,$$
and the two symmetric properties for $Y$ hold. Then we say the
Gromov--Hausdorff distance between $X$ and $Y$ is at most $\epsilon$.
The Gromov--Hausdorff distance $d_{GH}(X,Y)$ is the infinum of all
such $\epsilon$.

The Gromov--Hausdorff distance defines a topology on the set
of compact metric spaces, and hence a notion of convergence. It 
follows from results of Gromov (see e.g. [30], pg. 281, Cor. 1.11)
that the class of compact Ricci-flat manifolds with diameter $\le D$
is precompact. Thus in particular, if we have a sequence of
Calabi--Yau $n$-folds whose complex structure converges to
a large complex structure limit point (or any other boundary point for
that matter) and whose metrics have diameter bounded above,
then there is a convergent subsequence, and then the basic question is:
what is the limit? The conjecture which motivated the work
of this paper is the following:

\proclaim Conjecture 6.2. Let $\overline{\M}$ be a compactified moduli space
of complex deformations of a simply-connected Calabi--Yau $n$-fold
$X$ with holonomy group $SU(n)$, and let $p\in\overline{\M}$ be a 
large complex structure limit
point (see [27] for the precise Hodge-theoretic definition of 
this notion). Let $(X_i,g_i)$ be a sequence of Calabi--Yau manifolds
with Ricci-flat K\"ahler
metric which are complex deformations of $X$, with the sequence
$[X_i]\in \overline{\M}$ converging suitably to $p$, and $C_1\ge Diam(X_i)\ge C_2>0$
for all $i$. Then a subsequence of $(X_i,g_i)$ converges
to a metric space $(X_{\infty},d_{\infty})$, where $X_{\infty}$ is
homeomorphic to $S^n$. Furthermore, $d_{\infty}$ is induced by a Riemannian
metric on $X_{\infty}\setminus\Delta$, with $\Delta\subseteq X_{\infty}$
a set of codimension 2.

A similar conjecture was also made by Kontsevich, Soibelman and Todorov
(see [22], [25]).\medskip

\noindent {\bf Remark 6.3.} Conjecture 6.2
is obvious in the elliptic curve case (ignoring the fact that elliptic curves
are not simply-connected), no matter how the sequence of points approaches
the large complex structure limit point. However, in the K3 case, more care must
be taken. In this paper, we have considered limits mirror to points approaching
the large K\"ahler limit along a ray in the K\"ahler cone. However, if a sequence
of points approaching the large K\"ahler limit approaches the boundary
of the projectivized K\"ahler cone, we might expect further degeneracies
in the Gromov-Hausdorff limits. For example, a product of two elliptic
curves $E_1\times E_2= {\bf R}^4/\boldz^4$, with a metric
$\pmatrix{\epsilon^{-3}&0&0&0\cr 0&\epsilon&0&0\cr 0&0&\epsilon&0\cr 
0&0&0&\epsilon\cr}$ has a special Lagrangian fibration given by
projection on the the first and third factors, and has fibres of area $\epsilon$.
When we normalise the metrics to have diameter one, the sequence 
of Riemannian manifolds  
converges to an $S^1$  as $\epsilon\rightarrow 0$. As pointed out to
us by N.C. Leung, this construction descends to the corresponding
Kummer surfaces.  The limit of the Kummer surfaces is then a closed interval.

Thus we expect that the correct restriction on sequences of points in the complex
moduli space in Conjecture 6.2 should correspond in the mirror to
K\"ahler classes staying within a proper
subcone of the K\"ahler cone. We can now prove the conjecture for
the limits of K3 surfaces considered in this paper, where the K\"ahler class
tends to $\infty$ along a ray, which we have seen reduces to the following result.

\proclaim Theorem 6.4. Let $j:J\rightarrow B$ be an
elliptically fibred K3 surface with a section and singular fibres all of type $I_1$,
and let $f_i:X_i\rightarrow
B$ be a sequence of elliptically fibred K3 surfaces with jacobian $j$.
Let $\omega_i$ correspond to a Ricci-flat K\"ahler 
 metric on $X_i$ with  $\omega_i^2$ independent of $\, i$,
and with $\int_{f_i^{-1}(b)}\omega_i=\epsilon_i\rightarrow
0$ as $i\rightarrow\infty$.  Then the sequence of Riemannian manifolds 
$(X_i,\epsilon_i \omega_i)$ converges in the Gromov--Hausdorff sense
to $B$, the metric on $B$ being induced from 
the (singular) Riemannian metric given, in local coordinates, by $W_0^{-1}
dy\otimes d\bar y$, with $W_0$ as defined in \S 4.

\noindent {\it Proof.} As usual, after choosing a topological 
zero-section of each $X_i$,
we can identify $X_i$ with $J$ as a manifold.  We may then view the
$\omega _i$ as corresponding to a sequence of Riemannian metrics 
$g_i$ on $J$, and prove 
that the sequence $J_i = (J,\epsilon_i g_i)$ converges in the Gromov--Hausdorff sense
to $B$ (with the given metric).

Using Remark 4.6, we can choose the class ${\bf B}_i$ determining
$\omega_i$ in a fundamental domain for the $B$-field by making, for
each $i$, a judicious choice of zero-section $\sigma_0$. 

Consider now $B$ along with the metric $W_0^{-1}dy\otimes d\bar y$.
Near each singular fibre one can find a coordinate $y$ so that
$\tau_1=1$ and $\tau_2={1\over 2\pi i}\log y+h(y)$, for some
holomorphic function $h$, and from this one can see that each
point of $\Delta\subseteq B$ is at finite distance under this
metric, and thus $B$ becomes a compact metric space using geodesic
distance.

Now we need to show that for each $\delta>0$, $d_{GH}(J_i ,B)<\delta$
for sufficiently large $i$. We will apply Definition 6.1 to the maps
$f_i=j:J \rightarrow B$ and $\sigma_0:B\rightarrow J $.

Choose, using Corollary 3.7 and Lemma 5.3,
for each point $p_j\in\Delta$, a small disc $D_j$ around $p_j$
with the property that
\item{(1)} $Diam(D_j)<\delta/100$.
\item{(2)} $Diam(f^{-1}(D_j))<\delta/100$ for sufficiently small
$\epsilon_i$. 

Let $U=B\setminus \bigcup D_j$. Now let $x_1,x_2\in J$. Let
$\gamma$ be a path joining $x_1$ and $x_2$ such that, for a given $i$,
$$l_{\epsilon_ig_i}(\gamma)<d_{\epsilon_i g_i}(x_1,x_2)+\delta/100.$$
Here $l$ denotes length, and the subscript denotes the metric being
used. At the risk of increasing the length of $\gamma$ by $24\delta/100$,
we can assume that $\gamma$ enters and leaves each $f^{-1}(D_j)$ at
most once, and write $\gamma=\gamma_1+\gamma_2$, with $\gamma_1\subseteq
f^{-1}(U)$ and $\gamma_2\subseteq f^{-1}(\bigcup D_j)$, with 
$l_{\epsilon_ig_i}(\gamma_2)\le 24\delta/100$. Now if $f^{-1}(U)$ carried
a semi-flat metric, then $f^{-1}(U)\rightarrow U$ would in fact
be a Riemannian submersion, and distances decrease under submersions.
On the other hand, if $\epsilon_i$ is sufficiently small, it follows
from Theorem 5.6 that the metric $\epsilon_ig_i$ is close to
a semi-flat metric in the $C^0$ sense. Thus for sufficiently
large $i$, depending on $\delta$,
$$l_B(f(\gamma_1))\le l_{SF}(\gamma_1)\le l_{\epsilon_ig_i}(\gamma_1)+
C(\epsilon_i),$$
where $l_{SF}$ denotes length with respect to the suitably normalized
semi-flat metric close to $\epsilon_i g_i$, and
$C(\epsilon_i)$ is a constant depending on  $\epsilon_i$ (and $\delta$)
but independent of the path. Furthermore, $C(\epsilon_i)\rightarrow 0$ as
$\epsilon_i\rightarrow 0$. Thus, possibly replacing $f(\gamma_2)$
with a shorter path, we see that
$$\eqalign{d_B(f(x_1),f(x_2))&\le l_B(f(\gamma_1))+24\delta/100\cr
&\le l_{\epsilon_i g_i}(\gamma_1)+C(\epsilon_i)+24\delta/100.\cr}$$
Thus for sufficiently small $\epsilon_i$, we always have
$$d_B(f(x_1),f(x_2))< d_{\epsilon_i g_i}(x_1,x_2)+\delta.$$

Next, let $y_1,y_2\in B$, and let $\gamma$ be a path joining $y_1$ and
$y_2$ with
$$l_B(\gamma)<d_B(y_1,y_2)+\delta/100.$$
As before, we can assume that $\gamma$ enters and leaves each $D_i$
once, and write $\gamma=\gamma_1+\gamma_2$. Consider now the metric
on $\sigma_0(B)$; locally, this takes the form
$\epsilon_i(W^{-1}+W|b|^2)dy\otimes d\bar y$ for some $W$ and $b$.
Again, the metric on $f^{-1}(U)$ is close to a semi-flat metric, hence this
metric is close, in the $C^0$ sense, to
$(W_0^{-1}+\epsilon_i^2W_0|b_{SF}|^2)dy\otimes d\bar y$. Now the point
of choosing ${\bf B}_i$ to be in a fundamental domain for the $B$-field
is that $|b_{SF}|^2$ can then be uniformly bounded, independent of $i$.
Thus for small $\epsilon_i$, this metric is close to the given metric
on $B$. Thus there exists a constant $C(\epsilon_i)$ with
$C(\epsilon_i)\rightarrow 0$ as $\epsilon_i\rightarrow 0$ such that
$$l_{\epsilon_i g_i}(\sigma_0(\gamma_1))\le l_B(\gamma_1)+C(\epsilon_i).$$
Therefore $d_{\epsilon_i g_i}(\sigma_0(y_1),\sigma_0(y_2))\le l_B(\gamma)
+C(\epsilon_i)+24\delta/100$
so for sufficiently small $\epsilon_i$,
$$d_{\epsilon_i g_i}(\sigma_0(y_1),\sigma_0(y_2))\le d_B(y_1,y_2)+\delta.$$
Thus for sufficiently small $\epsilon$,
$$|d_B(y_1,y_2)-d_{\epsilon_i g_i}(\sigma_0(y_1),\sigma_0(y_2))|<\delta$$
for all $y_1,y_2\in B$.

If $x_1,x_2\in J$, similar arguments show that
$$d_{\epsilon_i g_i}(x_1,x_2)<d_B(f(x_1),f(x_2))+\delta$$
by joining $x_1$ and $x_2$ by a path which first connects $x_1$
to $\sigma_0(B)$ inside a fibre or inside $f^{-1}(D_j)$ for some $j$,
then follows a geodesic inside $\sigma_0(B)$ to the fibre containing
$x_2$, and then connects up to $x_2$ inside this fibre.
The inequality follows for sufficiently small $\epsilon_i$ since
the diameter with respect to $\epsilon_i g_i$
of any fibre $f^{-1}(y)$ for $y\in U$, for small
$\epsilon_i$, is bounded by $C\epsilon_i$, where $C$ depends
only on the periods over $U$.

This shows
$$|d_{\epsilon_i g_i}(x_1,x_2)-d_B(f(x_1),f(x_2))|<\delta$$
for sufficiently small $\epsilon_i$. Finally, similar methods show
$$|d_{\epsilon_i g_i}(x_1,x_2)-d_{\epsilon_i g_i}
(\sigma_0(f(x_1)),\sigma_0(f(x_2))|
<\delta$$
for all $x_1,x_2\in X$, and $\epsilon_i$ sufficiently small. $\bullet$\bigskip

\noindent {\bf Remark 6.5.} The metric on the base $B$ is McLean's metric
(see [26], [19], [14]) on the base of the special Lagrangian
$T^2$-fibration obtained by hyperk\"ahler rotation. In higher
dimensions we also expect this metric to appear in the limit, showing
a residual effect of the conjectural special Lagrangian fibration.
This metric would then be singular along some subset of the limit,
corresponding to the limit of the discriminant loci of the conjectural
special Lagrangian fibrations. We hope this will be codimension 2.
See [16] for further speculation along these lines.

Conversely, we hope that one approach to understanding the existence
of special Lagrangian fibrations would be to prove Conjecture 6.2, which
gives us insight into the behaviour of Ricci-flat metrics near large complex
structure limits. However, it is clear that any approach to prove 
Conjecture 6.2 in higher dimensions must be substantially different
to the one given here for K3 surfaces, where we have made use of the
existence of special Lagrangian fibrations as well as the hyperk\"ahler
trick to reduce to a question of K\"ahler degenerations.
\bigskip

{\hd Bibliography.}

\item{[1]} Anderson, M., ``The $L^2$ Structure of Moduli Spaces of Einstein
Metrics on $4$-manifolds,'' {\it Geom. Funct. Anal.}, {\bf 2}, (1992)
29--89.
\item{[2]} Anderson, M., Kronheimer, P., and Lebrun, C., 
``Complete Ricci-flat K\"ahler Manifolds of Infinite Topological type,''
{\it Comm. Math. Phys.} {\bf 125}, (1989) 637--642.
\item{[3]} Abramowitz, M., and Stegun, I., eds., {\it
Handbook of Mathematical Functions,} Dover Publications, 1965.
\item{[4]} Aspinwall, P., and Morrison, D., ``String Theory on K3
surfaces,'' in {\it Essays on Mirror Manifolds II},
Greene, B.R., Yau, S.-T. (eds.) Hong Kong, International Press 1996, 703--716.
\item{[5]} Barth, W., Peters, C., and van de Ven, A., {\it
Compact Complex Surfaces}, Springer-Verlag, 1984.
\item{[6]} Besse, A., {\it Einstein Manifolds}, Springer-Verlag, 1987.
\item{[7]} Cheng, S.Y., and Yau, S.-T., ``On the existence of a complete 
K\"ahler metric on non compact complex manifolds and the regularity of 
Fefferman's equation.''
{\it Comm. Pure Appl. Math.} {\bf 33}, (1980) 507--544.
\item{[8]} Chern, S.-s., {\it Complex Manifolds without Potential
Theory, 2nd edition}, Springer-Verlag, 1979.
\item{[9]} Croke, C., ``Some Isoperimetric Inequalities and
Eigenvalue Estimates,'' {\it Ann. Sci. \'Ecole Norm. Sup. (4)},
{\bf 13}, (1980) 419--435.
\item{[10]} Dolgachev, I., ``Mirror Symmetry for Lattice Polarized K3
surfaces,'' Algebraic Geometry, {\bf 4}, {\it J. Math. Sci.} {\bf 81},
(1996) 2599--2630. 
\item{[11]} Gilbarg, D., and Trudinger, N., {\it Elliptic Partial
Differential Equations of Second Order,} 2. ed., rev. 3. printing, 
Springer-Verlag, 1998.
\item{[12]} Greene, B., Shapere, A., Vafa, C., and Yau, S.-T., 
``Stringy Cosmic Strings and Noncompact Calabi--Yau Manifolds,''
{\it Nucl. Phys. B}, {\bf 337}, (1990) 1--36.
\item{[13]} Gross, M., ``Special Lagrangian Fibrations I: Topology,'' 
in {\it Integrable Systems and Algebraic Geometry,} eds. M.-H. Saito,
Y. Shimizu and K. Ueno, World Scientific, 1998, 156--193.
\item{[14]} Gross, M., ``Special Lagrangian Fibrations II: Geometry,''
{\it Surveys in Differential Geometry,} Somerville: MA, International Press,
1999, 341--403.
\item{[15]} Gross, M., ``Topological Mirror Symmetry,'' {\it Invent. Math.,}
{\bf 144}, (2001), 75--137.
\item{[16]} Gross, M., ``Examples of Special Lagrangian Fibrations,''
preprint, math.AG/0012002, (2000).
\item{[17]} Gross, M., and Wilson, P.M.H., ``Mirror Symmetry via 3-tori for
a Class of Calabi--Yau Threefolds,'' 
{\it Math. Ann.}, {\bf 309}, (1997) 505--531.
\item{[18]} Harvey, R., and Lawson, H.B. Jr.,  ``Calibrated Geometries,'' {\it
 Acta
Math.} {\bf 148}, 47-157 (1982).
\item{[19]} Hitchin, N., ``The Moduli Space of Special Lagrangian 
Submanifolds,'' {\it Ann. Scuola Norm. Sup. Pisa Cl. Sci. (4)}, {\bf 25}
(1997) 503--515.
\item{[20]} Kobayashi, R., ``Moduli of Einstein Metrics on
a $K3$ Surface and Degeneration of type ${\rm I}$,'' {\it K\"ahler metric
and moduli spaces,} {\it Adv. Stud. Pure Math., 18-II}, Academic
Press, (1990) 257--311.
\item{[21]} Kobayashi, R., and Todorov, A., ``Polarized Period Map
for Generalized K3 Surfaces and the Moduli of Einstein Metrics,'' {\it
Tohoku Math. J. (2)}, {\bf 39}, (1987), 341--363.
\item{[22]}  Kontsevich, M., and Soibelman, Y., ``Homological Mirror Symmetry
and Torus Fibrations,'' preprint,  math.SG/0011041, (2000).
\item{[23]} Li, P., ``On the Sobolev Constant and the $p$-spectrum
of a Compact Riemannian Manifold,'' {\it Ann. Sci. \'Ecole Norm. Sup.
(4)}, {\bf 13}, (1980) 451--468.
\item{[24]} Lu, P., ``K\"ahler-Einstein Metrics on Kummer Threefold
and special Lagrangian Tori,'' {\it Comm. Anal. Geom.}, {\bf 7},
(1999) 787--806.
\item{[25]} Manin, Yu., ``Moduli, Motives, Mirrors,''
preprint, math.AG/0005144.
\item{[26]} McLean, R.C., `` Deformations of Calibrated Submanifolds,'' 
{\it Comm. Anal. Geom.} {\bf 6}, (1998) 705--747.
\item{[27]} Morrison, D., ``Compactifications of Moduli Spaces Inspired
by Mirror Symmetry,'' 
In {\it Journ\'ees de G\'eometrie Alg\'ebrique d'Orsay, Juillet 1992,} 
Asterisque
{\bf 218}, 243-271 (1993).
\item{[28]} Nergiz, S., and Sa\c clio\~ glu, C., ``A Quasiperiodic
Gibbons--Hawking Metric and Spacetime Foam,'' {\it Phys. Rev. D (3)},
{\bf 53}, (1996) 2240--2243.
\item{[29]} Ooguri, H., and Vafa, C., ``Summing up Dirichlet Instantons,''
{\it Phys. Rev. Lett.}, {\bf 77}, (1996) 3296--3298.
\item{[30]} Petersen, P., {\it Riemannian Geometry}, Springer-Verlag,
1997.
\item{[31]} Schoen, R., and Yau, S.-T., {\it Lectures on Differential
Geometry}, Conference Proceedings and Lecture Notes in Geometry
and Topology, Volume I, International Press, 1994.
\item{[32]} Strominger, A., Yau, S.-T., and Zaslow, E.,  ``Mirror Symmetry is
T-Duality,'' {\it Nucl. Phys.} {\bf B479}, (1996) 243--259.
\item{[33]} Tian, G., and Yau, S.-T., ``Complete K\"ahler manifolds with
zero Ricci curvature, I.''
{\it J. Amer. Math. Soc.} {\bf 3}, (1990) 579--609.
\item{[34]} Tian, G., and Yau, S.-T., ``Complete K\"ahler manifolds with
zero Ricci curvature, II.''
{\it Invent. math.} {\bf 106}, (1991) 27--60.
\item{[35]} Yau, S.-T., ``On the Ricci Curvature of a Compact K\"ahler
Manifold and the Complex Monge-Amp\`ere Equation. I.''
{\it Comm. Pure Appl. Math.} {\bf 31}, (1978) 339--411.
\item{[36]} Yau, S.-T., ``Survey on partial differential equations in differential 
geometry,'' In {\it Seminar on Differential Geometry,} 
{\it Annals of Mathematics Studies, 102}, Princeton University Press, (1982) 3--71.
\end